\documentclass[12pt,reqno]{amsart}
\usepackage{amsmath}
\usepackage{amssymb}
\usepackage{verbatim}
\usepackage{mathrsfs}
\usepackage{graphicx} 
\usepackage{caption}
\usepackage{subcaption}
\usepackage{hyperref}

\usepackage[top=0.5in,bottom=0.5in,left=0.7in,right=0.7in]{geometry}

\newtheorem{lemma}{Lemma}

\newtheorem{cor}[lemma]{Corollary}
\newtheorem{theorem}[lemma]{Theorem}

\newtheorem{exmp}{Example}

\newtheorem{definition}[lemma]{Definition}

\numberwithin{equation}{section}
\numberwithin{lemma}{section}

\newcommand{\C}{\mathbb{C}}    
\newcommand{\N}{\mathbb{N}}    
\newcommand{\PL}{\mathbb{P}}   
\newcommand{\R}{\mathbb{R}}    
\newcommand{\Z}{\mathbb{Z}}    

\newcommand{\wh}{\widehat}
\renewcommand{\le}{\leqslant}
\renewcommand{\ge}{\geqslant}
\newcommand{\bs}{\backslash}
\newcommand{\ol}{\overline}
\newcommand{\la}{\langle}
\newcommand{\ra}{\rangle}
\newcommand{\bo}{\mathscr{O}} 

\newcommand{\gl}{\lambda}
\newcommand{\eps}{\epsilon}
\newcommand{\pp}{\mathsf{p}}
\newcommand{\mrv}{{\mathring{v}}}
\newcommand{\mrw}{{\mathring{w}}}
\newcommand{\mra}{{\mathring{a}}}
\newcommand{\mrb}{{\mathring{b}}}
\newcommand{\mrphi}{{\mathring{\phi}}}
\newcommand{\ID}{{\text{Id}}}
\newcommand{\DG}{{\text{Diag}}}

\DeclareMathOperator{\FF}{\textsf{F}}

\DeclareMathOperator{\UU}{\textbf{\textsf{U}}}
\newcommand{\er}{\eqref}

\newcommand{\bpa}{\breve{a}}
\newcommand{\bpb}{\breve{b}}

\newcommand{\bpv}{\breve{v}}

\newcommand{\bpphi}{\breve{\phi}}

\newcommand{\bp}{ \begin{proof} }
\newcommand{\ep}{\end{proof} }
\newcommand{\be}{ \begin{equation} }
\newcommand{\ee}{ \end{equation} }
\newcommand{\tp}{\mathsf{T}}
\newcommand{\dm}{\mathsf{d}} 
\newcommand{\vgu}{\upsilon} 
\newcommand{\Vgu}{\Upsilon} 

\newcommand{\sr}{\operatorname{sr}}  
\newcommand{\vmo}{\operatorname{vm}}
\newcommand{\bvmo}{\operatorname{bvm}}

\newcommand{\lp}[1]{l_{#1}(\mathbb{Z})}

\newcommand{\lrs}[3]{(l_{#1}(\mathbb{Z}))^{#2\times #3}}
\newcommand{\sq}{l(\mathbb{Z})} 

\newcommand{\setsp}{\;:\;}     
\newcommand{\Lp}[1]{L_{#1}(\mathbb{R})}
\newcommand{\sd}{\mathcal{S}}  
\newcommand{\tz}{\mathcal{T}}  
\newcommand{\cM}{\mathcal{M}}

\newcommand{\ctV}{\tilde{\mathcal{V}}}
\newcommand{\cW}{\mathcal{W}}

\newcommand{\bpo}{\operatorname{bo}} 

\begin{document}

\title[Quasi-tight Multiframelets with High Balancing Orders]{Compactly Supported Quasi-tight Multiframelets with High Balancing Orders and Compact Framelet Transforms}

\author{Bin Han and Ran Lu}

\address{Department of Mathematical and Statistical Sciences,
University of Alberta, Edmonton,\quad Alberta, Canada T6G 2G1.
\quad {\tt bhan@ualberta.ca \quad rlu3@ualberta.ca}\quad {\tt http://www.ualberta.ca/$\sim$bhan}
}
\thanks{Research was supported in part by the Natural Sciences and Engineering Research Council of Canada (NSERC).
}


\makeatletter \@addtoreset{equation}{section} \makeatother

\begin{abstract}
Framelets (a.k.a. wavelet frames) are of interest in both theory and applications. Quite often, tight or dual framelets with high vanishing moments are constructed through the popular oblique extension principle (OEP). Though OEP can increase vanishing moments for improved sparsity,
it has a serious shortcoming for scalar framelets:
the associated discrete framelet transform is often not compact and deconvolution is
unavoidable. Here we say that a framelet transform is compact if it can be implemented by convolution using only finitely supported filters.
On one hand, \cite[Theorem~1.3]{han09} proves that for any scalar dual framelet constructed through OEP from any pair of scalar spline refinable functions, if it has a compact discrete framelet transform, then it can have at most one vanishing moment.
On the other hand, in sharp contrast to the extensively studied scalar framelets, multiframelets (a.k.a. vector framelets) derived through OEP from refinable vector functions are much less studied and are far from well understood.
Also, most constructed multiframelets often lack balancing property which reduces sparsity.  In this paper, we are particularly interested in quasi-tight multiframelets, which are special dual multiframelets but behave almost identically as tight multiframelets.
From any compactly supported \emph{refinable vector function having at least two entries}, we prove that we can always construct through OEP a compactly supported quasi-tight multiframelet such that (1) its associated discrete framelet transform is compact and has the highest possible balancing order; (2) all compactly supported framelet generators have the highest possible order of vanishing moments, matching the approximation/accuracy order of its underlying refinable vector function. This result demonstrates great advantages of OEP for multiframelets (retaining all the desired properties) over scalar framelets.
The key ingredient of our proof relies on
a newly developed normal form of matrix-valued filters, which is of independent interest and importance for greatly reducing the difficulty of studying refinable vector functions and multiframelets/multiwavelets.
We also study discrete multiframelet transforms employing OEP-based multiframelet filter banks.
To illustrate our theoretical results, we provide a few examples of quasi-tight multiframelets with a compact discrete framelet transform and having the highest possible balancing orders and vanishing moments.
\end{abstract}


\maketitle

\pagenumbering{arabic}

\section{Introduction and Main Results}

\subsection{Backgrounds on Framelets}
Wavelets and framelets are of interest in applications such as image processing and numerical algorithms (see e.g. \cite{chs02,daubook,dgm86,dhrs03,ds13,hanbook} and references therein). Framelets generalize wavelets by allowing redundancy with flexibility.
Let us first recall some necessary notations and definitions.
In this paper, we mainly deal with vector functions with entries in $\Lp{2}$.  By $f\in (\Lp{2})^{r\times s}$ we mean that $f$ is an $r\times s$ matrix of functions in $\Lp{2}$, and we define
\[
\la f, g\ra:=\int_\R f(x) \ol{g(x)}^\tp dx,\qquad
f\in (\Lp{2})^{r\times t}, g\in (\Lp{2})^{s\times t}.
\]
Note that $\la f,g\ra$ is an $r\times s$ matrix of complex numbers.
In particular, we define $(\Lp{2})^r:=(\Lp{2})^{r\times 1}$, the space of all $r\times 1$ column vector of functions in $\Lp{2}$. For $f\in (\Lp{2})^r$, we shall adopt the following notation:
\[
f_{\gl;k}(x):=|\gl|^{1/2} f(\gl x -k), \qquad x, \gl,k\in \R.
\]

Let $\dm$ be a positive integer with $\dm\ge 2$.
For vector functions $\eta\in (\Lp{2})^r$ and $\psi\in (\Lp{2})^s$, we say that $\{\eta;\psi\}$ is \emph{a $\dm$-framelet} in $\Lp{2}$ if
there exist positive constants $C_1$ and $C_2$ such that
\be\label{framelet}
C_1\|f\|_{\Lp{2}}^2\le
\sum_{k\in \Z} |\la f, \eta(\cdot-k)\ra|^2+\sum_{j=0}^\infty \sum_{k\in \Z}
|\la f, \psi_{\dm^j;k}\ra|^2\le C_2\|f\|_{\Lp{2}}^2, \qquad \forall\, f\in \Lp{2},
\ee
where $|\la f, \psi_{\dm^j;k}\ra|^2=\|\la f, \psi_{\dm^j;k}\ra\|^2_{l_2}:=\la f, \psi_{\dm^j;k}\ra \la \psi_{\dm^j;k},f\ra$.
If $r=1$ (i.e., $\eta$ is a function), then $\{\eta;\psi\}$ is called \emph{a scalar $\dm$-framelet}. If $r>1$ (i.e., $\eta$ is a vector function), then $\{\eta; \psi\}$ is often called \emph{a $\dm$-multiframelet}. For simplicity, in this paper we refer both of them as a framelet.
We now recall the definition of a dual $\dm$-framelet. Let $\eta,\tilde{\eta}\in (\Lp{2})^r$ and $\psi,\tilde{\psi}\in (\Lp{2})^s$. We say that $( \{\eta; \psi\},\{\tilde{\eta}; \tilde{\psi}\})$ is \emph{a dual $\dm$-framelet} in $\Lp{2}$ if both $\{\eta; \psi\}$ and $\{\tilde{\eta}; \tilde{\psi}\}$ are $\dm$-framelets in $\Lp{2}$, and for all $f,g\in \Lp{2}$,
\be\label{framelet:expr}
\la f,g\ra=\sum_{k\in \Z} \la f,\eta(\cdot-k) \ra \la \tilde{\eta}(\cdot-k),g\ra+\sum_{j=0}^\infty \sum_{k\in \Z}
\la f, \psi_{\dm^j;k}\ra \la \tilde{\psi}_{\dm^j;k},g\ra
\ee
with the series converging absolutely. It is straightforward to see that $( \{\eta; \psi\},\{\tilde{\eta}; \tilde{\psi}\})$ is a dual $\dm$-framelet if and only if $(\{\tilde{\eta}; \tilde{\psi}\},\{\eta; \psi\})$ is a dual $\dm$-framelet in $\Lp{2}$. $\{\eta;\psi\}$ is called \emph{a tight $\dm$-framelet} in $\Lp{2}$ if $(\{\eta;\psi\}, \{\eta;\psi\})$ is a dual $\dm$-framelet.
Framelets are often constructed from refinable vector functions. By $\lp{0}$ we denote the set of all finitely supported sequences on $\Z$. For $a=\{a(k)\}_{k\in \Z}\in \lrs{0}{r}{s}$, we define
$$
\wh{a}(\xi):=\sum_{k\in \Z} a(k) e^{-ik\xi},\qquad \xi\in \R,
$$
which is an $r\times s$ matrix of $2\pi$-periodic trigonometric polynomials. A vector function $\phi\in (\Lp{2})^r$ is \emph{a $\dm$-refinable vector function} with a (matrix-valued) refinement filter/mask $a\in \lrs{0}{r}{r}$ if
\[
\wh{\phi}(\dm \xi)=\wh{a}(\xi)\wh{\phi}(\xi),\qquad \forall \xi\in \R.
\]
Here, the Fourier transform of $f\in \Lp{1}$ is defined to be $\wh{f}(\xi):=\int_\R f(x) e^{-ix\xi} dx$ for $\xi\in \R$, and the definition is naturally extended to $\Lp{2}$ functions and tempered distributions. $\wh{\phi}$ is the $r\times 1$ vector function obtained by taking entry-wise Fourier transform on $\phi$. If $r=1$, then $\phi$ is called a (scalar) refinable function.
One of the most important examples of refinable scalar functions are B-splines. For $m\in \N$, the \emph{B-spline function $B_m$ of order $m$} is defined by
\be \label{bspline}
B_1:=\chi_{[0,1]} \quad \mbox{and}\quad B_m:=B_{m-1}*B_1=\int_0^1 B_{m-1}(\cdot-t) dt.
\ee
The B-spline function $B_m$ is a piecewise polynomial function, belongs to $C^{m-1}(\R)$ with support $[0,m]$, and is $\dm$-refinable: $\wh{B_m}(\dm \xi)=\wh{a^B_{m,\dm}}(\xi)\wh{B_m}(\xi)$ with %
\be \label{aBm}
\wh{a^B_{m,\dm}}(\xi):=\dm^{-m} (1+e^{-i\xi}+\cdots+e^{-i(\dm-1)\xi})^m.
\ee

The most general way of deriving framelets from refinable (vector) functions is through the oblique extension principle (OEP), which was introduced in \cite{dhrs03} and \cite{chs02} for scalar framelets. Here let us recall a special version of the oblique extension principle stated in \cite[Theorem~6.4.1]{hanbook} (c.f. \cite[Theorem~1.1]{han09}) for constructing dual framelets from refinable vector functions.

\begin{theorem}\label{thm:df}
Let $\dm$ be a positive integer with $\dm\ge 2$.
Let $\theta,\tilde{\theta},a,\tilde{a}\in \lrs{0}{r}{r}$ and $\phi,\tilde{\phi}\in (\Lp{2})^r$ be compactly supported $\dm$-refinable vector functions satisfying $\wh{\phi}(\dm\xi)=\wh{a}(\xi)\wh{\phi}(\xi)$
and
$\wh{\tilde{\phi}}(\dm\xi)=\wh{\tilde{a}}(\xi)\wh{\tilde{\phi}}(\xi)$.
For finitely supported matrix-valued filters $b,\tilde{b}\in \lrs{0}{s}{r}$,
define
\[
\wh{\eta}(\xi):=\wh{\theta}(\xi)\wh{\phi}(\xi),\quad
\wh{\psi}(\xi):=\wh{b}(\xi/\dm)\wh{\phi}(\xi/\dm)
\quad \mbox{and}\quad
\wh{\tilde{\eta}}(\xi):=\wh{\tilde{\theta}}(\xi)\wh{\tilde{\phi}}(\xi),\quad
\wh{\tilde{\psi}}(\xi):=\wh{\tilde{b}}(\xi/\dm)\wh{\tilde{\phi}}(\xi/\dm).
\]
Then $(\{\eta;\psi\}, \{\tilde{\eta}; \tilde{\psi}\})$ is a dual $\dm$-framelet in $\Lp{2}$ if
\begin{enumerate}
\item[(1)] $\ol{\wh{\phi}(0)}^\tp \wh{\Theta}(0)\wh{\tilde{\phi}}(0)=1$ with $\wh{\Theta}(\xi):=\ol{\wh{\theta}(\xi)}^\tp \wh{\tilde{\theta}}(\xi)$;
\item[(2)] all entries in $\psi$ and $\tilde{\psi}$ have at least one vanishing moment, i.e., $\wh{\psi}(0)=\wh{\tilde{\psi}}(0)=0$.
\item[(3)] $(\{a;b\},\{\tilde{a};\tilde{b}\})_{\Theta}$ forms an OEP-based dual $\dm$-framelet filter bank satisfying
\be \label{dffb}
\ol{\wh{{a}}(\xi)}^\tp \wh{\Theta}(\dm \xi){\wh{\tilde{a}}(\xi)}+\ol{\wh{{b}}(\xi)}^\tp {\wh{\tilde{b}}(\xi)}=\wh{\Theta}(\xi) \quad \mbox{and}\quad
\ol{\wh{{a}}(\xi)}^\tp \wh{\Theta}(\dm \xi){\wh{\tilde{a}}(\xi+\tfrac{2\pi \gamma}{\dm})}+\ol{\wh{{b}}(\xi)}^\tp {\wh{\tilde{b}}(\xi+\tfrac{2\pi \gamma}{\dm})}=0
\ee
for all $\gamma=1,\ldots, \dm-1$ and $\xi\in \R$.
\end{enumerate}
\end{theorem}
Under the additional assumption that $\mbox{span}\{\wh{\phi}(\xi+2\pi k) \setsp k\in \Z\}=\C^r=\mbox{span}(\{\wh{\tilde{\phi}}(\xi+2\pi k) \setsp k\in \Z\}$ for infinitely many $\xi\in [-\pi,\pi]$, \cite[Theorem~6.4.1]{hanbook} also shows that the above items (1)--(3) are necessary conditions for $(\{\eta;\psi\}, \{\tilde{\eta}; \tilde{\psi}\})$ to be a dual $\dm$-framelet in $\Lp{2}$. For two smooth functions $f$ and $g$, for convenience we shall adopt the following notation: For $m\in \N$ and $\xi_0\in\R$,
\[
f(\xi)=g(\xi)+\bo(|\xi-\xi_0|^m),\quad \xi \to \xi_0 \quad \mbox{stands for}\quad f^{(j)}(\xi_0)=g^{(j)}(\xi_0), \quad\forall\, j=0,\ldots,m-1.
\]
The sparsity of the framelet representation in \eqref{framelet:expr} is largely due to the vanishing moments of the framelet generators $\psi$ and $\tilde{\psi}$.
A compactly supported vector function $\psi\in (\Lp{2})^s$ has \emph{$n$ vanishing moments} if $\wh{\psi}(\xi)=\bo(|\xi|^n)$ as $\xi \to 0$. We further define $\vmo(\psi):=n$ with $n$ being the largest such integer.
The vanishing moments of $\psi$ in Theorem~\ref{thm:df} is closely related to the approximation property of the refinable vector function $\phi$ and sum rules of its associated refinement filter/mask $a$.
For a finitely supported matrix-valued filter $a\in \lrs{0}{r}{r}$, we say that $a$ has \emph{$m$ sum rules with respect to the dilation factor $\dm$} with a matching filter $\vgu\in \lrs{0}{1}{r}$ if $\wh{\vgu}(0)\ne 0$ and
\be \label{sr}
\wh{\vgu}(\dm\xi)\wh{a}(\xi)=\wh{\vgu}(\xi)+
\bo(|\xi|^m)
\quad \mbox{and}
\quad
\wh{\vgu}(\dm\xi)
\wh{a}(\xi+\tfrac{2\pi \gamma}{\dm})=\bo(|\xi|^m),\quad \xi\to 0,\quad\gamma=1,\ldots, \dm-1.
\ee
In particular, we define $\sr(a,\dm):=m$ with $m$ being the largest possible integer in \eqref{sr}.
In the scalar case (i.e., $r=1$), \eqref{sr} is equivalent to $(1+e^{-i\xi}+\cdots+e^{-i(\dm-1)\xi})^m \mid \wh{a}(\xi)$ by taking $\wh{\vgu}(\xi)=1/\wh{\phi}(\xi)+\bo(|\xi|^m)$ as $\xi \to0$ and noting $\wh{\phi}(\xi):=\prod_{j=1}^\infty \wh{a}(\dm^{-j}\xi)$ with $\wh{a}(0)=1$.

For a dual $\dm$-framelet constructed in Theorem~\ref{thm:df},
regardless of the choice of filters $\theta$ and $\tilde{\theta}$, it is well known through a simple argument that $\vmo(\psi)\le \sr(\tilde{a},\dm)$ and $\vmo(\tilde{\psi})\le \sr(a,\dm)$. From a pair of B-spline filters $a^B_{m,\dm}$ and $a^B_{n,\dm}$ in \eqref{aBm}, any dual $\dm$-framelet, derived through Theorem~\ref{thm:df} with the trivial choice $\wh{\theta}(\xi)=\wh{\tilde{\theta}}(\xi)=1$,
has at most one vanishing moment, i.e., $\vmo(\psi)=\vmo(\tilde{\psi})=1$ (see \cite{dhrs03,rs97}), even though $\sr(a^B_{m,\dm},\dm)=m$ can be arbitrarily large.
As well explained in \cite{chs02,dhrs03} (also see \cite{dh04,han09,hm03}), the main advantage of OEP is to increase the orders of vanishing moments of $\psi$ and $\tilde{\psi}$ in Theorem~\ref{thm:df} by properly choosing the filters $\theta$ and $\tilde{\theta}$ so that $\vmo(\psi)= \sr(\tilde{a},\dm)$ and $\vmo(\tilde{\psi})= \sr(a,\dm)$. A lot of compactly supported scalar tight or dual framelets with the highest possible vanishing moments have been constructed in the literature, to mention only a few, see \cite{cs10,ch00,chs02,dh04,dhrs03,dhacha,dh18pp,han97,han09,hanbook,hm03,hm05,jqt01,js15,mothesis,rs97,sel01} and many references therein.
In particular, see Chapter~3 of \cite{hanbook} for comprehensive study and references on scalar tight or dual framelets.

Multiwavelets have certain advantages over scalar wavelets and have been initially studied in \cite{ghm94,glt93} and references therein.
In sharp contrast to the extensively studied OEP-based scalar framelets,
constructing multiframelets (a.k.a. vector framelets) through OEP in Theorem~\ref{thm:df} is much more difficult than scalar framelets. To our best knowledge, we are only aware of \cite[Chapter~2]{mothesis} for studying OEP-based tight multiframelets, and \cite{han09,hm03} for investigating OEP-based dual multiframelets with vanishing moments.
However, as observed in \cite{han09}, OEP-based dual framelets with OEP-based dual framelet filter banks in Theorem~\ref{thm:df} has a drawback: their associated discrete framelet transforms are not compact. Here we say that a discrete framelet transform is compact if it can be implemented by convolution using only finitely supported filters. Discrete framelet transforms have been mentioned in \cite{dhrs03} for scalar framelets and in \cite{han09} for multiframelets. For multiframelets with high vanishing moments and with multiplicity $r>1$, it is well known that they often have a much lower balancing order (e.g., see \cite{cj00,han09,han10,LV98,lv98,sel00}), leading to a significant loss of sparsity of their associated discrete multiframelet transform.
Since we are not aware of any detailed study on discrete multiframelet transforms using OEP-based filter banks in the literature,
we shall study in Section~\ref{sec:ffrt} discrete multiframelet transforms using OEP-based matrix-valued filter banks as well as their discrete vanishing moments and balancing property. Then we shall explain in detail in Section~\ref{sec:ffrt} why almost all constructed OEP-based dual framelets with OEP-based filter banks in the literature have non-compact discrete multiframelet transforms and cannot avoid nonstable deconvolutions. Such drawbacks of OEP-based scalar framelets and multiframelets make them much less attractive in applications.

\subsection{Our Contributions}

To avoid all the above-mentioned shortcomings of OEP for multiframelets (see Section~\ref{sec:ffrt} for details), in this paper we are particularly interested in quasi-tight multiframelets, which are special dual multiframelets, but behave almost identically as tight multiframelets.
Let $\eta\in (\Lp{2})^r$ and $\psi=[\psi^1,\ldots,\psi^s]^\tp \in (\Lp{2})^s$.
For $\eps_1,\ldots,\eps_s\in \{\pm 1\}$, we say that $\{\eta; \psi\}_{(\eps_1,\ldots,\eps_s)}$ is
\emph{a quasi-tight $\dm$-framelet} in $\Lp{2}$ if $(\{\eta; [\psi^1,\ldots,\psi^s]^\tp\},\{\eta; [\eps_1\psi^1,\dots,\eps_s\psi^s]^\tp\})$ is a dual $\dm$-framelet in $\Lp{2}$. Or equivalently, we say that
$\{\eta; \psi\}_{(\eps_1,\ldots,\eps_s)}$
is a quasi-tight $\dm$-framelet if $\{\eta; \psi\}$ is a $\dm$-framelet satisfying \eqref{framelet} and
\be \label{qtf:expr}
f=\sum_{k\in \Z} \la f, \eta(\cdot-k)\ra \eta(\cdot-k)+\sum_{j=0}^\infty
\sum_{\ell=1}^s \sum_{k\in \Z}
\eps_\ell \la f, \psi^\ell_{\dm^j;k}\ra \psi^\ell_{\dm^j;k},
\qquad \forall\,f\in \Lp{2}
\ee
with the above series converging unconditionally in $\Lp{2}$, where $\psi=[\psi^1,\ldots,\psi^s]^\tp$.
Obviously, if $\eps_1=\cdots=\eps_s=1$, then a quasi-tight $\dm$-framelet
is simply a tight $\dm$-framelet. One example of quasi-tight framelets appeared in \cite[Example~3.2.2]{hanbook}. Scalar quasi-tight framelets have been studied in \cite{dh18pp} for dimension one and in \cite{dhacha} for high dimensions with the added feature of directionality. Furthermore, if $\{\eta; \psi\}_{\eps_1,\ldots,\eps_s}$ with $\psi=[\psi^1,\ldots,\psi^s]^\tp$ is a quasi-tight $\dm$-framelet in $\Lp{2}$, by \cite[Proposition~5]{han12}, then $\{\psi\}_{(\eps_1,\ldots,\eps_s)}$ is \emph{a homogeneous quasi-tight $\dm$-framelet} in $\Lp{2}$ satisfying
\eqref{framelet} and
\[
f=\sum_{j\in \Z}
\sum_{\ell=1}^s \sum_{k\in \Z}
\eps_\ell \la f, \psi^\ell_{\dm^j;k}\ra \psi^\ell_{\dm^j;k},
\qquad \forall\,f\in \Lp{2}
\]
with the above series converging unconditionally in $\Lp{2}$.

We are now ready to state the main result of this paper on quasi-tight multiframelets having all the desired properties.

\begin{theorem}\label{thm:qtf}
	Let $\dm\ge 2$ be an integer and $\phi\in (\Lp{2})^r$ be
a compactly supported $\dm$-refinable vector function with a matrix-valued refinement filter/mask $a\in \lrs{0}{r}{r}$ satisfying $\wh{\phi}(\dm\xi)\wh{a}(\xi)\wh{\phi}(\xi)$.
	Suppose that the filter $a$ has $m$ sum rules with respect to the dilation factor $\dm$ satisfying \eqref{sr} with a matching filter $\vgu\in \lrs{0}{1}{r}$ such that $\wh{\vgu}(0)\wh{\phi}(0)=1$.
	If the multiplicity $r\ge 2$, then there exist filters $\theta\in \lrs{0}{r}{r}$, $b\in \lrs{0}{s}{r}$ and $\eps_1,\ldots,\eps_s\in \{\pm 1\}$ such that
	\begin{enumerate}
		\item[(1)] $\{\mrphi; \psi\}_{(\eps_1,\ldots,\eps_s)}$ is a compactly supported quasi-tight $\dm$-framelet in $\Lp{2}$, where $\wh{\mathring{\phi}}(\xi):=\wh{\theta}(\xi)\wh{\phi}(\xi)$ and $\wh{\psi}(\xi):=\wh{b}(\xi/\dm)\wh{\phi}(\xi/\dm)$. Furthermore, $\psi$ has $m$ vanishing moments.
		
		\item[(2)] $\wh{\theta}$ (or for simplicity just $\theta$) is \textbf{strongly invertible}, i.e., $\wh{\theta}^{-1}$ is also an $r\times r$ matrix of $2\pi$-periodic trigonometric polynomials. The filter bank $\{\mathring{a}; \mathring{b}\}_{(\eps_1,\ldots,\eps_s)}$ is a finitely supported quasi-tight $\dm$-framelet filter bank, i.e.,
		\begin{align}
		 &\ol{\wh{\mathring{a}}(\xi)}^\tp {\wh{\mathring{a}}(\xi)}+\ol{\wh{\mathring{b}}(\xi)}^\tp \DG(\eps_1,\ldots,\eps_s) {\wh{\mathring{b}}(\xi)}=I_r, \label{qtffb=1}\\
		 &\ol{\wh{\mathring{a}}(\xi)}^\tp {\wh{\mathring{a}}(\xi+\tfrac{2\pi \gamma}{\dm})}+\ol{\wh{\mathring{b}}(\xi)}^\tp \DG(\eps_1,\ldots,\eps_s){\wh{\mathring{b}}(\xi+\tfrac{2\pi \gamma}{\dm})}=0,\label{qtffb=0}
		\end{align}
		for all $\gamma=1,\ldots,\dm-1$ and for all $\xi \in \R$, where the finitely supported matrix-valued filters $\mathring{a}\in \lrs{0}{r}{r}$ and $\mathring{b}\in \lrs{0}{s}{r}$ are defined by
		\be \label{mab}
		 \wh{\mathring{a}}(\xi):=\wh{\theta}(\dm \xi) \wh{a}(\xi) \wh{\theta}(\xi)^{-1}
		\quad \mbox{and}\quad
		 \wh{\mathring{b}}(\xi):=\wh{b}(\xi) \wh{\theta}(\xi)^{-1}.
		\ee
		\item[(3)] The filter $\mathring{b}$ has $m$ balanced vanishing moments (see Section~\ref{sec:ffrt} for its definition).
		\item[(4)] The associated discrete multiframelet transform employing the quasi-tight $\dm$-framelet filter bank
		$\{\mathring{a}; \mathring{b}\}_{(\eps_1,\ldots,\eps_s)}$ is compact and has the highest possible balancing order, i.e., $\bpo(\{\mra;\mrb\},\dm)=m$.(see Section~\ref{sec:ffrt} for details).
	\end{enumerate}
	Moreover, the compactly supported vector functions $\mathring{\phi}$ and $\psi$ satisfy
	\be \label{mphi:mpsi}
	 \wh{\mathring{\phi}}(\dm\xi)=\wh{\mathring{a}}(\xi)\wh{\mathring{\phi}}(\xi)
	\quad \mbox{and}\quad
	 \wh{\psi}(\dm\xi)=\wh{\mathring{b}}(\xi)\wh{\mathring{\phi}}(\xi).
	\ee
\end{theorem}

The key ingredient to prove Theorem~\ref{thm:qtf} is a newly developed normal form of a matrix-valued filter. The main idea of the normal form theory is to transform the original filter $a$ to a new filter $\mra$ with desired features so that we can implement construction techniques as in the scalar case. Generalizing \cite[Theorem~2.1]{han09} but under much weaker conditions, we obtain Theorem~\ref{thm:normalform} below which is a general result on a normal form of a matrix-valued filter. As a special case of Theorem~\ref{thm:normalform}, we can achieve a very nice structure on the newly transformed filter (Theorem~\ref{thm:normalform:2} below), which plays a key role in our study of quasi-tight framelets in this paper.

We say that a function $f$ is \emph{smooth near the origin} if all the derivatives of $f$ at the origin exist.

\begin{theorem}\label{thm:normalform}
	Let $\dm\ge 2$ be an integer and $a\in \lrs{0}{r}{r}$ be a finitely supported matrix-valued filter.
	Let $\phi$ be an $r\times 1$ vector of
	compactly supported distributions satisfying $\wh{\phi}(\dm \xi)=\wh{a}(\xi)\wh{\phi}(\xi)$ with $\wh{\phi}(0)\ne 0$.
	Suppose that the filter $a$ has $m$ sum rules with respect to $\dm$ satisfying \eqref{sr} with a matching filter $\vgu\in \lrs{0}{1}{r}$ such that $\wh{\vgu}(0)\wh{\phi}(0)=1$.
	Let $\wh{\mathring{\vgu}}$ be a $1\times r$ row vector and $\wh{u_\phi}$ be an $r\times 1$ column vector
	such that all the entries of
	$\wh{\mathring{\vgu}}$ and $\wh{u_\phi}$ are functions which are smooth near the origin and
	\be \label{vguphi=1:new}
	 \wh{\mathring{\vgu}}(\xi)\wh{u_\phi}(\xi)=1+\bo(|\xi|^m), \qquad \xi \to 0.
	\ee
	If the multiplicity $r\ge 2$, then for any positive integer $n\in \N$, there exists a strongly invertible $r\times r$ matrix $\wh{U}$ of $2\pi$-periodic trigonometric polynomials such that
	\be \label{normalform:general}
	\wh{\vgu}(\xi) \wh{U}(\xi)^{-1}=\wh{\mathring{\vgu}}(\xi)+\bo(|\xi|^m)
	\quad \mbox{and}\quad
	\wh{U}(\xi) \wh{\phi}(\xi)=
	 \wh{u_\phi}(\xi)+\bo(|\xi|^n),\qquad \xi\to 0.
	\ee
	Define
	 $\wh{\mathring{\phi}}(\xi):=\wh{U}(\xi) \wh{\phi}(\xi)$
	and $\wh{\mathring{a}}(\xi):=\wh{U}(\dm\xi)\wh{a}(\xi) \wh{U}(\xi)^{-1}$.
	Then the following statements hold:
	\begin{enumerate}
		\item[(i)] The new vector function $\mathring{\phi}$ is a vector of compactly supported distributions satisfying the refinement equation
		$\wh{\mathring{\phi}}(\dm \xi)=\wh{\mathring{a}}(\xi)\wh{\mathring{\phi}}(\xi)$ for all $\xi\in \R$
		and $\wh{\mathring{\phi}}(\xi)=\wh{u_\phi}(\xi)+\bo(|\xi|^n)$ as $\xi \to 0$.
		\item[(ii)] The new finitely supported matrix filter/mask $\mathring{a}$ has $m$ sum rules with respect to $\dm$ with the matching filter $\mathring{\vgu}$
		such that $\wh{\mathring{\vgu}}(0)\wh{\mathring{\phi}}(0)=1$, i.e., \eqref{sr} holds with $a$ and $\vgu$ being replaced by $\mathring{a}$ and $\mathring{\vgu}$, respectively.
	\end{enumerate}
\end{theorem}

We shall show in Lemma~\ref{lem:vguphi} that the condition in \eqref{vguphi=1:new} of Theorem~\ref{thm:normalform} is a necessary condition. As a special case of Theorem~\ref{thm:normalform}, we have the following result, which is used in the proof of Theorem~\ref{thm:qtf}.

\begin{theorem}\label{thm:normalform:2}
	Let $\dm\ge 2$ be a positive integer and $a\in \lrs{0}{r}{r}$ be a finitely supported matrix-valued filter.
	Let $\phi$ be an $r\times 1$ vector of compactly supported distributions satisfying
	$\wh{\phi}(\dm \xi)=\wh{a}(\xi)\wh{\phi}(\xi)$ with $\wh{\phi}(0)\ne 0$.
	Suppose that the filter $a$ has $m$ sum rules with respect to $\dm$ satisfying \eqref{sr} with a matching filter $\vgu\in \lrs{0}{1}{r}$ and $\wh{\vgu}(0)\wh{\phi}(0)=1$.
	If the multiplicity $r\ge 2$, then for any positive integer $n\in \N$, there exists a strongly invertible $r\times r$ matrix $\wh{U}$ of $2\pi$-periodic trigonometric polynomials such that the following properties hold:
	
	\begin{enumerate}
		\item[(i)] $\wh{\mathring{a}}(\xi):=\wh{U}(\dm\xi) \wh{a}(\xi)\wh{U}(\xi)^{-1}$ takes the form
		\be \label{normalform}
		\left[
		\begin{matrix}
			 (1+e^{-i\xi}+\cdots+e^{-i(\dm-1)\xi})^m P_{1,1}(\xi) &(1-e^{-i\dm \xi})^m P_{1,2}(\xi)\\
			(1-e^{-i\xi})^n P_{2,1}(\xi) &P_{2,2}(\xi)\end{matrix}\right],
		\ee
		with
		\be \label{normalform:11}
		 \wh{\mathring{a}_{1,1}}(\xi):=(1+e^{-i\xi}+\cdots+e^{-i(\dm-1)\xi})^m P_{1,1}(\xi)=1+\bo(|\xi|^n),\quad \xi\to 0,
		\ee
		where $P_{1,1}, P_{1,2}, P_{2,1}$ and $P_{2,2}$ are some $1\times 1$, $1\times (r-1)$, $(r-1)\times 1$ and $(r-1)\times (r-1)$ matrices of $2\pi$-periodic trigonometric polynomials. Moreover, define
		 \be\label{mrv:1}\wh{\mathring{\vgu}}(\xi):=[\wh{\mathring{\vgu}_1}(\xi),\ldots,
		 \wh{\mathring{\vgu}_r}(\xi)]:=\wh{\vgu}(\xi)\wh{U}(\xi)^{-1},\ee
		 \be\label{mrphi:1}\wh{\mathring{\phi}}(\xi):=[\wh{\mathring{\phi}_1}(\xi),\ldots,
		 \wh{\mathring{\phi}_r}(\xi)]^\tp:=\wh{U}(\xi)\wh{\phi}(\xi),\ee
		we have $\wh{\mathring{\phi}}(\dm \xi)=\wh{\mathring{a}}(\xi)
		 \wh{\mathring{\phi}}(\xi)$ with
		\be \label{normalform:phi} \wh{\mathring{\phi}_1}(\xi)=1+\bo(|\xi|^n)
		\quad \mbox{and}\quad
		 \wh{\mathring{\phi}_\ell}(\xi)=\bo(|\xi|^n),\quad \xi\to 0, \ell=2,\ldots,r,
		\ee
		and $\mra$ has $m$ sum rules with respect to $\dm$ with the matching filter $\mathring{\vgu}$ satisfying
\be\label{normalform:vgu}\wh{\mathring{\vgu}_1}(\xi)=1+\bo(|\xi|^m)
		\quad \mbox{and}\quad
		 \wh{\mathring{\vgu}_\ell}(\xi)=\bo(|\xi|^m),\quad \xi\to 0, \ell=2,\ldots,r.\ee
		
		\item[(ii)] If in addition
		 \be\label{moment:special}\wh{\vgu}(\xi)=\frac{\ol{\wh{\phi}(\xi)}^{\tp}}{\|\wh{\phi}(\xi)\|^2}+\bo(|\xi|^m),\quad\xi\to 0,\ee
		then $\wh{U}$ in item (i) can be chosen such that the following ``almost orthogonal" structure holds:
		 \be\label{eq:ortho}\ol{\wh{U}(\xi)}^{-\tp}\wh{U}(\xi)^{-1}=\DG\left(\|\wh{\phi}(\xi)\|^2,\|\wh{u_2}(\xi)\|^2,\dots, \|\wh{u_r}(\xi)\|^2\right)+\bo(|\xi|^{\max(m,n)}),\quad\xi\to 0,\ee
		where $\wh{u_j}$ is the $j$-th column of $\wh{U}^{-1}$ for $j=2,\dots,r$.
	\end{enumerate}
Conversely, if there exists $\wh{U}$ such that item (i) and \er{eq:ortho} hold, then \er{moment:special} must hold.
\end{theorem}

We comment on some important features involved in our contributions.

\begin{enumerate}
	\item[(1)] Our main result Theorem~\ref{thm:qtf} demonstrates that we can construct quasi-tight multiframelets from any refinable vector functions. This is not like existing works in the literature that study tight framelets, which often require that the refinable vector function $\phi$ should have stable integer shifts. This condition guarantees the existence of $\Theta\in\lrs{0}{r}{r}$ (which is often not strongly invertible at all) such that $\cM_{a,\Theta}$ is positive semi-definite, where $\cM_{a,\Theta}$ is defined in \er{cond:oep:tf} (see \cite[Proposition 3.4 and Theorem 4.3]{mo06}). The positive semi-definiteness of $\cM_{a,\Theta}$ is a necessary condition for the existence of tight framelets.

	\item[(2)] Theorem~\ref{thm:qtf} demonstrates great advantages of OEP for multiframelets. In the scalar case ($r=1$), OEP can increase the order of vanishing moments on framelet generators, but quite often it is inevitable to sacrifice the compactness of the associated discrete framelet transform. For example, \cite[Theorem~1.3]{han09} proves that for any scalar dual framelet constructed through OEP from any pair of scalar spline refinable functions, if it has a compact framelet transform, then it can have at most one vanishing moment. Besides, most of the multiframelets constructed in existing literatures lack the balancing property, which reduces sparsity when implementing a multi-level discrete multiframelet transform. Theorem~\ref{thm:qtf} guarantees the existence of quasi-tight multiframelets with all desired properties: (i) high order balanced vanishing moments on framelet generators; (ii) a compact and balanced associated discrete multiframelet transform.

	\item[(3)] The normal form of a matrix-valued filter greatly facilitates the study and construction of multiframelets and multiwavelets. It allows us to study multiframelets and multiwavelets in almost the same way as what we do in the scalar case. The study of the normal form of a matrix-valued filter is of interest in its own right.
	
\end{enumerate}

\subsection{Paper Structure}

The structure of the paper is as follows. We shall prove Theorems~\ref{thm:normalform} and \ref{thm:normalform:2} on the normal form of a matrix-valued filter in
Section~\ref{sec:normalform}. We will demonstrate that the normal form theory makes the study of refinable vector functions and matrix-valued filters almost as easy as the scalar case (i.e., $r=1$). In other words, using the normal form of a matrix-valued filter allows one to adopt almost all techniques from the scalar case to study multiwavelets and refinable vector functions.
In Section~\ref{sec:ffrt}, we shall study the discrete multiframelet transform using an OEP-based dual multiframelet filter bank. Then we shall discuss its various properties including the balancing property and the notion of balanced vanishing moments.
We shall also explain in Section~\ref{sec:ffrt} the possible shortcomings for using OEP-based filter banks and how to overcome such shortcomings.
In Section~\ref{sec:qtf}, we shall prove our main result stated in Theorem~\ref{thm:qtf}. In Section~\ref{sec:momfilter}, we characterize all possible strongly invertible filters $\theta$ in Theorem~\ref{thm:qtf} so that a quasi-tight framelet and a quasi-tight framelet filter bank can be constructed and satisfy all the desired properties in items (1)--(4) of Theorem~\ref{thm:qtf}. Finally, in Section~\ref{sec:example}, we shall provide a few examples of compactly supported quasi-tight multiframelets with all the desired properties.

\section{Normal Form of a Matrix-valued Filter/Mask}\label{sec:normalform}

In this section we shall prove Theorem~\ref{thm:normalform} on the existence of a normal form of a matrix filter/mask, and then use this to further obtain an ideal normal form as stated in Theorem~\ref{thm:normalform:2}.

We first make some comments on the importance of the normal form of a matrix-valued filter in the study of refinable vector functions and multiwavelets/multiframelets. The normal form (also called the canonical form) of a matrix-valued filter was initially introduced in \cite[Theorem~2.2]{hm03} for dimension one and was further developed in \cite[Proposition~2.4]{han03} for high dimensions to study multivariate vector subdivision schemes and multivariate refinable vector functions.

In the scalar case (i.e., $r=1$), recall that a scalar filter $a$ has $m$ sum rules if and only if $(1+e^{-i\xi}+\cdots+e^{-i(\dm-1)\xi})^m \mid \wh{a}(\xi)$. That is, $\wh{a}(\xi)=(1-e^{-i\dm\xi})^m A(\xi) (1-e^{-i\xi})^{-m}=
(1+e^{-i\xi}+\cdots+e^{-i(\dm-1)\xi})^m A(\xi)$ for a unique $2\pi$-periodic trigonometric polynomial $A(\xi)$. Now consider the case $r>1$. If a filter $\mra$ takes the form \eqref{normalform} in item (i) of Theorem~\ref{thm:normalform:2}, we can factorize $\mra$ as
$\wh{\mra}(\xi)=B(\dm \xi) A(\xi) B(\xi)^{-1}$
with
\[
B(\xi):=\wh{U}(\xi)^{-1} \left[ \begin{matrix} (1-e^{-i\xi})^m &\\
&I_{r-1}\end{matrix}\right],\quad
A(\xi):=\left[ \begin{matrix} P_{1,1}(\xi) &P_{1,2}(\xi)\\
(1-e^{-i\xi})^{m+n}P_{2,1}(\xi) &P_{2,2}(\xi)\end{matrix}\right].
\]
The above factorization of a matrix-valued filter allows us to theoretically study and construct multiwavelets/multiframelets with high vanishing moments from refinable vector functions, in almost the same way as the scalar case using the popular factorization technique in the scalar case (i.e., $r=1$).

On the other hand, as we will see later in the proof of Theorem~\ref{thm:qtf}, the almost orthogonal structure introduced in item (ii) of Theorem~\ref{thm:normalform:2} is the key to achieve the balancing property (see Section~\ref{sec:ffrt}) of the associated discrete multiframelet transform.

To prove Theorems~\ref{thm:normalform} and~\ref{thm:normalform:2},
we need a few auxiliary results.

\begin{lemma}\label{lem:U}
	Let $\wh{v}=[\wh{v_1},\ldots,\wh{v_r}]$ and $\wh{u}=[\wh{u_1},\ldots,\wh{u_r}]$ be $1\times r$ vectors of functions which are smooth near the origin such that $\wh{v}(0)\ne 0$ and $\wh{u}(0)\ne 0$. If $r\ge 2$, then for any positive integer $n\in \N$, there exists a strongly invertible $r\times r$ matrix $\wh{U}$ of $2\pi$-periodic trigonometric polynomials such that
\be \label{U:uv}
	 \wh{u}(\xi)=\wh{v}(\xi)\wh{U}(\xi)+\bo(|\xi|^n),\qquad \xi\to 0.
	\ee
\end{lemma}

\bp We first prove the claim for the special case $\wh{u}(\xi)=[1,0,\ldots,0]+\bo(|\xi|^n)$ as $\xi\to 0$.
Since $\wh{v}(0)\ne 0$, by permuting the entries of $\wh{v}$, we can assume that $\wh{v_1}(0)\ne 0$. Moreover, since $\wh{v}$ is smooth near the origin, we can find a $1\times r$ vector $\wh{\mathring{v}}$ of $2\pi$-periodic trigonometric polynomials such that $\wh{v}(\xi)=\wh{\mathring{v}}(\xi)+\bo(|\xi|^n)$ as $\xi \to 0$. Hence, without loss of generality, we assume that $\wh{v}$ is a vector of $2\pi$-periodic trigonometric polynomials. Since $\wh{v_1}(0)\ne 0$,
there exist $2\pi$-periodic trigonometric polynomials $\wh{w_j}(\xi), j=2,\ldots,r$ such that
\[
\wh{w_j}(\xi)=-\wh{v_j}(\xi)/\wh{v_1}(\xi)+\bo(|\xi|^n),\qquad \xi \to 0, j=2,\ldots,r.
\]
Define
\[
\wh{U_1}(\xi):=\left[ \begin{matrix} 1 &\wh{w_2}(\xi) &\cdots &\wh{w_r}(\xi)\\
0 &1 &\cdots &0\\
&\vdots &\ddots &\vdots\\
0 &0 &\cdots &1\end{matrix}\right].
\]
Since $\det (\wh{U_1}(\xi))=1$,
$\wh{U_1}$ is strongly invertible and
\be \label{vU1}
\wh{v}(\xi) \wh{U_1}(\xi)=[\wh{v_1}(\xi),0,\ldots,0]+\bo(|\xi|^n),\qquad \xi \to 0.
\ee
Note that $\wh{v_1}$ is a $2\pi$-periodic trigonometric polynomial with $\wh{v_1}(0)\ne 0$.
We now adopt an idea in
the proof of \cite[Theorem~2.1]{han09} to prove the claim.
Because there is no non-trivial common factor between the two $2\pi$-periodic trigonometric polynomials $\wh{v_1}(\xi)$ and $(1-e^{-i\xi})^{2n}$,
there must exist $2\pi$-periodic trigonometric polynomials $\wh{c}$ and $\wh{d}$ such that
\be \label{cdU2}
\wh{v_1}(\xi) \wh{c}(\xi)+(1-e^{-i\xi})^{2n} \wh{d}(\xi)=1\qquad \forall\; \xi\in \R.
\ee
Due to our assumption $r\ge 2$, we can define
\[
\wh{U_2}(\xi)=\left[ \begin{matrix}
\wh{c}(\xi) &-(1-e^{-i\xi})^n &0\\
(1-e^{-i\xi})^n \wh{d}(\xi) &\wh{v_1}(\xi) &0\\
0 &0 &I_{r-2}\end{matrix}\right].
\]
Using \eqref{vU1} and \eqref{cdU2}, we trivially conclude that
\[
\wh{v}(\xi) \wh{U_1}(\xi)\wh{U_2}(\xi)=
[\wh{v_1}(\xi),0,\ldots,0]\wh{U_2}(\xi)+\bo(|\xi|^n)
=[\wh{v_1}(\xi)\wh{c}(\xi), 0,\ldots,0]+\bo(|\xi|^n)
\]
as $\xi \to 0$. Due to \eqref{cdU2}, we have $\wh{v_1}(\xi)\wh{c}(\xi)=1-(1-e^{-i\xi})^{2n} \wh{d}(\xi)=1+\bo(|\xi|^n)$ as $\xi \to 0$ and $\det(\wh{U_2}(\xi))=1$. Hence, $\wh{U_2}$ is strongly invertible and
$\wh{v}(\xi) \wh{U_1}(\xi)\wh{U_2}(\xi)=[1,0,\ldots,0]+\bo(|\xi|^n)$ as $\xi \to 0$. The proof is completed for the special case of $\wh{u}$ by taking $\wh{U}(\xi):=\wh{U_1}(\xi)\wh{U_2}(\xi)$.

Generally, by what has been proved, there exist strongly invertible matrices $\wh{U_v}$ and $\wh{U_u}$ such that
\[
\wh{v}(\xi)\wh{U_v}(\xi)=[1,0,\ldots,0]+\bo(|\xi|^n),\qquad
\wh{u}(\xi)\wh{U_u}(\xi)=[1,0,\ldots,0]+\bo(|\xi|^n),\qquad \xi\to 0.
\]
Define $\wh{U}(\xi):=\wh{U_v}(\xi) \wh{U_u}(\xi)^{-1}$. Then $\wh{U}$ is strongly invertible and \eqref{U:uv} holds.
\ep

Note that Lemma~\ref{lem:U} often fails for $r=1$, since \eqref{U:uv} holds for $r=1$ if and only if $\wh{u}(\xi)/\wh{v}(\xi)=ce^{-ik\xi}+\bo(|\xi|^n)$ as $\xi\to 0$ for some $c\ne 0$ and $k\in \Z$.

The following result shows that the condition in \eqref{vguphi=1:new} of Theorem~\ref{thm:normalform} is also a necessary condition.

\begin{lemma}\label{lem:vguphi}
	Let $\dm\ge 2$ be an integer and $a\in \lrs{0}{r}{r}$ be a finitely supported matrix-valued filter.
	Let $\phi$ be an $r\times 1$ vector of
	compactly supported distributions satisfying $\wh{\phi}(\dm \xi)=\wh{a}(\xi)\wh{\phi}(\xi)$ with $\wh{\phi}(0)\ne 0$.
	Suppose that the filter $a$ has $m$ sum rules with respect to $\dm$ satisfying \eqref{sr} with a matching filter $\vgu\in \lrs{0}{1}{r}$ such that $\wh{\vgu}(0)\wh{\phi}(0)=1$.
	Then
	\be \label{vguphi=1:m}
	 \wh{\vgu}(\xi)\wh{\phi}(\xi)=1+\bo(|\xi|^m),\quad
	\xi\to 0.
	\ee
\end{lemma}

\bp The claim is essentially known in \cite[Proposition~3.2]{han03}. Here we provide a simple proof.
Since the filter $a$ satisfies \eqref{sr}, we have $\wh{\vgu}(\dm\xi ) \wh{a}(\xi)=\wh{\vgu}(\xi)+\bo(|\xi|^m)$ as $\xi \to 0$. Now we deduce from $\wh{\phi}(\dm\xi)=\wh{a}(\xi)\wh{\phi}(\xi)$ that
\be \label{rel:vguphi}
\wh{\vgu}(\dm \xi)\wh{\phi}(\dm\xi)
=\wh{\vgu}(\dm \xi) \wh{a}(\xi)\wh{\phi}(\xi)=
\wh{\vgu}(\xi)\wh{\phi}(\xi)+\bo(|\xi|^m), \quad \xi\to 0.
\ee
Considering the Taylor series of the function $\wh{\vgu}(\xi)\wh{\phi}(\xi)$ at $\xi=0$, since we assumed $\wh{\vgu}(0)\wh{\phi}(0)=1$ and $\dm\ge 2$, we can straightforwardly deduce from the above relation in \eqref{rel:vguphi}
that \eqref{vguphi=1:m} must hold.
\ep

\begin{lemma}\label{lem:moment}
	Let $m\in \N$ be a positive integer.
	Let $\wh{v}$ be a $1\times r$ row vector and $\wh{u}$ be an $r\times 1$ column vector such that all the entries of
	$\wh{v}$ and $\wh{u}$ are functions which are smooth near the origin such that
\be \label{vu=1:m} \wh{v}(\xi)\wh{u}(\xi)=1+\bo(|\xi|^m),\quad \xi \to 0.
\ee
	For any positive integer $n$, there must exist $1\times r$ vector $\wh{\mathring{v}}$ of functions which are smooth near the origin such that
	\be \label{vu=1:general}	 \wh{\mathring{v}}(\xi)=\wh{v}(\xi)+\bo(|\xi|^m)
	\quad \mbox{and}\quad \wh{\mathring{v}}(\xi)\wh{u}(\xi)=1+\bo(|\xi|^n),\quad \xi\to 0.
	\ee
\end{lemma}

\bp If $n\le m$, then we can simply take $\wh{\mathring{v}}:=\wh{v}$ and it follows directly from our assumption in \eqref{vu=1:m} that
\eqref{vu=1:general} trivially holds. So, we assume $n>m$.
We consider two cases $r=1$ and $r>1$.
If $r=1$, then $\wh{u}(0)\ne 0$.
Taking $\wh{\mathring{v}}(\xi):=1/\wh{u}(\xi)$, we see that \eqref{vu=1:general} is satisfied.

Suppose that $r>1$.
By Lemma~\ref{lem:U}, there exists a strongly invertible $r\times r$ matrix $\wh{U}$ such that $\wh{\breve{u}}(\xi):=\wh{U}(\xi)\wh{u}(\xi)=[1,0,\ldots,0]^\tp+\bo(|\xi|^n)$ as $\xi \to 0$.
We define $\wh{\breve{v}}(\xi)=
[\wh{\breve{v}_1}(\xi),\ldots,\wh{\breve{v}_r}(\xi)]
:= \wh{v}(\xi)\wh{U}(\xi)^{-1}$. Then it follows from \eqref{vu=1:m} that
\[
\wh{\breve{v}_1}(\xi)=
\wh{\breve{v}}(\xi)\wh{\breve{u}}(\xi)+\bo(|\xi|^n)
=\wh{v}(\xi)\wh{U}(\xi)^{-1} \wh{U}(\xi) \wh{u}(\xi)+\bo(|\xi|^n)
=\wh{v}(\xi)\wh{u}(\xi)=1+\bo(|\xi|^m), \quad \xi \to 0.
\]
We define $\wh{\mathring{v}}(\xi):=[1, \wh{\breve{v}_2}(\xi),\ldots, \wh{\breve{v}_r}(\xi)]\wh{U}(\xi)$.
Then
\[
\wh{\mathring{v}}(\xi) \wh{u}(\xi)=
[1,\wh{\breve{v}_2}(\xi),\ldots, \wh{\breve{v}_r}(\xi)] \wh{U}(\xi) \wh{u}(\xi)
=[1,\wh{\breve{v}_2}(\xi),\ldots, \wh{\breve{v}_r}(\xi)] \wh{\breve{u}}(\xi)
=1+\bo(|\xi|^n)
\]
as $\xi \to 0$. By $\wh{\breve{v}_1}(\xi)=1+\bo(|\xi|^m)$ as $\xi \to 0$ and noting $\wh{\breve{v}}(\xi)=\wh{v}(\xi)\wh{U}(\xi)^{-1}$, we have
\[
\wh{\mathring{v}}(\xi)=\wh{\breve{v}}(\xi)\wh{U}(\xi)+\bo(|\xi|^m)=\wh{v}(\xi)+\bo(|\xi|^m), \quad \xi\to 0.
\]
This completes the proof.
\ep

We are now ready to prove Theorem~\ref{thm:normalform}, which includes all the results on the normal form of a matrix-valued filter in \cite{han03,han09,hanbook,hm03} as special cases for dimension one.
Following the lines of our proof for Theorem~\ref{thm:normalform} below,
we also point out that Theorem~\ref{thm:normalform} can be generalized without much difficulty to multidimensional matrix-valued filters.

\begin{proof}[\textbf{Proof of Theorem~\ref{thm:normalform}}]
	Obviously, it suffices to prove the claims for $n \ge m$.
	By Lemma~\ref{lem:vguphi}, we see that \eqref{vguphi=1:m} holds. By our assumption in \eqref{vguphi=1:new} and the fact that $\wh{\phi}$ is smooth at every $\xi\in\R$ (because $\phi$ is a vector of compactly supported distributions), using Lemma~\ref{lem:moment}, without loss of generality we can assume that
	\be \label{vguphi=1}
	 \wh{\vgu}(\xi)\wh{\phi}(\xi)=1+\bo(|\xi|^n) \quad \mbox{and}
	\quad \wh{\mathring{\vgu}}(\xi)\wh{u_\phi}(\xi)=1+\bo(|\xi|^n),\qquad \xi \to 0.
	\ee

	Define $\wh{\breve{\vgu}}(\xi):=[1,0,\ldots,0]$.
	Since $\wh{\mathring{\vgu}}(0)\ne 0$ and $\wh{\vgu}(0)\ne0$, by Lemma~\ref{lem:U}, there exist strongly invertible $r\times r$ matrices $\wh{U_1}$ and $\wh{U_2}$ of $2\pi$-periodic trigonometric polynomials such that
	\be \label{vgu:breve}
	\wh{\breve{\vgu}}(\xi)=
	 \wh{\mathring{\vgu}}(\xi)\wh{U_1}(\xi)+\bo(|\xi|^{n})
	\quad\mbox{and}\quad
	 \wh{\vgu}(\xi)=\wh{\breve{\vgu}}(\xi) \wh{U_2}(\xi)+\bo(|\xi|^{n}),\quad
	\xi\to 0.
	\ee
	Define
	\[
	\wh{\breve{u}_\phi}(\xi):=
	\wh{U_1}(\xi)^{-1} \wh{u_\phi}(\xi),
	\qquad
	 \wh{\breve{\phi}}(\xi):=\wh{U_2}(\xi)\wh{\phi}(\xi),
	\qquad\mbox{and}\qquad
	 \wh{\breve{a}}(\xi):=\wh{U_2}(\dm\xi)\wh{a}(\xi) \wh{U_2}(\xi)^{-1}.
	\]
	Then it is trivial to check that $\wh{\breve{\phi}}(\dm \xi)=\wh{\breve{a}}(\xi)\wh{\breve{\phi}}(\xi)$ and $\breve{a}$ has $m$ sum rules with the matching filter $\breve{\vgu}$.
	Write $\breve{u}_\phi=[\breve{u}_1,\ldots, \breve{u}_r]^\tp$.
	Using \eqref{vguphi=1} and \eqref{vgu:breve} as well as $\wh{\breve{\vgu}}(\xi)=[1,0,\ldots,0]$,
	we observe that
	\[
	\wh{\breve{u}_1}(\xi)
	=\wh{\breve{\vgu}}(\xi) \wh{\breve{u}_\phi}(\xi)=
	 \wh{\mathring{\vgu}}(\xi)\wh{U_1}(\xi)
	\wh{U_1}(\xi)^{-1} \wh{u_\phi}(\xi)+\bo(|\xi|^{n})
	=\wh{\mathring{\vgu}}(\xi) \wh{u_\phi}(\xi)=1+\bo(|\xi|^n),\quad \xi\to 0.
	\]
	Write $\breve{\phi}=[\breve{\phi}_1,\ldots, \breve{\phi}_r]^\tp$. Since $\wh{\breve{\vgu}}(\xi)=[1,0,\ldots,0]$,
	we deduce from \eqref{vguphi=1} and \eqref{vgu:breve} that
	\[
	 \wh{\breve{\phi}_1}(\xi)=\wh{\breve{\vgu}}(\xi)\wh{\breve{\phi}}(\xi)=
	\wh{\vgu}(\xi) \wh{U_2}(\xi)^{-1}
	\wh{U_2}(\xi) \wh{\phi}(\xi)=
	 \wh{\vgu}(\xi)\wh{\phi}(\xi)=1+\bo(|\xi|^{n}),\quad \xi \to 0.
	\]
	There exist $2\pi$-periodic trigonometric polynomials $\wh{w_\ell}, \ell=2,\ldots,r$ such that
	\[
	 \wh{w_\ell}(\xi)=\wh{\breve{u}_\ell}(\xi)-
	 \wh{\breve{\phi}_\ell}(\xi)+\bo(|\xi|^{n}),\quad \xi \to 0, \ell=2, \ldots,r.
	\]
	Define
	\[
	\wh{U_3}(\xi):=\left[ \begin{matrix}
	1 &0 &\cdots &0\\
	\wh{w_2}(\xi) &1 &\cdots &0\\
	\vdots &\vdots &\ddots &\vdots\\
	\wh{w_r}(\xi) &0 &\cdots &1\end{matrix}
	\right].
	\]
	Since $\det(\wh{U_3}(\xi))=1$,
	the matrix $\wh{U_3}$ is strongly invertible. Moreover, by the definition of $\wh{w_\ell}$, we have
	\be \label{breve:phi}
	 \wh{U_3}(\xi)\wh{\breve{\phi}}(\xi)=
	 \wh{\breve{u}_\phi}(\xi)+\bo(|\xi|^n),\qquad \xi \to 0,
	\ee
	where we also used $\wh{\breve{u}_1}(\xi)=1+\bo(|\xi|^n)$ and $\wh{\breve{\phi}_1}(\xi)=1+\bo(|\xi|^{n})$ as $\xi \to 0$.
	
	Define $\wh{U}(\xi):=\wh{U_1}(\xi)\wh{U_3}(\xi) \wh{U_2}(\xi)$. Then $\wh{U}$ is strongly invertible and we now prove that all the claims in Theorem~\ref{thm:normalform} are satisfied. We first check \eqref{vguphi=1:new}. Using \eqref{vgu:breve} and $n\ge m$, we have
	\begin{align*}
	\wh{\vgu}(\xi)\wh{U}(\xi)^{-1}
	&=\wh{\vgu}(\xi) \wh{U_2}(\xi)^{-1} \wh{U_3}(\xi)^{-1} \wh{U_1}(\xi)^{-1}
	=\wh{\breve{\vgu}}(\xi) \wh{U_3}(\xi)^{-1}\wh{U_1}(\xi)^{-1}+\bo(|\xi|^{n})\\
	&=\wh{\breve{\vgu}}(\xi) \wh{U_1}(\xi)^{-1}+\bo(|\xi|^{n})
	 =\wh{\mathring{\vgu}}(\xi)+\bo(|\xi|^{n})=\wh{\mathring{\vgu}}(\xi)+\bo(|\xi|^{m}),
	\end{align*}
	as $\xi \to 0$, since the first row of $\wh{U_3}(\xi)^{-1}$ is $[1,0,\ldots,0]$ and $\wh{\breve{\vgu}}(\xi)=[1,0,\ldots,0]$.
	Similarly, by $\wh{\breve{\phi}}(\xi)=\wh{U_2}(\xi)\wh{\phi}(\xi)$ and using \eqref{breve:phi}, as $\xi\to 0$, we have
	\[
	\wh{U}(\xi)\wh{\phi}(\xi)=
	\wh{U_1}(\xi)\wh{U_3}(\xi) \wh{U_2}(\xi)\wh{\phi}(\xi)
	=\wh{U_1}(\xi)\wh{U_3}(\xi) \wh{\breve{\phi}}(\xi)
	=\wh{U_1}(\xi) \wh{\breve{u}_\phi}(\xi)+\bo(|\xi|^n)
	=\wh{u_\phi}(\xi)+\bo(|\xi|^n),
	\]
	where in the last identity we used the definition $\wh{\breve{u}_\phi}(\xi)=\wh{U_1}(\xi)^{-1}\wh{u_\phi}(\xi)$.
	This proves \eqref{normalform:general}.
	
	We now check items (i) and (ii). By $\wh{\phi}(\dm\xi)=\wh{a}(\xi)\wh{\phi}(\xi)$, we obviously have
	\[
	\wh{\mathring{\phi}}(\dm \xi)=\wh{U}(\dm \xi)\wh{\phi}(\dm \xi)
	=\wh{U}(\dm\xi) \wh{a}(\xi)\wh{\phi}(\xi)=
	\wh{\mathring{a}}(\xi) \wh{\mathring{\phi}}(\xi).
	\]
	Now by \eqref{vguphi=1:new}, we have
	 $\wh{\mathring{\phi}}(\xi)=\wh{U}(\xi)\wh{\phi}(\xi)=
	\wh{u_\phi}(\xi)+\bo(|\xi|^n)$ as $\xi \to 0$. This proves item (i).
	
	Since $\wh{U}$ is strongly invertible, the filter $\mathring{a}$ must be finitely supported. Since $a$ satisfies \eqref{sr} and \eqref{vguphi=1:new} holds, for $\gamma=0,\ldots, \dm-1$, we have
	
$$\begin{aligned}
	\wh{\mathring{\vgu}}(\dm \xi) \wh{\mathring{a}}(\xi+\tfrac{2\pi \gamma}{\dm})&=
	\wh{\vgu}(\dm \xi) \wh{U}(\dm \xi)^{-1}
	\wh{U}(\dm\xi) \wh{a}(\xi+\tfrac{2\pi \gamma}{\dm}) \wh{U}(\xi+\tfrac{2\pi \gamma}{\dm})^{-1}
	=\wh{\vgu}(\dm\xi) \wh{a}(\xi+\tfrac{2\pi \gamma}{\dm}) \wh{U}(\xi+\tfrac{2\pi \gamma}{\dm})^{-1}\\
	&=\begin{cases} \wh{\vgu}(\xi)\wh{U}(\xi)^{-1}+\bo(|\xi|^m)=\wh{\mathring{\vgu}}(\xi)+\bo(|\xi|^m), &\gamma=0,\\[0.3cm]
	\bo(|\xi|^m) & \gamma\in \{ 1,\ldots,\dm-1\},\end{cases}\qquad \xi\to 0,\end{aligned}$$
which proves item (ii).
\end{proof}

Finally, we prove Theorem~\ref{thm:normalform:2}, which plays an important role in our proof of Theorem~\ref{thm:qtf}.

\bp[\textbf{Proof of Theorem~\ref{thm:normalform:2}}]
We first prove item (i).
By Theorem~\ref{thm:normalform},
there exists a strongly invertible $r\times r$ matrix $\wh{U}$ of $2\pi$-periodic trigonometric polynomials such that all the claims of Theorem~\ref{thm:normalform} hold with $\wh{\mathring{\vgu}}(\xi)=
	[1,0,\ldots,0]$ and $\wh{u_\phi}(\xi)=[1,0,\ldots,0]^\tp$. Now by item (ii) of Theorem~\ref{thm:normalform}, we conclude that
	\be \label{a11}
	 \wh{\mathring{a}_{1,1}}(\xi+\tfrac{2\pi \gamma}{\dm})=\bo(|\xi|^m),\qquad \xi\to 0, \quad\gamma=1,\ldots,\dm-1
	\ee
	and
	\be \label{a12}
	 \wh{\mathring{a}_{1,2}}(\xi+\tfrac{2\pi \gamma}{\dm})=\bo(|\xi|^m),\qquad \xi\to 0, \quad\gamma=0,\ldots,\dm-1.
	\ee
\eqref{a11} is equivalent to $(1+e^{-i\xi}+\cdots+e^{-i(\dm-1)\xi})^m \mid \wh{\mra_{1,1}}(\xi)$, and \eqref{a12} is equivalent to
	$(1-e^{-i\dm \xi})^m \mid \wh{\mathring{a}_{1,2}}(\xi)$.
On the other hand, we have $\wh{\mathring{\phi}}(\xi)=\wh{U}(\xi)\wh{\phi}(\xi)=\wh{u_\phi}(\xi)=[1,0,\ldots,0]^\tp+\bo(|\xi|^n)$ as $\xi \to 0$, which is simply \eqref{normalform:phi}.
	Observing that $\wh{\mathring{\phi}}(\dm\xi)=\wh{\mathring{a}}(\xi)\wh{\mathring{\phi}}(\xi)$,
we conclude from \eqref{normalform:phi} that $\wh{\mathring{a}_{1,1}}(\xi)=1+\bo(|\xi|^n)$
	and $\wh{\mathring{a}_{2,1}}(\xi)=\bo(|\xi|^n)$ as $\xi\to 0$.
Thus \eqref{normalform} and \eqref{normalform:11} hold, and this proves item (i).
	
Next, we prove item (ii).  By Theorem~\ref{thm:normalform}, there exists
	a strongly invertible filter $V\in\lrs{0}{r}{r}$ such that
	 $$\wh{\vgu}(\xi)\wh{V}(\xi)=\wh{\mathring{\vgu}}(\xi)=[1,0,\dots,0]+\bo(|\xi|^{m}),\quad \wh{V}(\xi)^{-1}\wh{\phi}(\xi)=\wh{\mathring{\phi}}(\xi)=[1,0,\dots,0]^{\tp}+\bo(|\xi|^{\tilde{n}}),\quad \xi\to 0,$$
	where $\tilde{n}=\max(m,n)$. It follows from \er{moment:special} and the above identities that
	 \be\label{moment:special:2}[1,0,\dots,0]\wh{V}^{-1}(\xi)=\frac{\ol{\wh{\phi}(\xi)}^{\tp}}{\|\wh{\phi}(\xi)\|^2}+\bo(|\xi|^{m}),\quad [1,0,\dots,0]\ol{\wh{V}(\xi)}^{\tp}=\ol{\wh{\phi}(\xi)}^{\tp}+\bo(|\xi|^{\tilde{n}}),\quad \xi\to 0.\ee
	For $j=1,\dots,r$, denote $\wh{V_j}$ the $j$-th column of $\wh{V}$. It is easy to see from \eqref{moment:special:2} that $\wh{V_1}(\xi)=\wh{\phi}(\xi)+\bo(|\xi|^{\tilde{n}})$ as $\xi\to 0$.
	Set $\wh{u_1}(\xi):=\wh{V_1}(\xi)$ and choose $g_1\in \lp{0}$ such that $\wh{g_1}(\xi)=\frac{1}{\|\wh{\phi}(\xi)\|^2}+\bo(|\xi|^{\tilde{n}})$ as $\xi\to 0$. For $j=2,\dots,r$, define $u_j\in\lrs{0}{r}{1}$ and choose $g_{j}\in\lp{0}$ recursively via
\be\label{fj}\wh{u_j}(\xi)=\wh{V_j}(\xi)-\sum_{l=1}^{j-1}\wh{V_j}(\xi)^{\tp}\ol{\wh{u_l}(\xi)}\wh{g}_l(\xi)\wh{u_l}(\xi),
\ee
\be\label{gj}\wh{g_{j}}(\xi)=\frac{1}{\|\wh{u_{j}}(\xi)\|^2}+\bo(|\xi|^{\tilde{n}}),\quad\xi\to 0.\ee
Define
\be\label{W:ortho}\wh{U}(\xi)^{-1}:=[\wh{u_1}(\xi),\wh{u_2}(\xi),\dots,\wh{u_r}(\xi)]=[\wh{\phi}(\xi)+\bo(|\xi|^{\tilde{n}}),\wh{u_2}(\xi),\dots,\wh{u_r}(\xi)],\quad\xi\to 0.\ee
 By our construction, we have $\det(\wh{U}^{-1})=\det(\wh{V})$. This implies that $\wh{U}^{-1}$ is strongly invertible, and so is $\wh{U}$. For $j=1,\dots,r$, we have
	 $$\wh{u_j}(\xi)=\left(\wh{V_j}(\xi)-\sum_{l=1}^{j-1}\wh{V_j}(\xi)^{\tp}\ol{\wh{u_l}(\xi)}\frac{\wh{u_l}(\xi)}{\|\wh{u_l}(\xi)\|^2}\right)+\bo(|\xi|^{\tilde{n}}),\quad\xi\to 0.$$
	This means whenever $j\neq k$, we have
	 \be\label{W:ortho2}\ol{\wh{u_j}(\xi)}^{\tp}\wh{u_k}(\xi)=\bo(|\xi|^{\tilde{n}}),\quad\xi\to 0.\ee
	Note that the first column of $\wh{U}^{-1}$ is $\wh{V_1}$. It follows that
	 $$\ol{\wh{U}(\xi)}^{-\tp}\wh{U}(\xi)^{-1}=\DG\left(\|\wh{\phi}(\xi)\|^2, \|\wh{u_2}(\xi)\|^2,\dots,\|\wh{u_r}(\xi)\|^2\right)+\bo(|\xi|^{\tilde{n}}),\quad\xi\to 0.$$
Hence \er{eq:ortho} holds since $\tilde{n}=\max(m,n)$. Hence item (ii) is proved.
	
	Conversely, suppose that item (i) and \er{eq:ortho} hold. As $\wh{U}$ is strongly invertible, we see from \er{eq:ortho} that $\|\wh{\phi}(0)\|^2\neq 0$. Now use item (i), \er{eq:ortho} and $\max(m,n)\ge m$, we have
	 $$\begin{aligned}\wh{\vgu}(\xi)&=[1,0,\dots,0]\wh{U}(\xi)+\bo(|\xi|^m)=\frac{1}{\|\wh{\phi}(\xi)\|^2}[1,0,\dots,0]\ol{\wh{U}(\xi)}^{-\tp}\wh{U}(\xi)^{-1}\wh{U}(\xi)+\bo(|\xi|^m)\\
	 &=\frac{1}{\|\wh{\phi}(\xi)\|^2}[1,0,\dots,0]\ol{\wh{U}(\xi)}^{-\tp}+\bo(|\xi|^m)=\frac{\ol{\wh{\phi}(\xi)}^{\tp}}{\|\wh{\phi}(\xi)\|^2}+\bo(|\xi|^m),\quad\xi\to 0.
	\end{aligned}$$
	This proves \er{moment:special}. This completes the proof of Theorem~\ref{thm:normalform:2}.
\ep

\section{Properties of Multiframelet Transform for OEP-based Filter Banks}
\label{sec:ffrt}

In this section, we shall systematically study the discrete multiframelet transform using an OEP-based dual $\dm$-multiframelet filter bank. Then we study its balancing property and possible shortcomings.

To state a discrete multiframelet transform, let us recall some necessary definitions.
By $(\sq)^{s\times r}$ we denote the linear space of all sequences $v: \Z \rightarrow \C^{s\times r}$. For a positive integer $\dm$ and a filter $a\in \lrs{0}{r}{r}$, the \emph{subdivision operator} $\sd_{a,\dm}$ and the \emph{transition operator} $\tz_{a,\dm}$ are defined to be
\[
[\sd_{a,\dm} v](n):=\sqrt{\dm} \sum_{k\in \Z} v(k) a(n-\dm k),\quad
[\tz_{a,\dm} v](n):=\sqrt{\dm} \sum_{k\in \Z} v(k)
\ol{a(k-\dm n)}^\tp,\quad n\in \Z,\quad v\in (\sq)^{s\times r}.
\]
Moreover, we define a new filter $a^\star$ by $\wh{a^\star}(\xi):=\ol{\wh{a}(\xi)}^\tp$, that is,
$a^\star(k):=\ol{a(-k)}^\tp$ for all $k\in \Z$. The convolution $v*a$ is defined to be
\[
[v*a](n):=\sum_{k\in \Z} v(k) a(n-k),\quad n\in \Z.
\]
Hence, it is easy to verify that $\sd_{a,\dm} v=\sqrt{\dm}[v(\dm\cdot)]*a$ and $\tz_{a,\dm} v=\sqrt{\dm}[v*a^\star](\dm \cdot)$. Therefore, both operators can be implemented efficiently using convolutions.

\subsection{Multi-level Discrete Multiframelet Transform and the Balancing Property}
Let $a,\tilde a, \theta, \tilde{\theta}\in \lrs{0}{r}{r}$ and
$b,\tilde{b}\in \lrs{0}{s}{r}$ be finitely supported filters.
Define a filter $\Theta:=\theta^\star*\tilde{\theta}$, i.e., $\wh{\Theta}(\xi):=\ol{\wh{\theta}(\xi)}^\tp \wh{\tilde{\theta}}(\xi)$.
We now state the discrete multiframelet transform using the finitely supported matrix-valued filter bank $(\{a;b\},\{\tilde{a};\tilde{b}\})_{\Theta}$. For $J\in \N$, the $J$-level discrete multiframelet decomposition is defined by
\be\label{frame:decomp}
v_j:=\tz_{a,\dm} v_{j-1},\quad w_j:=\tz_{b,\dm} v_{j-1}, \qquad j=1,\ldots, J,
\ee
where $v_0\in (\sq)^{1\times r}$ is a (vector-valued) input signal/data. The $J$-level discrete framelet reconstruction procedure is as follows:
\begin{enumerate}
\item[Step 1.] Compute $\tilde{v}_J:=v_J*\Theta$, where the convolution is
	$v_J*\Theta:=\sum_{k\in \Z} v_J(k)\Theta(\cdot-k)$.
\item[Step 2.] Recursively compute $\tilde{v}_j, j=J,\ldots,1$ by
	\be\label{sd:v0}
	\tilde{v}_{j-1}:= \sd_{\tilde{a},\dm} \tilde{v}_j+\sd_{\tilde{b},\dm} w_{j}, \qquad j=J,\ldots,1.
	\ee
\item[Step 3.] Recover $\mathring{v}_0$ through the deconvolution from $\tilde{v}_0=\mathring{v}_0*\Theta$.
\end{enumerate}

We shall address
several important issues on a (multi-level) discrete multiframelet transform such as the perfect reconstruct property, the balancing property, and balanced vanishing moments.
Note that Step 3 in the $J$-level discrete multiframelet reconstruction may have infinitely many solutions, a unique solution, or no solution at all.
We say that a discrete multiframelet transform has \emph{the generalized perfect reconstruction property} if any original input signal $v_0$ can be reconstructed as one of the solutions $\mathring{v}_0$ of $\tilde{v}_0=\mathring{v}_0*\Theta$ in
Step 3.

To analyze a $J$-level discrete multiframelet transform using the filter bank $(\{a;b\},\{\tilde{a};\tilde{b}\})_{\Theta}$, we define:
\begin{itemize}
	
	\item The \emph{$J$-level discrete multiframelet analysis/decomposition operator}:
	$$\cW_J:(\sq)^{1\times r}\rightarrow (\sq)^{1\times (sJ+r)},\quad \cW_J(v)=(\tz_{{b},\dm} v,\tz_{{b},\dm}\tz_{{a},\dm}v,\dots,\tz_{{b},\dm}\tz_{{a},\dm}^{J-1}v,\tz_{{a},\dm}^Jv).$$
	Define $\cW:=\cW_1$ as the one-level analysis/decomposition operator.
	
	\item The \emph{$J$-level discrete multiframelet synthesis/reconstruction operator}:
	$$\ctV_J:(\sq)^{1\times (sJ+r)}\rightarrow (\sq)^{1\times r},\quad \ctV_J(\mrw_1,\mrw_2,\dots,\mrw_J,\mrv_J)=\tilde{v}_0,$$
	for all $\mrw_j\in(\sq)^{1\times s}$ and $\mrv_J\in(\sq)^{1\times r}$, where $\tilde{v}_{j-1},j=J,\dots,1$ are recursively computed via
	$$\tilde{v}_{j-1}:= \sd_{\tilde{a},\dm} \tilde{v}_j+\sd_{\tilde{b},\dm} \mrw_{j}, \qquad j=J,\ldots,1,$$
	with $\tilde{v}_J:=\mrv_J$. Define $\ctV:=\ctV_1$ as the one-level synthesis/reconstruction operator.

	\item The \emph{$J$-level convolution operator}:
	$$C_{\Theta;J}:(\sq)^{1\times (sJ+r)}\rightarrow (\sq)^{1\times (sJ+r)},\quad C_{\Theta;J}(\mrw_1,\mrw_2,\dots,\mrw_J,\mrv_J)=(\mrw_1,\mrw_2,\dots,\mrw_J,\mrv_J*\Theta),$$
	for all $\mrw_j\in(\sq)^{1\times s}$ and $\mrv_J\in(\sq)^{1\times r}$.
\end{itemize}

We observe that the $J$-level discrete multiframelet transform using the filter bank $(\{a;b\},\{\tilde{a};\tilde{b}\})_{\Theta}$
has the generalized perfect reconstruction property for an input signal $v\in (\sq)^{1\times r}$ if and only if
\be\label{gpr:1}\ctV_J C_{\Theta;J} \cW_J(v)=C_{\Theta}(v),\ee
where $C_{\Theta}$ is the \emph{convolution operator} $C_{\Theta}(v)=v*\Theta$. Moreover, by
$$
\cW_J=(\ID_{(\sq)^{1\times s(J-1)}}\otimes \cW)\cdots(\ID_{(\sq)^{1\times s}}\otimes \cW)\cW
$$
and
$$
\ctV_J=\ctV(\ID_{(\sq)^{1\times s}}\otimes \ctV)\cdots(\ID_{(\sq)^{1\times s(J-1)}}\otimes \ctV),
$$
we see that the $J$-level multiframelet transform
has the generalized perfect reconstruction property for $v$ if and only if the one-level multiframelet transform does, that is,
\be \label{pr:1}
\sd_{\tilde{a},\dm}([\tz_{a,\dm} v]*\Theta)+
\sd_{\tilde{b},\dm} (\tz_{b,\dm} v)=v*\Theta.
\ee
%


Following the approach in \cite{han13,hanbook} for scalar framelets,
we now provide the necessary and sufficient conditions for the generalized perfect reconstruction property of a discrete multiframelet transform as follows:

\begin{theorem}\label{thm:ffrt:pr:general}
	Let $a,\tilde a, \theta, \tilde{\theta}\in \lrs{0}{r}{r}$ and
	$b,\tilde{b}\in \lrs{0}{s}{r}$ be finitely supported filters. Define $\Theta:=\theta^\star*\tilde{\theta}$. The following statements are equivalent to each other:
	
	\begin{enumerate}
		\item[(i)] For any $J\in\N$, the $J$-level discrete multiframelet transform using the filter bank $(\{a;b\},\{\tilde{a};\tilde{b}\})_{\Theta}$ has the generalized perfect reconstruction property for any input signal $v\in (\sq)^{1\times r}$ (or for any input signal $v\in \lrs{0}{1}{r}$).
		
		\item[(ii)] $(\{a;b\},\{\tilde{a};\tilde{b}\})_{\Theta}$ is an OEP-based dual $\dm$-framelet filter bank satisfying \eqref{dffb}.	 \end{enumerate}
\end{theorem}

\bp
(i) $\Rightarrow$ (ii): The perfect reconstruction property of the one-level discrete multiframelet transform for $v_0\in \lrs{0}{1}{r}$ is equivalent to \eqref{pr:1}.
For $v\in \lrs{0}{1}{r}$, we observe
\[
\wh{\sd_{a,\dm} v}(\xi)=\dm^{1/2} \wh{v}(\dm \xi) \wh{a}(\xi) \quad \mbox{and}\quad
\wh{\tz_{a,\dm} v}(\dm \xi)=\dm^{-1/2}
\sum_{\gamma=0}^{\dm-1} \wh{v}(\xi+\tfrac{2\pi \gamma}{\dm})
\ol{\wh{a}(\xi+\tfrac{2\pi\gamma}{\dm})}^\tp.
\]
Therefore, in the frequency domain, \eqref{pr:1} is equivalent to
\be \label{pr:2}
\sum_{\gamma=0}^{\dm-1}
\wh{v}(\xi+\tfrac{2\pi \gamma}{\dm})
\left[\ol{\wh{a}(\xi+\tfrac{2\pi\gamma}{\dm})}^\tp\wh{\Theta}(\dm \xi)
\wh{\tilde{a}}(\xi)+\ol{\wh{b}(\xi+\tfrac{2\pi\gamma}{\dm})}^\tp
\wh{\tilde{b}}(\xi)\right]=\wh{v}(\xi)\wh{\Theta}(\xi).
\ee
Using the same argument as in \cite[Theorem~2.1]{han13} and \cite[Theorem~1.1.1]{hanbook} by selecting $v$ as a sequence of Dirac sequences in \eqref{pr:2}, we deduce from \eqref{pr:2} that \eqref{pr:1} implies \eqref{dffb}, that is, $(\{a;b\},\{\tilde{a};\tilde{b}\})_{\Theta}$  must be an OEP-based dual $\dm$-framelet filter bank.

(ii) $\Rightarrow$ (i): Suppose that $(\{a;b\},\{\tilde{a};\tilde{b}\})_{\Theta}$  is an OEP-based dual $\dm$-framelet filter bank satisfying \eqref{dffb}. Then \eqref{pr:2} must hold for all $v\in \lrs{0}{1}{r}$. In the time domain, \eqref{pr:2} is equivalent to \eqref{pr:1}, which further implies \eqref{gpr:1}. This proves the generalized perfect reconstruction property for any $v\in\lrs{0}{1}{r}$. Now for arbitrary $v\in(\sq)^{1\times r}$, one can use the locality of the subdivision and transition operators (see the proof of \cite[Theorem~2.1]{han13} and \cite[Theorem~1.1.1]{hanbook}) to prove that \eqref{pr:2} holds for all $v\in (\sq)^{1\times r}$. This completes the proof.
\ep

However, if the deconvolution in Step 3 of the $J$-level discrete multiframelet reconstruction has infinitely many solutions or no solution at all, without any extra information on the input signal, then one cannot exactly recover the original input signal $v_0$ from Step 3. Hence, we say that a discrete multiframelet transform has \emph{the perfect reconstruction property} if any original input signal $v_0$ can be reconstructed as the unique solution $\mathring{v}_0$ of $\tilde{v}_0=\mathring{v}_0*\Theta$ in
Step 3.

To study the perfect reconstruction property of a discrete multiframelet transform, we need the following auxiliary result.

\begin{lemma}\label{lem:conv}
	Let $\Theta\in \lrs{0}{r}{r}$ be a finitely supported filter. Define the convolution operator $C_\Theta$ by $C_\Theta (v):=v*\Theta$ for any sequences $v\in (\sq)^{1\times r}$. Then
	\begin{enumerate}
		\item[(1)] The mapping $C_\Theta: (l_\infty(\Z))^{1\times r}\rightarrow (l_\infty(\Z))^{1\times r}$ is injective (or bijective) if and only if $\det(\wh{\Theta}(\xi))\ne 0$ for all $\xi \in \R$.
		\item[(2)] The mapping $C_\Theta: (l_{si}(\Z))^{1\times r}\rightarrow (l_{si}(\Z))^{1\times r}$ is injective (or bijective) if and only if $\det(\wh{\Theta}(\xi))\ne 0$ for all $\xi \in \R$, where $l_{si}(\Z)$ denotes the space of all slowly increasing sequences, i.e., $v\in l_{si}(\Z)$ if $(1+|\cdot|^2)^{-m} v\in l_\infty(\Z)$ for some $m\in \N$.
		\item[(3)] The mapping $C_\Theta: (\sq)^{1\times r}\rightarrow (\sq)^{1\times r}$ is injective (or bijective) if and only if $\det(\wh{\Theta}(\xi))$ is a nontrivial monomial (i.e., $\det(\wh{\Theta}(\xi))=ce^{-im\xi}$ for some $m\in \Z$ and $c\in \C\bs\{0\}$).
	\end{enumerate}
\end{lemma}

\bp We first prove items (1) and (2) simultaneously. Let $V_0$ be either $(l_\infty(\Z))^{1\times r}$ or $(l_{si}(\Z))^{1\times r}$. Suppose that $C_\Theta:V_0\rightarrow V_0$ is injective, but $\det(\wh{\Theta}(\xi_0))=0$ for some $\xi_0\in\R$. We start with the case $r=1$. In this case, we have
$$0=\wh{\Theta}(\xi_0)=\sum_{k\in\Z}\Theta(k)e^{-ik\xi_0}.$$
Let $v\in  l_\infty(\Z)$ be defined by
\be\label{v0}
v(k)=e^{-ik\xi_0},\quad\forall k\in\Z.
\ee
By definition, we have
$$(v*\Theta)(n)=\sum_{k\in\Z}v(k)\Theta(n-k)=e^{-in\xi_0}\sum_{k\in\Z}e^{-i(k-n)\xi_0}\Theta(n-k)=e^{-in\xi_0}\wh{\Theta}(\xi_0)=0,\quad\forall n\in\Z.$$
So we find a non-zero sequence $v$ with $v*\Theta=0$. Hence $C_\Theta$ is not injective, which is a contradiction. So we must have $\wh{\Theta}(\xi)\neq 0$ for all $\xi\in\R$.

Now we consider $r>1$. As $\det(\wh{\Theta}(\xi_0))=0$, we can find an invertible $r\times r$ matrix $Q$ such that all elements in the first row of $Q\wh{\Theta}(\xi_0)$ are zero. Let $v\in  l_\infty(\Z)$ be defined as \er{v0}, and let $u\in (l_\infty(\Z))^{1\times r}$ be defined by
$$
u(k)=[v(k),0,\dots,0]Q,\quad\forall k\in\Z.
$$
It follows immediately that $u*\Theta=0$, which again contradicts our assumption that $C_\Theta$ is injective. Hence, $\det(\wh{\Theta}(\xi))\neq 0$ for all $\xi\in\R$.

Conversely, suppose that $\det(\wh{\Theta}(\xi))\neq 0$ for all $\xi\in\R$.
Then the filter $\Theta^{-1}$, which is defined by $\wh{\Theta^{-1}}(\xi):=\wh{\Theta}(\xi)^{-1}$, must be well defined and has exponential decaying coefficients. Consequently, we can deduce that
\be\label{conv:sq}
(v*\Theta)*\Theta^{-1}=v*(\Theta*\Theta^{-1})=v=v*(\Theta^{-1}*\Theta)=(v*\Theta^{-1})*\Theta
\ee
and $v*\Theta^{-1}\in V_0$ for all $v\in V_0$. Hence, $C_\Theta$ must be bijective.
This proves items (1) and (2).

Finally, we prove item (3). Suppose that $C_\Theta:(\sq))^{1\times r}\rightarrow (\sq)^{1\times r}$ is injective, but $\det(\wh{\Theta}(\xi))$ is not a non-trivial monomial. Then there exist $\xi_0\in\C$ such that $\det(\wh{\Theta}(\xi_0))=0$. Then by applying the same argument as in the proof of item (1), we conclude that $C_\Theta$ is not injective, which is a contradiction.

Conversely, if $\det(\wh{\Theta}(\xi))$ is a non-trivial monomial, then $\Theta$ is strongly invertible and $\Theta^{-1}\in\lrs{0}{r}{r}$. Consequently, \eqref{conv:sq} must hold for all $v\in (\sq)^{1\times r}$. Hence,
$C_\Theta$ is bijective.
\ep

Now Lemma~\ref{lem:conv} yields the following characterization on the perfect reconstruction property of a discrete multiframelet transform employing an OEP-based dual framelet filter bank.

\begin{theorem}\label{thm:ffrt:pr}
	Let $a,\tilde a, \theta, \tilde{\theta}\in \lrs{0}{r}{r}$ and
	$b,\tilde{b}\in \lrs{0}{s}{r}$ be finitely supported filters. Define $\Theta:=\theta^\star*\tilde{\theta}$.
	Let $V_0=(l_{si}(\Z))^{1\times r}$ (or respectively, $V_0=(\sq)^{1\times r}$).
	For any $J\in \N$, the $J$-level discrete multiframelet transform using the filter bank $(\{a;b\},\{\tilde{a};\tilde{b}\})_{\Theta}$ has the perfect reconstruction property for any input signal from $V_0$ if and only if
	\begin{enumerate}
		\item[(i)] $(\{a;b\},\{\tilde{a};\tilde{b}\})_{\Theta}$ is an OEP-based dual $\dm$-framelet filter bank satisfying \eqref{dffb};
		
		\item[(ii)] $\det(\wh{\Theta}(\xi))\ne 0$ for all $\xi \in \R$ (or respectively, $\det(\wh{\Theta}(\xi))$ is a non-trivial monomial), where $\wh{\Theta}(\xi):=\ol{\wh{\theta}(\xi)}^\tp \wh{\tilde{\theta}}(\xi)$.
	\end{enumerate}
\end{theorem}

Except for the examples in \cite{han09}, all constructed OEP-based dual framelet filter banks with non-trivial $\Theta$ do not satisfy item (ii) of Theorem~\ref{thm:ffrt:pr}. For the convenience of the reader,
we now present two concrete examples of tight framelet filter banks such that
item (ii) of Theorem~\ref{thm:ffrt:pr} fails.

\begin{exmp}
Let $\dm=2$ and
consider the B-spline filter $a^B_{2,2}\in \lp{0}$:
	 $$\wh{a^B_{2,2}}(\xi)=\frac{1}{4}(1+e^{-i\xi})^2,\quad\xi\in\R.$$
	It is well known that $a^B_{2,2}$ is the mask associated to the refinable function $\wh{B_2}(\xi)=\left(\frac{1-e^{-i\xi}}{i\xi}\right)^2$. That is, $\wh{B_2}(2\xi)=\wh{a^B_2}(\xi)\wh{B_2}(\xi)$ for all $\xi\in\R$. With
$$
\wh{\theta}(\xi)=(1+e^{-i\xi})/2,\quad \wh{\Theta}(\xi)=\ol{\wh{\theta}(\xi)}\wh{\theta}(\xi)=(2+e^{-i\xi}+e^{i\xi})/4,
$$
one can construct a tight $2$-framelet filter bank $\{a^B_{2,2};b\}_\Theta$ satisfying
$$
\ol{\wh{a^B_{2,2}}(\xi)}\wh{\Theta}(2\xi)\wh{a^B_{2,2}}(\xi)+\ol{\wh{b}(\xi)}^{\tp}\wh{b}(\xi)=\wh{\Theta}(\xi),\quad \ol{\wh{a^B_{2,2}}(\xi)}\wh{\Theta}(2\xi)\wh{a^B_{2,2}}(\xi+\pi)+\ol{\wh{b}(\xi)}^{\tp}\wh{b}(\xi+\pi)=0$$
	for all $\xi\in\R$, where $b:=[b_1,b_2]^{\tp}\in\lrs{0}{2}{1}$ is given by
	 $$\wh{b_1}(\xi)=\frac{\sqrt{2}}{8}(e^{-i\xi}-1)(e^{-i\xi}+1)^3,\quad \wh{b_2}(\xi)=\frac{\sqrt{2}}{4}(e^{-2i\xi}-1),\quad\xi\in\R.$$
Define
$\wh{\eta}(\xi):=\wh{\theta}(\xi)\wh{B_2}(\xi)$ and
$\wh{\psi}(\xi):=\wh{b}(\xi/2)\wh{B_2}(\xi/2)$ for all $\xi\in\R$.
Note that $\wh{\eta}(0)=\wh{\theta}(0)\wh{B_2}(0)=1$ and $\psi$ has one vanishing moment.
By Theorem~\ref{thm:df},
$\{\eta;\psi\}$ forms a compactly supported tight $2$-framelet in $\Lp{2}$. However, because $\wh{\Theta}(\pi)=0$, Theorem~\ref{thm:ffrt:pr} tells us that the tight $2$-framelet filter bank $\{a^B_{2,2};b\}_{\Theta}$ cannot have the perfect reconstruction property for certain input signals.
\end{exmp}

\begin{exmp}Let $\phi_1(x)=B_2(\cdot-1)=\max(1-|x|,0)$ for all $x\in\R$. Then $\phi:=[\phi_1,0]^{\tp}$ is a $2$-refinable vector of compactly supported functions satisfying $\wh{\phi}(2\xi)=\wh{a}(\xi)\wh{\phi}(\xi)$ with $$
\wh{a}(\xi)=\frac{1}{4}\begin{bmatrix}e^{-i\xi}+2+e^{i\xi} &0\\
	0 & 1\end{bmatrix},\quad\xi\in\R.$$
	With	 $$\wh{\theta}(\xi)=\begin{bmatrix}(1+e^{-i\xi})/2 & 0\\
	0 & 0\end{bmatrix},\quad \wh{\Theta}(\xi)=\ol{\wh{\theta}(\xi)}^{\tp}\wh{\theta}(\xi),\quad\xi\in\R,$$
	we can construct a tight $2$-multiframelet filter bank $\{a;b\}_\Theta$ satisfying
	 $$\ol{\wh{a}(\xi)}^{\tp}\wh{\Theta}(2\xi)\wh{a}(\xi)+\ol{\wh{b}(\xi)}^{\tp}\wh{b}(\xi)=\wh{\Theta}(\xi),\quad \ol{\wh{a}(\xi)}^{\tp}\wh{\Theta}(2\xi)\wh{a}(\xi+\pi)+\ol{\wh{b}(\xi)}^{\tp}\wh{b}(\xi+\pi)=0$$
	for all $\xi\in\R$, where $b\in\lrs{0}{2}{2}$ is given by
	 $$\wh{b}(\xi)=\frac{\sqrt{2}}{8}\begin{bmatrix}
	(1-e^{-i\xi})(1+e^{-i\xi})^3 & 0\\
	2e^{-i\xi}(e^{-2i\xi}-1) & 0
	\end{bmatrix},\quad\xi\in\R.$$
Define $\wh{\eta}(\xi):=\wh{\theta}(\xi)\wh{\phi}(\xi)$ and $\wh{\psi}(\xi)=\wh{b}(\xi/2)\wh{\phi}(\xi/2)$.
Note that $\|\wh{\eta}(0)\|^2=\ol{\wh{\phi}(0)}^\tp\wh{\Theta}(0)\wh{\phi}(0)=1$ and $\psi$ has one vanishing moment.
By Theorem~\ref{thm:df},
$\{\eta;\psi\}$ is a tight $2$-multiframelet in $\Lp{2}$. However, $\det(\wh{\Theta})$ is identically zero, which clearly does not satisfy item (ii) of Theorem~\ref{thm:ffrt:pr}. As a consequence, the tight $2$-multiframelet filter bank $\{a;b\}_{\Theta}$ does not have the perfect reconstruction property.
\end{exmp}

Next, we discuss the discrete vanishing moments and the balancing property of a discrete multiframelet transform.
``Smooth'' signals are often modelled by polynomial sequences.
For $m\in \N$, by $\PL_{m-1}$ we denote the space of all polynomial sequences of degree less than $m$. The sparsity of a discrete multiframelet transform is described by its ability to have zero framelet coefficients $w_j$ for polynomial input data.
The input to a discrete multiframelet transform is a vector sequence in $(\sq)^{1\times r}$, while most data in applications are scalar-valued, i.e., in $\sq$. Hence, we have to convert a scalar sequence into a vector sequence by using the \emph{standard vector conversion operator}
\[
\mathring{E}: \sq\rightarrow (\sq)^{1\times r}\quad \mbox{with}\quad
[\mathring{E}v](k):=[v(rk),v(rk+1),\ldots, v(rk+(r-1))],\qquad k\in \Z.
\]
Note that $\mathring{E}$ is a linear bijective mapping.
Let $(\{a;b\},\{\tilde{a};\tilde{b}\})_\Theta$ be an OEP-based dual $\dm$-framelet filter bank and $(\{\eta;\psi\},\{\tilde{\eta};\tilde{\psi}\})$ be its corresponding dual $\dm$-framelet.
Define $m:=\sr(a,\dm)$ be its sum rule order.
Ideally, since the multi-level discrete multiframelet transform is recursive, to have sparsity of a multiframelet transform, we hope that
\be \label{bp:highpass}
\tz_{b,\dm} \mathring{E} (\pp)=0, \qquad \forall\; \pp\in \PL_{m-1}
\ee
and
\be \label{bp:lowpass}
\tz_{a,\dm} \mathring{E}(\pp)\in \mathring{E}(\PL_{m-1}),\qquad \forall\; \pp\in \PL_{m-1}.
\ee
The condition in \eqref{bp:lowpass} guarantees that the output signal $\tz_{a,\dm} \mathring{E}(\pp)$ is still in the space $\mathring{E}(\PL_{m-1})$ for any input data $\pp\in \mathring{E}(\PL_{m-1})$.
The condition in \eqref{bp:highpass} preserves sparsity for all levels, that is, the framelet coefficients $w_j:=\tz_{b,\dm} \tz_{a,\dm}^{j-1} \mathring{E}(\pp)=0$ for all $\pp\in \PL_{m-1}$ and $j\in \N$.

Hence, we say that a filter $b$ has \emph{$m$ balanced vanishing moments} if \eqref{bp:highpass} holds. Moreover, we define $\bvmo(b,\dm):=m$ with $m$ being the largest possible integer such that \eqref{bp:highpass} holds.
Similarly, we say that a discrete multiframelet transform (using the filter bank $(\{a;b\},\{\tilde{a};\tilde{b}\})_{\Theta}$) or a filter bank $\{a;b\}$ has \emph{$m$ balancing order} with respect to the dilation factor $\dm$ if both \eqref{bp:highpass} and \eqref{bp:lowpass} hold. In particular, we define $\bpo(\{a;b\},\dm):=m$ with $m$ being the largest such integer satisfying both  \eqref{bp:highpass} and \eqref{bp:lowpass}. We observe that $\bvmo(b)\le \vmo(\psi)$ and $\bpo(\{a;b\},\dm)\le \bvmo(b,\dm)$.
For the case $r=1$, we always have
$\bpo(\{a;b\},\dm)=\bvmo(b,\dm)=\vmo(\psi)$.
But for $r>1$, it was first observed in \cite{lv98} that $\bpo(\{a,b\},\dm)<\vmo(\psi)$ often happens. This reduced sparsity hurdles the applications of multiwavelets and multiframelets. How to remedy this shortcoming has been extensively studied in the function setting in \cite{cj00,sel00} and in the setting of discrete multiframelet transforms in \cite{han09,han10,hanbook}.

The following result on characterizing the balanced vanishing moments and the balancing property is known (e.g., see \cite[Lemma 7.6.3]{hanbook} or \cite[Theorem 4.4]{han09}).

\begin{theorem}\label{prop:bp}
	Let $\dm\ge 2$ be a positive integer.
	Let $a\in \lrs{0}{r}{r}$ and $b\in \lrs{0}{s}{r}$. Define
	\be \label{vgu:special}
	\wh{\Vgu}(\xi):=[1, e^{i\xi/r},\ldots, e^{i (r-1)\xi/r}].
	\ee
	Then
	\begin{enumerate}
		\item[(1)] The filter $b$ has $m$ balanced vanishing moments satisfying \eqref{bp:highpass} if and only if
		\be \label{cond:bvmo}
		 \wh{\Vgu}(\xi)\ol{\wh{b}(\xi)}^\tp=\bo(|\xi|^m),\qquad \xi \to 0.
		\ee
		\item[(2)] The filter bank $\{a;b\}$ has $m$ balancing order satisfying both \eqref{bp:highpass} and \eqref{bp:lowpass} if and only if \eqref{cond:bvmo} holds and there exists $c\in\lp{0}$ with $\wh{c}(0)\neq 0$ such that
		
		\be \label{cond:ao}
		 \wh{c}(\xi)\wh{\Vgu}(\xi)\ol{\wh{a}(\xi)}^\tp=\wh{\Vgu}(\dm\xi)+\bo(|\xi|^m),\quad\xi\to 0.
		\ee
		
	\end{enumerate}
\end{theorem}

\subsection{Difficulties in the Construction of OEP-based Multiframelets.}
In the following we discuss the difficulties involved in constructing multiframelets through OEP. For the scalar case $r=1$,
it is well known (e.g., see \cite[Proposition~3.3.1]{hanbook}) that a tight $\dm$-framelet $\{\eta;\psi\}$ constructed through OEP in Theorem~\ref{thm:df} must satisfy
\[
\vmo(\psi)=\min(\sr(a,\dm),\tfrac{1}{2}\vmo(u_{a,\Theta})) \quad \mbox{with}\quad
\wh{u_{a,\Theta}}(\xi):=\wh{\Theta}(\xi)-\wh{\Theta}(\dm\xi) |\wh{a}(\xi)|^2,
\]
where $n:=\vmo(u_{a,\Theta})$ is the largest integer satisfying $\wh{u_{a,\Theta}}(\xi)=\bo(|\xi|^n)$ as $\xi \to 0$.
For the B-spline filter $a^B_{m,\dm}$, we have
$$\wh{a^B_{m,\dm}}(\xi)=\dm^{-m}(1+e^{-i\xi}+\dots+e^{-i(\dm-1)\xi})^m,\quad m\in\N.$$
It is easy to see that $\sr(a^B_{m,\dm},\dm)=m$ and $\vmo(1-|\wh{a^B_{m,\dm}}(\xi)|^2)=2$. Consequently, any tight $\dm$-framelet derived from a B-spline refinable function in Theorem~\ref{thm:df} with the trivial choice $\wh{\Theta}(\xi)=1$ has no more than one vanishing moment.
The main purpose of OEP (see
\cite{chs02,dh04,dhrs03})
is to increase vanishing moments of the framelet generator $\psi$ by properly choosing $\Theta$ such that
$\vmo(u_{a,\Theta})=2\sr(a,\dm)$.
A lot of compactly supported scalar tight framelets and scalar dual framelets with the highest possible vanishing moments (i.e., $\vmo(\psi)=\sr(a,\dm)$) have been constructed through the OEP in the literature, to mention only a few, see \cite{chs02,dh04,dhrs03,dh18pp,ds13,han15, hanbook,hm05} and references therein.
In particular, see Chapter~3 of \cite{hanbook} for a comprehensive study of scalar tight or dual framelets.
Except the examples in \cite{han09}, however, except the trivial choice $\wh{\Theta}(\xi)=I_r$, the filters $\Theta$ in all known constructions of OEP-based framelets are not strongly invertible and therefore, their associated discrete framelet transforms are not compact, which seriously hinders their applications.

On the other hand, constructing tight framelets through OEP in Theorem~\ref{thm:df} is much more difficult when $r>1$. A necessary condition to construct an OEP-based tight multiframelet is the positive semi-definiteness of the matrix $\cM_{a,\Theta}(\xi)$ for all $\xi\in\R$, where
{\footnotesize
	\be \label{cond:oep:tf}
	\cM_{a,\Theta}(\xi):=
	\left[ \begin{matrix}
		\wh{\Theta}(\xi) & &\\
		&\ddots &\\
		& &\wh{\Theta}(\xi+\tfrac{2\pi (\dm-1)}{\dm})\end{matrix}\right]
	-\left[ \begin{matrix}
		\ol{\wh{a}(\xi)}^\tp\\
		\vdots\\
		 \ol{\wh{a}(\xi+\tfrac{2\pi(\dm-1)}{\dm})}^\tp
	\end{matrix}\right] \wh{\Theta}(\dm \xi) \left[ \begin{matrix} \wh{a}(\xi), &\ldots, &\wh{a}(\xi+\tfrac{2\pi (\dm-1)}{\dm})\end{matrix}\right].
	\ee
}
For $r>1$, it is much harder to find $\Theta$ which makes $\cM_{a,\Theta}$ positive semi-definite at every $\xi\in\R$. Moreover, $\Theta$ has to satisfy additional complicated conditions to increase vanishing moments of $\psi$. The necessary condition $\cM_{a,\Theta}(\xi)\ge 0$ often fails even with the trivial choice $\wh{\Theta}(\xi):=I_r$ for many matrix-valued filters $a$. For example, consider the widely used Hermite cubic splines:
	 \be\label{HCS}\phi_1(x)=\begin{cases}
		(1-x)^2(1+2x), &x\in[0,1]\\
		(1+x)^2(1-2x), &x\in[-1,0)\\
		0, & \text{otherwise},
	\end{cases}\quad \phi_2(x)=\begin{cases}
		(1-x)^2x, &x\in[0,1]\\
		(1+x)^2x, &x\in[-1,0)\\
		0, & \text{otherwise}.
	\end{cases}
	\ee
	Then $\phi=[\phi_1,\phi_2]^{\tp}$ is a $2$-refinable vector of compactly supported $\Lp{2}$ functions satisfying $\wh{\phi}(2\xi)=\wh{a}(\xi)\wh{\phi}(\xi)$ for all $\xi\in\R$, where $a\in\lrs{0}{2}{2}$ is the Hermite interpolatory filter:
	 \be\label{mask:HCS:0}\wh{a}(\xi)=\frac{1}{16}\begin{bmatrix}
		4e^{i\xi}+8+4e^{-i\xi} & 6(e^{i\xi}-e^{-i\xi})\\
		-(e^{i\xi}-e^{-i\xi}) &-e^{i\xi}+4-e^{-i\xi}
	\end{bmatrix},\quad\xi\in\R.\ee
	With the trivial choice $\wh{\Theta}(\xi)=I_2$, we have
	 $$\det(I_2-\ol{\wh{a}(\xi)}^{\tp}\wh{a}(\xi))=\frac{1}{1024}(\cos(\xi)-1)(\cos^3(\xi)+5\cos^2(\xi)+331\cos(\xi)-113)<0,$$
	whenever $\xi\in(-\frac{\pi}{4},\frac{\pi}{4})$. Therefore, the necessary condition $\cM_{a,I_2}(\xi)\ge 0$ cannot hold.

Therefore, finding a suitable choice of $\Theta$ is the first difficulty in the construction of OEP-based tight multiframelets. To our best knowledge, OEP-based tight multiframelets was investigated in \cite{mothesis,mo06} and
OEP-based dual multiframelets have been studied and constructed in \cite{hm03,han09}.

Another difficulty is the deconvolution in Step 3 of the reconstruction process. The conditions in Theorem~\ref{thm:df} guarantee that the original signal $v_0$ must be a solution of the deconvolution problem $\tilde{v}_0=v_0*\Theta$. However, deconvolving $\tilde{v}_0=v_0*\Theta$ is inefficient and not stable.
The essence of OEP is to replace the original $\phi$ by another desired $\dm$-refinable (vector) function $\eta$ satisfying $\wh{\eta}(\dm\xi)=\wh{\mra}(\xi)\wh{\eta}(\xi)$ with     $\wh{\mra}(\xi):=\wh{\theta}(\dm\xi)\wh{a}(\xi)\wh{\theta}(\xi)^{-1}$ (recall that $\wh{\Theta}=\ol{\wh{\theta}}^\tp\wh{\tilde{\theta}}$).
The above obstacle is due to the fact that the refinement mask/filter of $\mra$ often has infinite support, even though $\eta$ has compact support.
In the multiframelet case ($r>1$), the determinant of $\wh{\Theta}$ could even be identically zero and therefore, the solution of the deconvolution problem is not unique at all. Hence, it appears impossible for multiframelets constructed through OEP to achieve both high vanishing moments and an efficient framelet transform simultaneously.
The first breakthrough to knock down this dead-end for OEP is probably \cite{han09} showing the real advantage of OEP for $r>1$. If $\Theta$ is strongly invertible, then the solution $\mathring{v}_0$ to the deconvolution problem $\tilde{v}_0=\mathring{v}_0*\Theta$ is simply given by $\mathring{v}_0=\tilde{v}_0*\Theta^{-1}$ (here $\Theta^{-1}\in\lrs{0}{r}{r}$ is the filter with $\wh{\Theta^{-1}}=\wh{\Theta}^{-1}$) and the trouble of deconvolution is completely gone. Indeed, as proved in \cite[Theorem~1.2]{han09}, if $r>1$, then one can always construct a dual $\dm$-framelet through OEP in Theorem~\ref{thm:df}
from any pair of matrix-valued filters such that the dual framelet has the highest possible vanishing moments and both $\theta$ and $\tilde{\theta}$ are strongly invertible (consequently, $\Theta$ is strongly invertible).
It is the purpose of this paper to
show in Theorem~\ref{thm:qtf} that we can always construct quasi-tight multiframelets, which are much stronger than simply dual multiframelets, with all desired properties being kept.

\section{Proof of Theorem~\ref{thm:qtf}}
\label{sec:qtf}
In this section, we prove our main result Theorem~\ref{thm:qtf}. For the convenience of later presentation, we need the following notations:

\begin{enumerate}
	
	\item[(1)] For $\gamma\in \Z$ and $u\in (\sq)^{s\times r}$, the \emph{$\gamma$-coset sequence} of $u$ with respect to the dilation factor $\dm$ is the sequence $u^{[\gamma;\dm]}\in (\sq)^{s\times r}$ given by
	 $$u^{[\gamma;\dm]}(k)=u(\gamma+\dm k),\quad k\in\Z.$$
	It is straightforward to check that
	 \be\label{coset}\wh{u}(\xi)=\sum_{\gamma=0}^{\dm-1}\wh{u^{[\gamma;\dm]}}(\dm\xi)e^{-i\gamma\xi},\quad\forall u\in\lrs{0}{s}{r},\quad \xi\in\R.\ee
Thus by letting $\FF_{r;\dm}(\xi)$ to be the $(\dm r)\times (\dm r)$ matrix defined via
\be\label{Fourier}\FF_{r;\dm}(\xi):=\left(e^{-i(l-1)(\xi+2\pi\frac{k-1}{\dm})}I_r\right)_{1\le l,k\le \dm},\ee	
we have
\be\label{coset:1}\left[\wh{u}(\xi),\wh{u}(\xi+\tfrac{2\pi}{\dm}),\dots,\wh{u}(\xi+\tfrac{2\pi (\dm-1)}{\dm})\right]=\left[\wh{u^{[0;\dm]}}(\dm\xi),\wh{u^{[1;\dm]}}(\dm\xi),\dots,\wh{u^{[\dm-1;\dm]}}(\dm\xi)\right]\FF_{r;\dm}(\xi).
\ee
Observe that $\ol{\FF_{r;\dm}(\xi)}^{\tp}\FF_{r;\dm}(\xi)=\dm I_{\dm r}$, \er{coset:1} is equivalent to
\be\label{coset:2}\left[\wh{u}(\xi),\wh{u}(\xi+\tfrac{2\pi}{\dm}),\dots,\wh{u}(\xi+\tfrac{2\pi (\dm-1)}{\dm})\right]\ol{\FF_{r;\dm}(\xi)}^{\tp}=\dm\left[\wh{u^{[0;\dm]}}(\dm\xi),\wh{u^{[1;\dm]}}(\dm\xi),\dots,\wh{u^{[\dm-1;\dm]}}(\dm\xi)\right].
\ee

	\item[(2)] For $j\in\{1,\dots,\dm\}$ and $u\in\lrs{0}{r}{r}$, let $D_{u;\dm}(\xi)$ be the $(\dm r)\times (\dm r)$ block diagonal matrix defined via
	 \be\label{Duni}D_{u;\dm}(\xi):=\begin{bmatrix}\wh{u}(\xi) & & &\\
		& \wh{u}(\xi+\tfrac{2\pi}{\dm}) & & \\
	&	&\ddots &\\
	&	& &\wh{u}(\xi+\tfrac{2\pi (\dm-1)}{\dm})\end{bmatrix},\ee
 and let $E_{u;\dm}(\xi)$ be the $(\dm r)\times(\dm r)$ block matrix, whose $(l,k)$-th $r\times r$ block is
	 \be\label{Euni}(E_{u;\dm}(\xi))_{l,k}:=\wh{u^{[k-l;\dm]}}(\xi),\ee
	for $1\le l,k\le \dm$. Then direct calculation yields
	 \be\label{DEF}\FF_{r;\dm}(\xi)D_{u;\dm}(\xi)\ol{\FF_{r;\dm}(\xi)}^{\tp}=\dm E_{u;\dm}(\dm\xi).\ee

\end{enumerate}

The following theorem plays a key role in the proof of Theorem~\ref{thm:qtf}.

\begin{theorem}\label{thm:qtf:nf}Let $\dm\ge 2$ and $r\ge2$ be positive integers and let $a\in\lrs{0}{r}{r}$. Suppose that $a$ has $m$ sum rules with respect to $\dm$ with a matching filter $\vgu\in\lrs{0}{1}{r}$ satisfying $\wh{\vgu}(\xi)=[1,0,\dots,0]+\bo(|\xi|^m)$ as $\xi \to 0$. Further suppose that $\phi$ is an $r\times 1$ vector of compactly supported functions in $\Lp{2}$ satisfying $\wh{\phi}(\dm\xi)=\wh{a}(\xi)\wh{\phi}(\xi)$ and $\wh{\phi}(\xi)=[1,0,\dots,0]^{\tp}+\bo(|\xi|^n)$ as $\xi\to 0$ for some $n\ge 2m$. Then for any strongly invertible $U\in\lrs{0}{r}{r}$ satisfying \er{eq:ortho}, there exist $b\in \lrs{0}{s}{r}$ and $\epsilon_1,\dots\epsilon_s\in\{\pm1\}$ for some $s\in\N$ such that
	
\begin{enumerate}
\item[(i)]$\{a;b\}_{\UU;(\eps_1,\dots,\eps_s)}$ is an OEP-based quasi-tight $\dm$-multiframelet filter bank, i.e.,
\be\label{nf:pr:1}\ol{\wh{a}(\xi)}^{\tp}\wh{\UU}(\dm\xi)\wh{a}(\xi)+\ol{\wh{b}(\xi)}^{\tp}\DG(\epsilon_1,\dots,\epsilon_s) \wh{b}(\xi)=\wh{\UU}(\xi),\ee
\be\label{nf:pr:2}\ol{\wh{a}(\xi)}^{\tp}\wh{\UU}(\dm\xi)\wh{a}(\xi+2\pi\tfrac{\gamma}{\dm})+\ol{\wh{b}(\xi)}^{\tp}\DG(\epsilon_1,\dots,\epsilon_s) \wh{b}(\xi+2\pi\tfrac{\gamma}{\dm})=0,\quad\gamma=1,\dots,\dm-1,\ee
 where $\wh{\UU}(\xi)=\ol{\wh{U}(\xi)}^{-\tp}\wh{U}(\xi)^{-1}$ for all $\xi\in\R$.

\item[(ii)] $\{\eta;\psi\}_{(\eps_1,\dots,\eps_s)}$ is a compactly supported quasi-tight $\dm$-framelet in $\Lp{2}$ such that all the entries of $\psi$ have $m$ vanishing moments, where
\be\label{qtf0}\wh{\eta}(\xi)=\wh{U}(\xi)^{-1}\wh{\phi}(\xi),\quad \wh{\psi}(\xi)=\wh{b}(\xi/\dm)\wh{\phi}(\xi/\dm),\quad\xi\in\R.\ee
\end{enumerate}
\end{theorem}

\bp By our assumptions on $\wh{\vgu}$ and $\wh{\phi}$, we see that $\wh{a}$ must take the form in \er{normalform}, i.e.,
$$\wh{a}(\xi)=\begin{bmatrix} \wh{a_{1,1}}(\xi) & \wh{a_{1,2}}(\xi)\\
\wh{a_{2,1}}(\xi)& P_{2,2}(\xi)\end{bmatrix},$$
with
\begin{align*}
&\wh{a_{1,1}}(\xi)=(1+e^{-i\xi}+\dots+e^{-i(\dm-1)\xi})^mP_{1,1}(\xi)=1+\bo(|\xi|^n),\quad \xi \to 0,\\
&\wh{a_{1,2}}(\xi)=(1-e^{-i\dm\xi})^mP_{1,2}(\xi),\quad \wh{a_{2,1}}(\xi)=(1-e^{-i\xi})^nP_{2,1}(\xi),
\end{align*}
where $P_{1,1}, P_{1,2}, P_{2,1}$ and $P_{2,2}$ are some $1\times 1$, $1\times (r-1)$, $(r-1)\times 1$ and $(r-1)\times (r-1)$ matrices of $2\pi$-periodic trigonometric polynomials.
Define
$$\wh{a_1}(\xi):=\wh{\UU}(\xi)-\ol{\wh{a}(\xi)}^\tp\wh{\UU}(\dm\xi)\wh{a}(\xi),$$
$$\wh{a_j}(\xi):=-\ol{\wh{a}(\xi)}^\tp\wh{\UU}(\dm\xi)\wh{a}(\xi+\tfrac{2\pi (j-1)}{\dm}),\quad j=2,\dots,\dm.$$

For $j=1$, using \er{eq:ortho} and the fact that $\|\wh{\phi}(\xi)\|^2=1+\bo(|\xi|^n)$ as $\xi\to 0$, we have
$$
\begin{bmatrix}p_1(\xi) & p_2(\xi) \\
p_3(\xi)& p_4(\xi)\end{bmatrix}:=
\wh{a_1}(\xi)=\begin{bmatrix}1 & \\
& \wh{C}(\xi)\end{bmatrix}-\ol{\wh{a}(\xi)}^\tp\begin{bmatrix}1 & \\
& \wh{C}(\dm \xi)\end{bmatrix}\wh{a}(\xi)+\bo(|\xi|^n),\quad\xi\to 0,
$$
where $\wh{C}(\xi)=\DG\left(\|\wh{u_2}(\xi)\|^2,\dots,\|\wh{u_r}(\xi)\|^2\right)$ and $\wh{u_j}$ denote the $j$-th column of $\wh{U}^{-1}$ for all $j=1,\dots,r$. Then
\begin{itemize}
	\item $p_1$ is a $2\pi$-periodic trigonometric polynomial satisfying
$$
p_1(\xi)=1-\left(|\underbrace{\wh{a_{1,1}}(\xi)}_{=1+\bo(|\xi|^n)}|^2+\underbrace{\ol{\wh{a_{2,1}}(\xi)}^{\tp}}_{=\bo(|\xi|^n)}\wh{C}(\dm\xi)\wh{a_{2,1}}(\xi)\right)+\bo(|\xi|^n)=\bo(|\xi|^n),\quad\xi\to 0.
$$
	
	\item $p_2$ is a $1\times (r-1)$ vector of $2\pi$-periodic trigonometric polynomials satisfying
$$
p_2(\xi)=-\ol{\wh{a_{1,1}}(\xi)}\underbrace{\wh{a_{1,2}}(\xi)}_{=\bo(|\xi|^m)}-\underbrace{\ol{\wh{a_{2,1}}(\xi)}^{\tp}}_{=\bo(|\xi|^n)}\wh{C}(\dm\xi)P_{2,2}(\xi)+\bo(|\xi|^n)=\bo(|\xi|^m),\quad\xi\to 0.
$$
	
\item  $p_3$ and $p_4$ are $(r-1)\times 1$ and $(r-1)\times (r-1)$ matrices of $2\pi$-periodic trigonometric polynomials satisfying $p_3(\xi)=\ol{p_2(\xi)}^{\tp}=\bo(|\xi|^m),\quad\xi\to 0$.
\end{itemize}
Since $n\ge 2m$, we conclude that
$\wh{a_1}$ admits the following factorization:
\be\label{fac:mom:2}\wh{a_1}(\xi)=\ol{\begin{bmatrix}(1-e^{-i\xi})^m &\\
		& I_{r-1}\end{bmatrix}}^{\tp}\wh{B_1}(\xi)
\begin{bmatrix}(1-e^{-i\xi})^m &\\
	& I_{r-1}\end{bmatrix},\ee
for some $B_1\in\lrs{0}{r}{r}$.

For $j=2,\dots,\dm$, we have
$$
\begin{bmatrix}
p_{j,1}(\xi) & p_{j,2}(\xi)\\
p_{j,3}(\xi) & p_{j,4}(\xi)
\end{bmatrix}:=\wh{a_j}(\xi)=-\ol{\wh{a}(\xi)}^{\tp}\begin{bmatrix}1 &\\
& \wh{C}(\dm\xi)\end{bmatrix}\wh{a}(\xi+2\pi\tfrac{ (j-1)}{\dm})+\bo(|\dm\xi|^n),\quad \xi \to 0.
$$
By $n\ge 2m$, we observe that
\begin{itemize}
	\item $p_{j,1}$ is a $2\pi$-periodic trigonometric polynomial satisfying
	 $$\begin{aligned}p_{j,1}(\xi)&=-\left(\underbrace{\ol{\wh{a_{1,1}}(\xi)}}_{=\bo(|\xi+2\pi (j-1)/\dm|^m)}
\underbrace{\wh{a_{1,1}}(\xi+2\pi \tfrac{(j-1)}{\dm})}_{=\bo(|\xi|^m)}+
\underbrace{\ol{\wh{a_{2,1}}(\xi)}^{\tp}}_{=\bo(|\xi|^n)}\wh{C}(\dm\xi)\underbrace{\wh{a_{2,1}}(\xi+2\pi\tfrac{ (j-1)}{\dm})}_{=\bo(|\xi+2\pi (j-1)/\dm|^n)}\right)+\bo(|\dm\xi|^n)\\	 &=(1-e^{i\xi})^m(1-e^{-i(\xi+2\pi\frac{(j-1)}{\dm})})^m\wh{F_j}(\xi),
	\end{aligned}$$
	where $\wh{F_j}(\xi)$ denotes some $2\pi$-periodic trigonometric polynomial.	
	
	\item $p_{j,2}$ is a $1\times (r-1)$ vector of $2\pi$-periodic trigonometric polynomials satisfying
$$
\begin{aligned}&p_{j,2}(\xi)=-\ol{\wh{a_{1,1}}(\xi)}\underbrace{\wh{a_{1,2}}(\xi+2\pi\tfrac{ (j-1)}{\dm})}_{=\bo(|\xi|^m)}-\underbrace{\ol{\wh{a_{2,1}}(\xi)}^{\tp}}_{=\bo(|\xi|^n)}\wh{C}(\dm\xi)P_{2,2}(\xi+2\pi\tfrac{ (j-1)}{\dm})+\bo(|\dm\xi|^n)\\
	&=\bo(|\xi|^m),\quad\xi\to 0.
	\end{aligned}$$
	
	\item $p_{j,3}$ is a $(r-1)\times 1$ vector of $2\pi$-periodic trigonometric polynomials satisfying
$$
\begin{aligned}p_{j,3}(\xi)&=-\underbrace{\ol{\wh{a_{1,2}}(\xi)}}
_{=\bo(|\xi+2\pi(j-1)/\dm|^m)}\wh{a_{1,1}}(\xi+2\pi\tfrac{ (j-1)}{\dm})-\ol{P_{2,2}(\xi)}^{\tp}\wh{C}(\dm\xi)\underbrace{\wh{a_{2,1}}(\xi+2\pi\tfrac{ (j-1)}{\dm})}_{=\bo(|\xi+2\pi(j-1)/\dm|^m)}+\bo(|\dm\xi|^n)\\ &=\bo(|\xi+2\pi\tfrac{(j-1)}{\dm}|^m),\quad\xi\to 0.
	\end{aligned}
$$

\item $p_{j,4}$ is some $(r-1)\times(r-1)$ matrix of $2\pi$-periodic trigonometric polynomials.
\end{itemize}
Thus $\wh{a_j}$ admits the following factorization:
\be\label{fac:mom:3}\wh{a_j}(\xi)=\ol{\begin{bmatrix}(1-e^{-i\xi})^m &\\
		& I_{r-1}\end{bmatrix}}^{\tp}\wh{B_j}(\xi)
\begin{bmatrix}(1-e^{-i(\xi+2\pi\tfrac{(j-1)}{\dm})})^m &\\
	& I_{r-1}\end{bmatrix},\quad j=2,\dots,\dm,\ee
for some $B_j\in\lrs{0}{r}{r}$. Hence by letting
\be\label{Deltam}\wh{\Delta_m}(\xi):=\begin{bmatrix}(1-e^{-i\xi})^m & 0\\
	0 & I_{r-1}\end{bmatrix},\ee
and
\be \label{cond:oep:tf:2}
\begin{aligned}\cM_{a,\UU}(\xi):=&
\left[ \begin{matrix}
	\wh{\UU}(\xi) & &\\
	&\ddots &\\
	& &	\wh{\UU}(\xi+\tfrac{2\pi (\dm-1)}{\dm})\end{matrix}\right]-\left[ \begin{matrix}
	\ol{\wh{a}(\xi)}^\tp\\
	\vdots\\
	 \ol{\wh{a}(\xi+\tfrac{2\pi(\dm-1)}{\dm})}^\tp
\end{matrix}\right] 	\wh{\UU}(\dm \xi) \left[ \begin{matrix} \wh{a}(\xi), &\ldots, &\wh{a}(\xi+\tfrac{2\pi (\dm-1)}{\dm})\end{matrix}\right],
\end{aligned}\ee
it follows from \er{fac:mom:2} and \er{fac:mom:3} that
\be
\label{fac:maW}\cM_{a,\UU}(\xi)=\ol{D_{\Delta_m;\dm}(\xi)}^{\tp}\cM(\xi)D_{\Delta_m;\dm}(\xi),
\ee
where $\cM$ is some $(\dm r)\times (\dm r)$ Hermitian matrix of $2\pi$-periodic trigonometric polynomials, and $D_{\Delta_m;\dm}$ is defined via \er{Duni} with $u=\Delta_m$.

It follows from \er{DEF} and \er{fac:maW} that
\be\label{pr:coset}\begin{aligned}
	 &\dm^{-2}\FF_{r;\dm}(\xi)\cM_{a,\UU}(\xi)\ol{\FF_{r;\dm}(\xi)}^{\tp}\\
	 =&\dm^{-4}\left(\FF_{r;\dm}(\xi)\ol{D_{\Delta_m;\dm}(\xi)}^{\tp}\ol{\FF_{r;\dm}(\xi)}^{\tp}\right)\left(\FF_{r;\dm}(\xi)\cM(\xi)\ol{\FF_{r;\dm}(\xi)}^{\tp}\right)\left(\FF_{r;\dm}(\xi)D_{\Delta_m;\dm}(\xi)\ol{\FF_{r;\dm}(\xi)}^{\tp}\right)\\
	 =&\ol{E_{\Delta_m;\dm}(\dm\xi)}^{\tp}\tilde{\cM}(\xi)E_{\Delta_m;\dm}(\dm\xi),
\end{aligned}\ee
where $\tilde{\cM}(\xi)=\dm^{-2}\FF_{r;\dm}(\xi)\cM(\xi)\ol{\FF_{r;\dm}(\xi)}^{\tp}$ and  $E_{\Delta_m;\dm}$ is defined as in \er{Euni} with $u=\Delta_m$.

On the other hand, using\er{coset:2} and \er{DEF}, we see that
\be\label{pr:coset:1}
\dm^{-2}\FF_{r;\dm}(\xi)\cM_{a,\UU}(\xi)\ol{\FF_{r;\dm}(\xi)}^{\tp}=\dm^{-1}E_{\UU;\dm}(\dm\xi)-\left[ \begin{matrix}
	 \ol{\wh{a^{[0;\dm]}}(\dm\xi)}^\tp\\
	\vdots\\
		 \ol{\wh{a^{[\dm-1;\dm]}}(\dm\xi)}^\tp
\end{matrix}\right] 	\wh{\UU}(\dm \xi) \left[ \begin{matrix} \wh{a^{[0;\dm]}}(\dm\xi), &\ldots, &\wh{a^{[\dm-1;\dm]}}(\dm\xi)\end{matrix}\right].
\ee
Hence $\tilde{\cM}(\xi)$ only depends on $\dm\xi$, say $\tilde{\cM}(\xi)=\mathring{\cM}(\dm\xi),$ where $\mathring{\cM}$ is some $(\dm r) \times (\dm r)$ Hermitian matrix of $2\pi$-periodic trigonometric polynomials.
We now claim that $\mathring{\cM}$ can be factorized in the following way:
\be\label{Herm:fac}\mathring{\cM}(\xi)=\ol{\tilde{U}(\xi)}^{\tp}\DG(I_{s_1},-I_{s_2})\tilde{U}(\xi),\ee
for some $(s_1+s_2)\times (\dm r)$ matrix of $2\pi$-periodic trigonometric polynomials $\tilde{U}(\xi)$. In fact, there always exist $(\dm r)\times (\dm r)$ matrices of $2\pi$-periodic trigonometric polynomials $\mathring{\cM}_1(\xi)$ and $\mathring{\cM}_2(\xi)$ such that
$$\mathring{\cM}(\xi)=\ol{\mathring{\cM}_1(\xi)}^{\tp}\mathring{\cM}_1(\xi)-\ol{\mathring{\cM}_2(\xi)}^{\tp}\mathring{\cM}_2(\xi).$$
For example, take $\mathring{\cM}_1(\xi)=I_{\dm r}+\frac{1}{4}\mathring{\cM}(\xi)$ and $\mathring{\cM}_2(\xi)=I_{\dm r}-\frac{1}{4}\mathring{\cM}(\xi).$ Then simply choose $\tilde{U}=[\mathring{\cM}_1^{\tp},\mathring{\cM}_2^{\tp}]^{\tp}$, we see that \er{Herm:fac} holds with $s_1=s_2=\dm r$.
Once we have factorized $\mathring{\cM}$ as in \er{Herm:fac}, define  $b\in\lrs{0}{(s_1+s_2)}{r}$ and $\epsilon_1,\dots,\epsilon_{s_1+s_2}\in\{\pm1\}$ via
\be\label{bpb}\wh{b}(\xi):=\tilde{U}(\dm\xi)\FF_{r;\dm}(\xi)\begin{bmatrix}
	I_r\\
	\pmb{0}_{(\dm-1)r \times r}
\end{bmatrix}\wh{\Delta_m}(\xi),
\ee
\be\label{epsilon}\epsilon_1=\dots=\epsilon_{s_1}=1,\quad \epsilon_{s_1+1}=\dots=\epsilon_{s_1+s_2}=-1,\ee

where $\pmb{0}_{q\times t}$ denotes the $q\times t$ zero matrix. Using \er{pr:coset}, \er{bpb} and \er{epsilon}, we have
\be\label{pr:oep}\begin{aligned}&\begin{bmatrix}
		\ol{\wh{b}(\xi)}^{\tp}\\
		\vdots\\
		 \ol{\wh{b}(\xi+2\pi\tfrac{\dm-1}{\dm})}^{\tp}
	 \end{bmatrix}\DG(\epsilon_1,\dots,\epsilon_{s_1+s_2})\begin{bmatrix}\wh{b}(\xi)&\dots&\wh{b}(\xi+2\pi\tfrac{\dm-1}{\dm})\end{bmatrix}\\
	 =&\ol{D_{\Delta_m;\dm}(\xi)}^{\tp}\ol{\FF_{r;\dm}(\xi)}^{\tp}\ol{\tilde{U}(\dm\xi)}^{\tp}\DG(I_{s_1},-I_{s_2})\tilde{U}(\dm\xi)\FF_{r;\dm}(\xi)D_{\Delta_m;\dm}(\xi)\\
	 =&\ol{D_{\Delta_m;\dm}(\xi)}^{\tp}\ol{\FF_{r;\dm}(\xi)}^{\tp}\mathring{\cM}(\dm\xi)\FF_{r;\dm}(\xi)D_{\Delta_m;\dm}(\xi)\\
	 =&\ol{\FF_{r;\dm}(\xi)}^{\tp}\ol{E_{\Delta_m;\dm}(\dm\xi)}^{\tp}\mathring{\cM}(\dm\xi)E_{\Delta_m;\dm}(\dm\xi)\FF_{r;\dm}(\xi)\\
	=&\cM_{a,\UU}(\xi).
\end{aligned}\ee
Note that \er{pr:oep} is equivalent to \er{nf:pr:1} and \er{nf:pr:2}. This proves item (i).

Define $\eta$ and $\psi$ as in \er{qtf0}.
By $\wh{\phi}(\xi)=[1,0,\dots,0]^{\tp}+\bo(|\xi|^n)$ as $\xi\to 0$ and $n\ge 2m$, we have
$$\begin{aligned}\wh{\psi}(\dm\xi)&=\wh{b}(\xi)\wh{\phi}(\xi)=\tilde{U}(\dm\xi)\FF_{r;\dm}(\xi)\begin{bmatrix}
I_r\\
\pmb{0}_{(\dm-1)r \times r}
\end{bmatrix}\underbrace{\wh{\Delta_m}(\xi)\wh{\phi}(\xi)}_{=\bo(|\xi|^m)}=\bo(|\xi|^m),\quad\xi\to 0,\end{aligned}
$$
which means that all the entries of $\psi$ have $m$ vanishing moments.
Note that $\ol{\wh{\phi}(0)}^\tp \wh{\UU}(0)\wh{\phi}(0)=1$.
Now by Theorem~\ref{thm:df}, $\{\eta;\psi\}_{(\eps_1,\dots,\eps_s)}$ is a quasi-tight $\dm$-framelet in $\Lp{2}$. This proves item (ii).
\ep

Now we are ready to prove Theorem~\ref{thm:qtf}.

\begin{proof}[\textbf{Proof of Theorem~\ref{thm:qtf}}]
Let $n\ge 2m$ be a positive integer. By Theorem~\ref{thm:normalform}, there exists a strongly invertible filter $\theta\in\lrs{0}{r}{r}$ such that
\be\label{nf1}\wh{\mathring{\vgu}}(\xi):=\wh{\vgu}(\xi)\wh{\theta}(\xi)^{-1}=r^{-\frac{1}{2}}\wh{\Vgu}(\xi)+\bo(|\xi|^m),\quad \wh{\mathring{\phi}}(\xi):=\wh{\theta}(\xi)\wh{\phi}(\xi)=r^{-\frac{1}{2}}\ol{\wh{\Vgu}(\xi)}^{\tp}+\bo(|\xi|^n),\quad\xi\to 0,\ee
where $\wh{\Vgu}(\xi):=[1, e^{i\xi/r},\ldots, e^{i (r-1)\xi/r}]$ (also see \eqref{vgu:special}).
Now by Theorem~\ref{thm:normalform:2}, there exists a strongly invertible $U\in\lrs{0}{r}{r}$ such that
 \be\label{W:ortho:0}\begin{aligned}\ol{\wh{U}(\xi)}^{-\tp}\wh{U}(\xi)^{-1}&=\DG\left(\|\wh{\mrphi}(\xi)\|^2,\|\wh{u_2}\|^2,\dots,\|\wh{u_r}\|^2\right)
 +\bo(|\xi|^n)\\
 	 &=\DG\left(1,\|\wh{u_2}\|^2,\dots,\|\wh{u_r}\|^2\right)
 +\bo(|\xi|^n),\quad\xi\to 0,
 \end{aligned}\ee
 where $\wh{u_j}$ denotes the $j$-th column of $\wh{U}^{-1}$, and
\be\label{Wvgu}\wh{\bpv}(\xi):=\wh{\mathring{\vgu}}(\xi)\wh{U}(\xi)^{-1}=[1,0,\dots,0]+\bo(|\xi|^{m}),\quad \wh{\bpphi}(\xi):=\wh{U}(\xi)\wh{\mathring{\phi}}(\xi)=[1,0,\dots,0]^{\tp}+\bo(|\xi|^{n}),\quad \xi\to 0.\ee
Define $\mra,\bpa\in\lrs{0}{r}{r}$ via
\be\label{mra}\wh{\mra}(\xi)=\wh{\theta}(\dm\xi)\wh{a}(\xi)\wh{\theta}(\xi)^{-1},\quad \wh{\bpa}(\xi)=\wh{U}(\dm\xi)\wh{\mra}(\xi)\wh{U}(\xi)^{-1}.\ee
It is trivial that $\bpa$ has $m$ sum rules with the matching filter $\bpv$ and $\wh{\bpphi}(\dm\xi)=\wh{\bpa}(\xi)\wh{\bpphi}(\xi)$ holds for all $\xi\in\R$. Thus by item (i) of Theorem~\ref{thm:qtf:nf}, there exist $\bpb\in \lrs{0}{s}{r}$ and $\epsilon_1,\dots,\epsilon_s\in\{\pm1\}$ such that \er{nf:pr:1} and \er{nf:pr:2} hold with $a$ and $b$ being replaced by $\bpa$ and $\bpb$ respectively. Hence by defining $b,\mrb\in \lrs{0}{s}{r}$ via
\be\label{highpass}\wh{b}(\xi)=\wh{\bpb}(\xi)\wh{U}(\xi)\wh{\theta}(\xi),\quad \wh{\mrb}(\xi)=\wh{b}(\xi)\wh{\theta}(\xi)^{-1},\ee
we see that item (2) follows right away.\\

Next, define
$$\wh{\psi}(\xi):=\wh{b}(\xi/\dm)\wh{\phi}(\xi/\dm)=\wh{\mrb}(\xi/\dm)\wh{\mrphi}(\xi/\dm)=\wh{\bpb}(\xi/\dm)\wh{\bpphi}(\xi/\dm),$$
for all $\xi\in\R$. By item (ii) of Theorem~\ref{thm:qtf:nf}, $\{\mrphi;\psi\}_{(\eps_1,\dots,\eps_s)}$ is a quasi-tight $\dm$-framelet in $\Lp{2}$ with $\wh{\psi}(\xi)=\bo(|\xi|^m)$ as $\xi\to 0$. This proves item (1). Moreover, \er{mphi:mpsi} holds trivially.\\

Finally, we have
$$\wh{\Vgu}(\xi)\ol{\wh{\mrb}(\xi)}^{\tp}=\wh{\Vgu}(\xi)\ol{\wh{\theta}(\xi)}^{-\tp}\ol{\wh{b}(\xi)}^{\tp}=r^{\frac{1}{2}}\ol{\wh{\phi}(\xi)}^{\tp}\ol{\wh{b}(\xi)}^{\tp}+\bo(|\xi|^n)=r^{\frac{1}{2}}\ol{\wh{\psi}(\dm\xi)}^{\tp}+\bo(|\xi|^n)=\bo(|\xi|^m),$$
and
$$\begin{aligned}\wh{\Vgu}(\xi)\ol{\wh{\mra}(\xi)}^{\tp}&=\wh{\Vgu}(\xi)\ol{\wh{\theta}(\xi)}^{-\tp}\ol{\wh{a}(\xi)}^{\tp}\ol{\wh{\theta}(\dm\xi)}^{\tp}=r^{\frac{1}{2}}\ol{\wh{\phi}(\xi)}^{\tp}\ol{\wh{a}(\xi)}^{\tp}\ol{\wh{\theta}(\dm\xi)}^{\tp}+\bo(|\xi|^n)\\
&=r^{\frac{1}{2}}\ol{\wh{\phi}(\dm\xi)}^{\tp}\ol{\wh{\theta}(\dm\xi)}^{\tp}+\bo(|\xi|^n)=\wh{\Vgu}(\dm\xi)+\bo(|\xi|^m),
\end{aligned}$$
as $\xi \to 0$. Hence by Theorem~\ref{prop:bp}, items (3) and (4) must hold. The proof is now complete.
\end{proof}

\section{Characterization of the Filters $\theta$ in Theorem~\ref{thm:qtf}}
\label{sec:momfilter}

To construct a quasi-tight framelet in Theorem~\ref{thm:qtf} from a matrix-valued filter having multiplicity greater than one, a desired filter $\theta$  in Theorem~\ref{thm:qtf} plays a key role and is guaranteed to exist by Theorem~\ref{thm:qtf}.
In this section, we characterize all desired filters $\theta$ in Theorem~\ref{thm:qtf}. Such a characterization will allow us to construct all possible quasi-tight framelet filter banks and quasi-tight framelets in Theorem~\ref{thm:qtf} having all the desired properties.

For $m\in \N$, we define a sequence $\nabla^m\pmb{\delta}\in \lp{0}$ through
$\wh{\nabla^m\pmb{\delta}}(\xi):=(1-e^{-i\xi})^m$.
Before proceeding further, we need the following technical lemma, which provides an equivalent way of interpreting the balanced vanishing moments condition.

\begin{lemma}\label{vm:bo:0}Let $r\ge 2$ and $s\in\N$ be positive integers. For any $m\in \N$, a filter $b\in\lrs{0}{s}{r}$ has $m$ balanced vanishing moments (i.e., \eqref{cond:bvmo} holds) if and only if
\be\label{coset:bbq} \wh{b}(\xi)=\left[\wh{q^{[0;r]}}(\xi),\wh{q^{[1;r]}}(\xi),\dots,\wh{q^{[r-1;r]}}(\xi)\right]
E_{\nabla^m \pmb{\delta};r}(\xi),\quad\xi\in\R,\ee
	for some $q\in\lrs{0}{s}{1}$,
	where $E_{\nabla^m\pmb{\delta};r}$ is defined in \er{Euni} with $\dm=r$ and
$u=\nabla^m\pmb{\delta}$.
\end{lemma}

\bp Suppose that $b$ has $m$ balanced vanishing moments, i.e., \er{cond:bvmo} holds. We deduce that
\be\label{vm:b:r}\wh{b_1}(r\xi)+e^{-i\xi}\wh{b_2}(r\xi)+\dots+e^{-i\xi(r-1)}\wh{b_r}(r\xi)=\wh{\nabla^m\pmb{\delta}}(\xi)\wh{q}(\xi),\ee
for some $q\in\lrs{0}{s}{1}$, where $\wh{b_j}$ denotes the $j$-th column of $\wh{b}$.  Since
$\wh{u}(\xi)=\sum_{\gamma=0}^{\dm-1} \wh{u^{[\gamma;\dm]}}(\dm \xi) e^{-i\gamma \xi}$ in
\er{coset}, we have
\be\label{vm:b:2}\begin{aligned}	 &\wh{b_1}(r\xi)+e^{-i\xi}\wh{b_2}(r\xi)+\dots+e^{-i(r-1)\xi}\wh{b_r}(r\xi)=\left(\sum_{j=0}^{r-1}e^{-ij\xi}\wh{\nabla^m\pmb{\delta}^{[j;r]}}(r\xi)\right)\left(\sum_{k=0}^{r-1}e^{-ik\xi}\wh{q^{[k;r]}}(r\xi)\right)\\ =&\sum_{j=0}^{r-1}e^{-ij\xi}\sum_{k=0}^j\wh{\nabla^m\pmb{\delta}^{[j-k;r]}}(r\xi)\wh{q^{[k;r]}}(r\xi)+\sum_{j=0}^{r-2}e^{-ij\xi}e^{-ir\xi}\sum_{k=j+1}^{r-1}\wh{\nabla^m\pmb{\delta}^{[j+r-k;r]}}(r\xi)\wh{q^{[k;r]}}(r\xi).
\end{aligned}\ee
Hence
\be\label{br}\wh{b_{r}}(\xi)=\sum_{k=0}^{r-1}\wh{\nabla^m\pmb{\delta}^{[r-1-k;r]}}(\xi)\wh{q^{[k;r]}}(\xi),\ee
and
\be\label{bj1}\begin{aligned}
	 \wh{b_{j+1}}(\xi)&=\sum_{k=0}^j\wh{\nabla^m\pmb{\delta}^{[j-k;r]}}(\xi)\wh{q^{[k;r]}}(\xi)+e^{-i\xi}\sum_{k=j+1}^{r-1}\wh{\nabla^m\pmb{\delta}^{[j+r-k;r]}}(\xi)\wh{q^{[k;r]}}(\xi)\\
	 &=\sum_{k=0}^{r-1}\wh{\nabla^m\pmb{\delta}^{[j-k;r]}}(\xi)\wh{q^{[k;r]}}(\xi),\quad j=0,\dots,r-2,\end{aligned}\ee
for all $\xi\in\R$, where the last line of \er{bj1} follows from the fact that
\be\label{coset:neg:b}e^{-i\xi}\wh{u^{[r-l;r]}}(\xi)=\wh{u^{[-l;r]}}(\xi),\quad\forall u\in\lrs{0}{t}{r},\quad l\in\Z.\ee
Thus \er{coset:bbq} follows right away from \er{br} and \er{bj1}.

Conversely, suppose that \er{coset:bbq} holds. Then \er{br} and \er{bj1} must hold. Thus we deduce that \er{vm:b:2} and \er{vm:b:r} hold. Now \er{cond:bvmo} follows trivially.
\ep

The following result provides a characterization for all the desired filters $\theta$ in Theorem~\ref{thm:qtf}.

\begin{theorem}\label{thm:momfilter}
Let $\dm\ge2$ and $r\ge 2$ be integers and $a\in\lrs{0}{r}{r}$ be a finitely supported matrix-valued filter. Let $\phi\in (\Lp{2})^r$ be a compactly supported vector function satisfying $\wh{\phi}(\dm\xi)=\wh{a}(\xi)\wh{\phi}(\xi)$. Assume that $a$ has $m$ sum rules with respect to $\dm$ with a matching filter $\vgu\in\lrs{0}{1}{r}$ such that $\wh{\vgu}(0)\wh{\phi}(0)=1$. Define $\wh{\Vgu}(\xi):=[1, e^{i\xi/r},\ldots, e^{i (r-1)\xi/r}]$ as in \er{vgu:special}.
Let $\theta\in\lrs{0}{r}{r}$ be a strongly invertible filter. If the filter $\theta$ satisfies the following two conditions:
	\begin{enumerate}
		\item[(i)]    There exist $2\pi$-periodic trigonometric polynomials $\wh{c}$ and $\wh{d}$ with $|\wh{c}(0)|=|\wh{d}(0)|=\|\wh{\Vgu}(0)\|^{-1}=r^{-\frac{1}{2}}$ such that
\be\label{vgu:bo}
\wh{\mathring{\vgu}}(\xi):=\wh{\vgu}(\xi)\wh{\theta}(\xi)^{-1}
=\wh{c}(\xi)\wh{\Vgu}(\xi)+\bo(|\xi|^m),
\quad \xi\to 0,
\ee
\be\label{phi:bo}
\wh{\mrphi}(\xi):=\wh{\theta}(\xi)\wh{\phi}(\xi)
=\wh{d}(\xi)\ol{\wh{\Vgu}(\xi)}^\tp+\bo(|\xi|^m),
\quad \xi \to 0;
\ee

\item[(ii)] All the entries of the following two matrices are $2\pi$-periodic trigonometric polynomials:
		 \be\label{m1}M_0(\xi):=\ol{E_{\nabla^m\pmb{\delta};r}(\xi)}^{-\tp}\left(I_r-\ol{\wh{\mra}(\xi)}^{\tp}\wh{\mra}(\xi)\right)E_{\nabla^m\pmb{\delta};r}(\xi)^{-1},\ee
		 \be\label{m2}M_j(\xi):=-\ol{E_{\nabla^m\pmb{\delta};r}(\xi)}^{-\tp}\ol{\wh{\mra}(\xi)}^{\tp}\wh{\mra}(\xi+\tfrac{2\pi j}{\dm})E_{\nabla^m\pmb{\delta};r}(\xi+\tfrac{2\pi j}{\dm})^{-1},\quad j=1,\dots,\dm-1,\ee
		where $\wh{\mra}(\xi)=\wh{\theta}(2\xi)\wh{a}(\xi)\wh{\theta}(\xi)^{-1}$
		and $E_{\nabla^m\pmb{\delta};r}$ is defined in \er{Euni} with $\dm=r$ and
$u=\nabla^m\pmb{\delta}$,
\end{enumerate}
then there must exist $b\in \lrs{0}{s}{r}$ and $\eps_1,\ldots,\eps_s\in \{\pm 1\}$ such that all the items (1)--(4) of Theorem~\ref{thm:qtf} are satisfied.
Conversely, if all the items (1)--(4) of Theorem~\ref{thm:qtf} are satisfied for some $b\in \lrs{0}{s}{r}$ and $\eps_1,\ldots,\eps_s\in \{\pm 1\}$, then the filter $\theta$ must satisfy item (ii) above and if additionally
\be\label{eigen:a}\text{$1$ is a simple eigenvalue of $\wh{a}(0)$ and $\det(\dm^j I_r-\wh{a}(0))\ne 0$ for all $j\in\Z\setminus\{0\}$ with $j\ge 1-\dm$,}
\ee
then $\theta$ must also satisfy item (i) above.
\end{theorem}

\bp For simplicity of presentation, we define $\wh{E_{m;r}}(\xi):=E_{\nabla^m\pmb{\delta};r}(\xi)$. Then $E_{m;r}$ is a finitely supported matrix-valued filter.
First observe from \er{DEF} and the fact $\FF_{r;r}(\xi)\ol{\FF_{r;r}}^{\tp}=rI_{r^2}$ that
$$\det(\wh{E_{m;r}}(\xi))=\det(E_{\nabla^m\pmb{\delta};r}(\xi))=\det(D_{\nabla^m\pmb{\delta};r}(\xi/r))=\prod_{j=0}^{r-1}(1-e^{-i(\xi+2\pi j)/r})^m=(1-e^{-i\xi})^m,\quad\xi\in\R.$$
Therefore, $\wh{E_{m;r}}(\xi)$ is invertible for all $\xi\in\R\setminus 2\pi\Z$, and thus all
the matrices $M_0,M_1,\dots,M_{\dm-1}$ in item (ii) are well defined for all $\xi\in\R\setminus 2\pi\Z$.

Suppose that items (i) and (ii) hold. Define
$$\cM_{\mra}(\xi):=I_{\dm r}-\begin{bmatrix}\ol{\wh{\mra}(\xi)}^{\tp}\\
\vdots\\
\ol{\wh{\mra}(\xi+\tfrac{2\pi (\dm-1)}{\dm})}^{\tp}\end{bmatrix}\left[\wh{\mra}(\xi),\dots, \wh{\mra}(\xi+\tfrac{2\pi (\dm-1)}{\dm})\right].$$
By item (ii), $\cM_{\mra}$ admits the following factorization:
\be\label{fac:cm:emr}\cM_{\mra}(\xi)=\ol{D_{E_{m;r};\dm}(\xi)}^{\tp}M(\xi)D_{E_{m;r};\dm}(\xi),\quad\xi\in\R,\ee
where $M(\xi)$ is some $(\dm r)\times (\dm r)$ matrix of $2\pi$-periodic trigonometric polynomials. Applying the same argument as in the proof of Theorem~\ref{thm:qtf:nf}, we have
$$\frac{1}{\dm^2}\FF_{r;\dm}(\xi)\cM_{\mra}(\xi)\ol{\FF_{r;\dm}(\xi)}^{\tp}=\ol{E_{E_{m;r};\dm}(\dm\xi)}^{\tp}\tilde{M}(\dm\xi)E_{E_{m;r};\dm}(\dm\xi)^{\tp},$$
where $\tilde{M}(\xi)$ is some $\dm r\times \dm r$ Hermitian matrix of $2\pi$-periodic trigonometric polynomials. Thus there exists an $s\times r$ matrix $\tilde{U}(\xi)$ of $2\pi$-periodic trigonometric polynomials such that
$$\tilde{M}(\xi)=\ol{\tilde{U}(\xi)}^{\tp}\DG(I_{s_1},-I_{s_2})\tilde{U}(\xi),\quad\xi\in\R,$$
for some $s_1,s_2\in\N_0$ satisfying $s_1+s_2=s$. Define $\mrb,b\in\lrs{0}{s}{r}$ and $\epsilon_1,\dots,\epsilon_s\in\{\pm1\}$ via
\be\label{mrb:emr}
\wh{\mrb}(\xi):=\tilde{U}(\dm\xi)\FF_{r;\dm}(\xi)\begin{bmatrix}I_r\\
\pmb{0}_{r(\dm-1)\times r}\end{bmatrix}\wh{E_{m;r}}(\xi),\quad \wh{b}(\xi):=\wh{\mrb}(\xi)\wh{\theta}(\xi),\ee
\be\label{eps:emr}\epsilon_1=\dots=\epsilon_{s_1}=1,\quad \epsilon_{s_1+1}=\dots=\epsilon_{s}=-1.\ee
Use \er{fac:cm:emr}, \er{mrb:emr} and \er{eps:emr}, it is straightforward to check that item (2) of Theorem~\ref{thm:qtf} holds. Next, by letting $q\in\lrs{0}{s}{1}$ be such that
$$\left[\wh{q^{[0;r]}}(\xi),\wh{q^{[1;r]}}(\xi),\dots,\wh{q^{[r-1;r]}}(\xi)\right]=\tilde{U}(\dm\xi)\FF_{r;\dm}(\xi)\begin{bmatrix}I_r\\
\pmb{0}_{r(\dm-1)\times r}\end{bmatrix},$$
we see that items (3) and (4) of Theorem~\ref{thm:qtf} follow immediately from Lemma~\ref{vm:bo:0}. On the other hand, since item (i) holds, it follows that $\|\wh{\mrphi}(0)\|^2=1$ and thus Theorem~\ref{thm:df} yields that $\{\mrphi;\psi\}_{(\eps_1,\dots,\eps_s)}$ is a quasi-tight $\dm$-framelet in $\Lp{2}$, where $\wh{\psi}(\xi)=\wh{b}(\xi/\dm)\wh{\phi}(\xi/\dm)$. Moreover, \er{phi:bo} and item (3) of Theorem~\ref{thm:qtf} guarantee that $\psi$ has $m$ vanishing moments. This proves item (1) of Theorem~\ref{thm:qtf}. Hence all the claims of Theorem~\ref{thm:qtf} hold.

Conversely, suppose that $\theta$ is a strongly invertible filter and
all the claims in Theorem~\ref{thm:qtf} holds. By item (3) of Theorem~\ref{thm:qtf}, \er{cond:bvmo}holds. Thus by Lemma~\ref{vm:bo:0}, there exists $q\in\lrs{0}{s}{1}$ such that
\be\label{coset:bq}
\wh{\mrb}(\xi)=\left[\wh{q^{[0;r]}}(\xi),\wh{q^{[1;r]}}(\xi),\dots,\wh{q^{[r-1;r]}}(\xi)\right]\wh{E_{m;r}}(\xi),\quad\xi\in\R.\ee
By \er{qtffb=1} and \er{qtffb=0}, we have
\be\label{m0:1}I_r-\ol{\wh{\mra}(\xi)}^{\tp}\wh{\mra}(\xi)=\ol{\wh{E_{m;r}}(\xi)}^{\tp}\begin{bmatrix}\ol{\wh{q^{[0;r]}}(\xi)}^{\tp}\\
	\vdots\\
	 \ol{\wh{q^{[r-1;r]}}(\xi)}^{\tp}\end{bmatrix}\DG(\epsilon_1,\dots,\epsilon_s)\left[\wh{q^{[0;r]}}(\xi),\dots\wh{q^{[r-1;r]}}(\xi)\right]\wh{E_{m;r}}(\xi),\ee
and
\be\label{mj:1}\begin{aligned}&-\ol{\wh{\mra}(\xi)}^{\tp}\wh{\mra}(\xi+\tfrac{2\pi j}{\dm})\\
	 =&\ol{\wh{E_{m;r}}(\xi)}^{\tp}\begin{bmatrix}\ol{\wh{q^{[0;r]}}(\xi)}^{\tp}\\
		\vdots\\
		 \ol{\wh{q^{[r-1;r]}}(\xi)}^{\tp}\end{bmatrix}\DG(\epsilon_1,\dots,\epsilon_s)\left[\wh{q^{[0;r]}}(\xi+\tfrac{2\pi j}{\dm}),\dots\wh{q^{[r-1;r]}}(\xi+\tfrac{2\pi j}{\dm})\right]\wh{E_{m;r}}(\xi+\tfrac{2\pi j}{\dm}),\end{aligned}\ee
for all $\xi\in\R$ and $j=1,\dots,\dm-1$. Define $M_0$ as \er{m1} and $M_j$ as \er{m2} for $j=1,\dots,\dm-1$. By \er{m0:1} and \er{mj:1}, it is trivial that item (ii) holds. If in addition that \er{eigen:a} holds, by \cite[Theorem 4.4]{han09}, \er{vgu:bo} and \er{phi:bo} must hold for some $2\pi$-periodic trigonometric polynomials $\wh{c}$ and $\wh{d}$ with $\wh{c}(0)\neq 0$ and $\wh{d}(0)\neq 0$. Furthermore, by \cite[Theorem 4.1.10 or Theorem 6.4.1]{hanbook}, we have $\|\wh{\mrphi}(0)\|^2=\ol{\wh{\phi}(0)}^{\tp}\ol{\wh{\theta}(0)}^{\tp}\wh{\theta}(0)\wh{\phi}(0)=1$. This means $|\wh{d}(0)|=\|\wh{\Vgu}(0)\|^{-1}=r^{-\frac{1}{2}}$. Finally, noting that $\wh{\mathring{\vgu}}(0)\wh{\mrphi}(0)=\wh{\vgu}(0)\wh{\phi}(0)=1$, we conclude that $|\wh{c}(0)|=r^{-\frac{1}{2}}$. This proves item (i).
\ep

The following corollary is an immediate consequence of Theorem~\ref{thm:momfilter}.

\begin{cor}Let $\dm\ge2$ and $r\ge 2$ be integers and $a\in\lrs{0}{r}{r}$ be a finitely supported matrix-valued filter. Let $\phi\in (\Lp{2})^r$ be a compactly supported vector function satisfying $\wh{\phi}(\dm\xi)=\wh{a}(\xi)\wh{\phi}(\xi)$ and define $\wh{\Vgu}(\xi):=[1, e^{i\xi/r},\ldots, e^{i (r-1)\xi/r}]$ as in \er{vgu:special}. Assume that $a$ has $m$ sum rules with respect to $\dm$ with a matching filter $\vgu\in\lrs{0}{1}{r}$ such that $\wh{\vgu}(0)\wh{\phi}(0)=1$. If $\theta\in\lrs{0}{r}{r}$ is a strongly invertible filter such that item (ii) of Theorem~\ref{thm:momfilter} holds and
\be\label{mrphi:vgu}\wh{\mrphi}(\xi):=\wh{\theta}(\xi)\wh{\phi}(\xi)=\wh{c}(\xi)\ol{\wh{\Vgu}(\xi)}^{\tp}+\bo(|\xi|^m),\quad \xi\to 0,\ee
	for some $2\pi$-periodic trigonometric polynomial $\wh{c}$ with $\wh{c}(0)\neq 0$, then
\be\label{phi:moment}\|\wh{\mrphi}(\xi)\|^2=
\|\wh{\mrphi}(0)\|^2
+\bo(|\xi|^{2m}),\quad\xi\to 0.\ee
\end{cor}

\bp Let $\theta\in\lrs{0}{r}{r}$ be a strongly invertible filter such that all above assumptions are satisfied. From the proof of Theorem~\ref{thm:momfilter}, we deduce that there exist $b\in\lrs{0}{s}{r}$ and $\epsilon_1,\dots,\epsilon_s\in\{\pm1\}$ such that item (2) of Theorem~\ref{thm:qtf} and the following condition hold:
\be\label{ref:mr:mom}\wh{\mrphi}(\dm\xi)=\wh{\mra}(\xi)\wh{\mrphi}(\xi)\quad\mbox{and}\quad \wh{\psi}(\xi):=\wh{\mrb}(\xi/\dm)\wh{\mrphi}(\xi/\dm)=\bo(|\xi|^m)\quad\mbox{as}\quad\xi\to 0,\ee
where $\wh{\mrphi}:=\wh{\theta}\wh{\phi}$, $\wh{\mra}:=\wh{\theta}(\dm\cdot)\wh{a}\wh{\theta}^{-1}$ and $\wh{\mrb}:=\wh{b}\wh{\theta}^{-1}$.
By multiplying $\ol{\wh{\mrphi}(\xi)}^{\tp}$ to the left and $\wh{\mrphi}(\xi)$ to the right to both sides of \er{qtffb=1} and using \er{ref:mr:mom}, we have
\be\label{PR:phi}\|\wh{\mrphi}(\dm\xi)\|^2+\ol{\wh{\psi}(\dm\xi)}^{\tp}\DG(\epsilon_1,\dots,\epsilon_s)\wh{\psi}(\dm\xi)=\|\wh{\mrphi}(\xi)\|^2,\ee
which yields
\be\label{PR:phi:2}\|\wh{\mrphi}(\dm\xi)\|^2= \|\wh{\mrphi}(\xi)\|^2+\bo(|\xi|^{2m}),\quad\xi\to 0.\ee
Hence \er{phi:moment} follows from \er{PR:phi:2} and $\dm\ge 2$.
\ep

The condition in \eqref{phi:moment} a key for vanishing moments of derived framelets $\psi$ from $\mathring{\phi}$.
To construct examples of quasi-tight framelets with all desired properties in Theorem~\ref{thm:qtf}, we first construct a desired filter $\theta$ satisfying items (i) and (ii) of Theorem~\ref{thm:momfilter}. Then Theorem~\ref{thm:momfilter} guarantees the existence of a filter $b\in \lrs{0}{s}{r}$ and $\eps_1,\ldots,\eps_s\in \{\pm 1\}$ satisfying all the items (1)--(4) of Theorem~\ref{thm:qtf}.

\section{Some Examples of Spline Quasi-tight Framelets with High Balancing Orders}
\label{sec:example}

In this section, we present some examples to illustrate our main result Theorem~\ref{thm:qtf}. Based on Theorem~\ref{thm:momfilter}, we provide some guidelines for constructing quasi-tight framelets satisfying all desired properties in Theorem~\ref{thm:qtf}.

Let $\dm\ge2$ and $r\ge 2$ be integers and $a\in\lrs{0}{r}{r}$ be a finitely supported matrix-valued filter. Let $\phi\in (\Lp{2})^r$ be a compactly supported refinable vector function satisfying $\wh{\phi}(\dm\xi)=\wh{a}(\xi)\wh{\phi}(\xi)$. Assume that $a$ has $m$ sum rules with respect to $\dm$ with a matching filter $\vgu\in\lrs{0}{1}{r}$ such that $\wh{\vgu}(0)\wh{\phi}(0)=1$. Let $\wh{\Vgu}(\xi):=[1, e^{i\xi/r},\ldots, e^{i (r-1)\xi/r}]$ as in \er{vgu:special}. The general construction steps are as follows:

\begin{enumerate}
	\item[(1)] Construct a strongly invertible filter $\theta\in \lrs{0}{r}{r}$ with short support satisfying items (i) and (ii) of Theorem~\ref{thm:momfilter}.
	
	\item[(2)] Construct a filter $b\in\lrs{0}{s}{r}$ such that \er{qtffb=1} and \er{qtffb=0} are satisfied (where $\mra$ and $\mrb$ are given by \er{mphi:mpsi}) for some $\epsilon_1,\dots,\epsilon_s\in\{\pm1\}$, and $\wh{b}(\xi)\wh{\phi}(\xi)=\bo(|\xi|^m)$ as $\xi\to 0$.
The existence of such $b\in\lrs{0}{s}{r}$ and $\epsilon_1,\dots,\epsilon_s\in\{\pm1\}$
is guaranteed by Theorem~\ref{thm:momfilter}.
\end{enumerate}
Define $\wh{\psi}(\xi)=\wh{b}(\xi/\dm)\wh{\phi}(\xi/\dm)$. Then $\{\mrphi;\psi\}_{(\epsilon_1,\dots,\epsilon_s)}$ is a compactly supported quasi-tight $\dm$-framelet in $\Lp{2}$ satisfying all the desired properties in Theorem~\ref{thm:qtf}.

\begin{figure}[htb]
		\centering
		 \begin{subfigure}[b]{0.3\textwidth}
			\centering
			 \includegraphics[width=\textwidth]{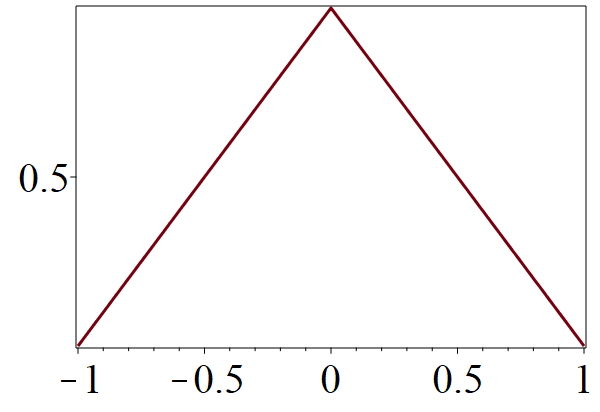}
			\caption{$\phi$}
			\label{fig:PhiOrg:B20}
		\end{subfigure}
		 \begin{subfigure}[b]{0.3\textwidth}
			\centering
			 \includegraphics[width=\textwidth]{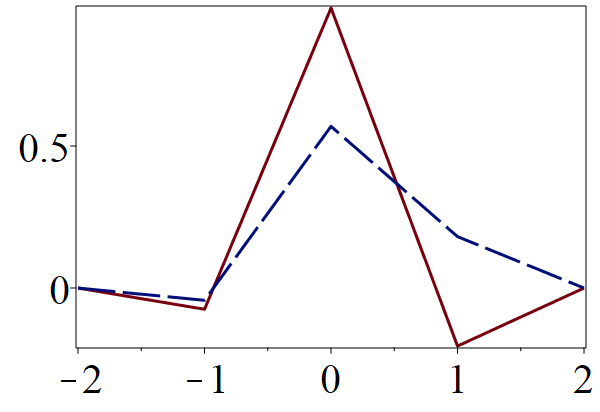}
			\caption{$\mrphi$}
			\label{fig:PhiNew:B20}
		\end{subfigure}
		 \begin{subfigure}[b]{0.3\textwidth}
			\centering
			 \includegraphics[width=\textwidth]{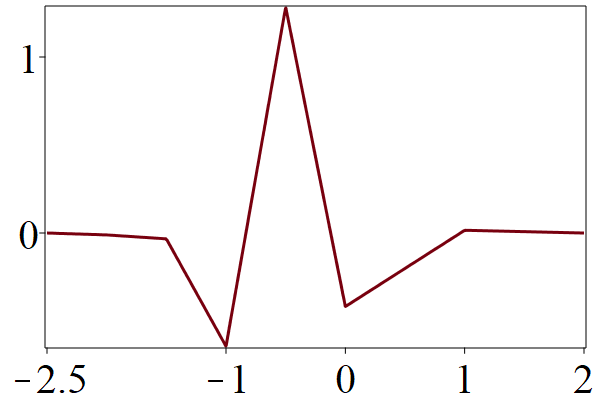}
			\caption{$\psi^1$}
			\label{fig:Psi1:B20}
		\end{subfigure}
		 \begin{subfigure}[b]{0.3\textwidth}
			\centering
			 \includegraphics[width=\textwidth]{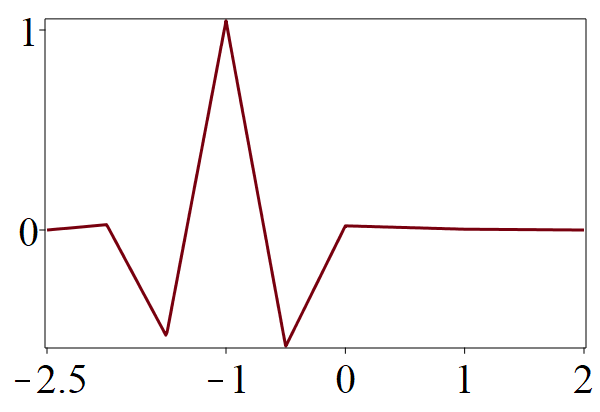}
			\caption{$\psi^2$}
			\label{fig:Psi2:B20}
		\end{subfigure}
		 \begin{subfigure}[b]{0.3\textwidth}
			\centering
			 \includegraphics[width=\textwidth]{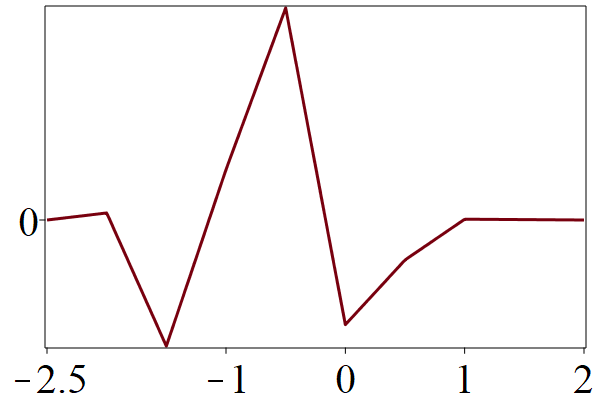}
			\caption{$\psi^3$}
			\label{fig:Psi3:B20}
		\end{subfigure}
		 \begin{subfigure}[b]{0.3\textwidth}
			\centering
			 \includegraphics[width=\textwidth]{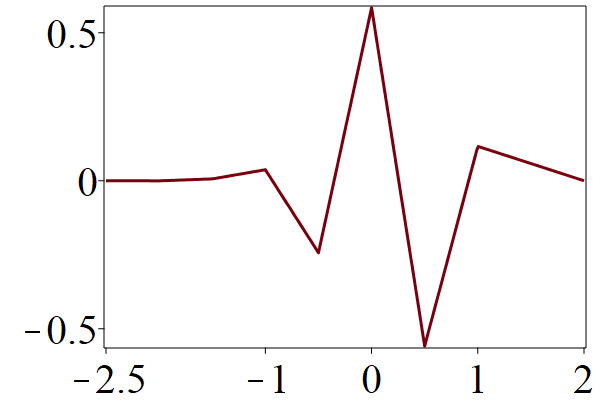}
			\caption{$\psi^4$}
			\label{fig:Psi4:B20}
		\end{subfigure}
		
		\caption{Graphs of $\phi=[B_2(\cdot-1),0]^{\tp}$ and a new refinable vector function $\mrphi$, together with graphs of $\psi^1,,\dots,\psi^4$ constructed from $\phi$ in Example~\ref{B20}. A graph with a solid (resp. dash) line denotes the first (resp. second) component of a vector function.
$\{\mathring{\phi}; (\psi^1,\ldots,\psi^4)^\tp\}_{(1,1,1,-1)}$ is a compactly supported quasi-tight $2$-framelet in $\Lp{2}$ with balanced vanishing moments $2$.
}
		\label{fig:psi:B20}
	\end{figure}

\begin{exmp}\label{B20}
Let $\dm=r=2$, and consider $\phi=[B_2(\cdot-1),0]^{\tp}$, where $B_2$ is the B-spline of order $2$ in \er{bspline}. Then $\phi$ satisfies $\wh{\phi}(2\xi)=\wh{a}(\xi)\wh{\phi}(\xi)$ with a filter $a\in\lrs{0}{2}{2}$ being given by
$$
\wh{a}(\xi)=\begin{bmatrix}\wh{A}(\xi) & 0\\
	0 &p(\xi)\end{bmatrix},
$$
with
	 \be\label{B2A}\wh{A}(\xi)=\frac{1}{4}(e^{i\xi}+2+e^{-i\xi}),\ee
	and $p(\xi)$ is any $2\pi$-periodic trigonometric polynomial. Note that $\sr(a,2)=2$ with any matching filter $\vgu\in\lrs{0}{1}{2}$ satisfying $\wh{\vgu}(\xi)=[1,0]+\bo(|\xi|^2)$ as $\xi\to0$.
We obtain a strongly invertible filter $\theta\in \lrs{0}{2}{2}$ satisfying items (i) and (ii) of Theorem~\ref{thm:momfilter} as follows:
$$
\wh{\theta}(\xi)=\frac{1}{6\sqrt{2}(1+\sqrt{3})}\begin{bmatrix}
	 -(3+\sqrt{3})e^{-i\xi}+8\sqrt{3}+9-\sqrt{3}e^{i\xi} & (7+13\sqrt{3})e^{-i\xi}-11+(4+3\sqrt{3})e^{i\xi}\\[0.5cm]
	 (3\sqrt{3}-1)e^{-i\xi}+3\sqrt{3}+8-e^{i\xi} & -\left(\tfrac{23\sqrt{3}}{3}+9\right)e^{-i\xi}-\tfrac{11\sqrt{3}}{3}+\left(\frac{4\sqrt{3}}{3}+3\right)e^{i\xi}\end{bmatrix}.
$$
Direct computation shows that \er{vgu:bo} and \er{phi:bo} hold with $m=2$ and	 $$\wh{c}(\xi)=-\frac{\sqrt{2}}{2}+\frac{\sqrt{2}(\sqrt{3}-1)i}{8}\xi+\bo(|\xi|^2),\quad \wh{d}(\xi)=-\frac{\sqrt{2}}{2}-\frac{\sqrt{2}(\sqrt{3}-1)i}{8}\xi+\bo(|\xi|^2),\quad\xi\to 0.$$


Here we have the freedom to choose $p(\xi)$ such that the degree of $\wh{\mra}(\xi):=\wh{\theta}(2\xi)\wh{a}(\xi)\wh{\theta}(\xi)^{-1}$ is as small as possible for simple presentation. By choosing
$p(\xi)=\frac{1}{4}e^{-i\xi}+\frac{7}{2}+\frac{1}{4}e^{i\xi}$,
we obtain $b\in \lrs{0}{4}{2}$
	such that $\{\mra;\mrb\}_{(\epsilon_1,\epsilon_2,\epsilon_3, \epsilon_4)}$,
with
\[
\wh{\mathring{a}}(\xi):=
\wh{\theta}(\dm\xi)\wh{a}(\xi)\wh{\theta}(\xi)^{-1}
\quad
\mbox{and}\quad
\wh{\mathring{b}}(\xi):=\wh{b}(\xi)\wh{\theta}(\xi)^{-1}
\]
and $\epsilon_1=\epsilon_2=\epsilon_3=1$ and $\epsilon_4=-1$,
is a finitely supported quasi-tight $2$-multiframelet filter bank with $2$ balancing orders, where
$$
\wh{b}(\xi):=e^{i\xi}D\left(U_1+\sqrt{3}U_2\right)\wh{U_3}(2\xi)\wh{F}(\xi)\wh{\theta}(\xi),$$
with
	\begin{itemize}
		\item $D=\tfrac{1}{9+4\sqrt{3}}\DG\left(d_1\lambda_1,d_2\sqrt{46}\sqrt{\lambda_2+\lambda_3\sqrt{3}},d_3\sqrt{927889}\sqrt{\lambda_4+\lambda_5\sqrt{3}},d_4\sqrt{\lambda_6+\lambda_7\sqrt{3}}\right)$ where
		 $$d_1=-\frac{1}{25180085734704776},\quad d_2=\frac{1}{170664162417838019952956},$$
		 $$d_3=\frac{1}{281700714176366254998791127171400862653282671261438405782186152367677},$$
		 $$d_4=\frac{1}{5009062155049388954350661899051556951670668};$$
		 $$\lambda_1=\sqrt{5437900978131564+3120526414024498\sqrt{3}},$$
		 $$\lambda_2=3719615046635084853,\quad \lambda_3=2147522348686212558,$$
		 $$\lambda_4=2270305904207568012940508624742913001613984744660994673284621339421892193611324,$$
		 $$\lambda_5=1310761724937036116861992517021585500520429420925813977012891122860726054086295,$$
		 $$\lambda_6=378533294810068098618941042771044135712181,$$
		 $$\lambda_7=218546299655829430553984760745834819062292.$$
		
		\item $U_1$ and $U_2$ are the $4\times 4$ constant matrices given by
		 $$U_1=\small\begin{bmatrix}3392119873 &-2778324120 & -4025874315 &8672724540\\
		0&-8288658986045588 &-10141655898429575 &61508560391015634\\
		0& 0 & 0 & \lambda_8\\
		0 & 0 &\lambda_9 & \lambda_{10}\end{bmatrix},$$
		
		 $$U2=\small\begin{bmatrix}0 &670619101 & 457342305 &-11821258848\\
		0&481466912421912 &5873639680107160 &-35445004131731402\\
		0& 0 & 0 & \lambda_{11}\\
		0 & 0 &\lambda_{12} & \lambda_{13}\end{bmatrix},$$
		where
		 $$\lambda_8=655638898967488291661954282237504457186876,$$
		 $$\lambda_9=1173830888361736998172384986078,$$
		 $$\lambda_{10}=-1055722690591344263872331100,$$
		 $$\lambda_{11}=-378533294810068098618941042771044135712181$$
		 $$\lambda_{12}=-677711074380894544629813729391,$$
		 $$\lambda_{13}=609520590739704131515104420.$$
		
		\item $\wh{U_3}$ is the $4\times 4$ matrix of $2\pi$-periodic trigonometric polynomials given by $\wh{U_3}=\wh{U_{3,1}}\wh{U_{3,2}}$ where
		 $$\wh{U_{3,1}}(\xi)=\begin{bmatrix}
		1 & 0 & 0 & 0\\[0.4cm]
		0 & 1 & 0 & 0\\[0.4cm]
		 \left(\dfrac{782849\sqrt{3}}{61240674}+\dfrac{1521092}{10206779}\right)e^{-i\xi} &\dfrac{699997-272171\sqrt{3}}{10206779}e^{-i\xi}& 1 & 0\\[0.4cm]
		 \left(\dfrac{17348101\sqrt{3}}{122481348}-\dfrac{9216904}{30620337}\right)e^{-i\xi} & \left(\dfrac{4881421\sqrt{3}}{61240674}-\dfrac{5507049}{40827116}\right)e^{-i\xi} &0 & 1\end{bmatrix},$$
		 $$\wh{U_{3,2}}(\xi)=\begin{bmatrix}1 & 0 & 0 & 0\\[0.4cm]
		 \left(\dfrac{\sqrt{3}}{6}+\dfrac{1}{2}\right)e^{-i\xi} &  1 & -\dfrac{\sqrt{3}}{3} & \left(\dfrac{2\sqrt{3}}{3}-3\right)e^{i\xi}+\dfrac{5\sqrt{3}}{6}+\dfrac{3}{2}\\[0.4cm]
		0 & 0 & 0 & \dfrac{1}{2}\\[0.4cm]
		-\dfrac{1}{6}e^{-i\xi} & 0 & \dfrac{1}{3} & -\dfrac{5}{6}-\dfrac{2}{3}e^{i\xi}\end{bmatrix}.$$
		
		\item $\wh{F}(\xi)$ is the $4\times 2$ matrix of $2\pi$-periodic trigonometric polynomials given by
		 \be\label{Fbo}\wh{F}(\xi)=\begin{bmatrix}2 & -1-e^{i\xi}\\
			-1-e^{-i\xi} & 2\\
			2e^{-i\xi} & -1-e^{-i\xi}\\
			-e^{-2i\xi}-e^{-i\xi} & 2e^{-i\xi}\end{bmatrix}.\ee
	\end{itemize}
	
The filter $b$ is supported on $[-4,3]$. Define $\psi=[\psi^1,\psi^2,\psi^3,\psi^4]^{\tp}$ via $\wh{\psi}(\xi)=\wh{b}(\xi/2)\wh{\phi}(\xi/2)$.
Define a new refinable vector function $\wh{\mathring{\phi}}(\xi):=\wh{\theta}(\xi) \wh{\phi}(\xi)$.
Then $\|\wh{\mathring{\phi}}(\xi)\|^2=1+\bo(|\xi|^4)$ as $\xi\to 0$ and
$\{\mrphi;\psi\}_{(\epsilon_1,\epsilon_2,\epsilon_3, \epsilon_4)}$ is a compactly supported quasi-tight $2$-framelet in $\Lp{2}$ such that all the items (1)--(4) of Theorem~\ref{thm:qtf} are satisfied with $m=2$. Note that $\psi$ has $2$ vanishing moments.
See Figure~\ref{fig:psi:B20} for graphs of $\phi,\mathring{\phi}, \psi^1,\dots,\psi^4.$
\end{exmp}

\begin{figure}[htb]
		\centering
		 \begin{subfigure}[b]{0.3\textwidth}
			\centering
			 \includegraphics[width=\textwidth]{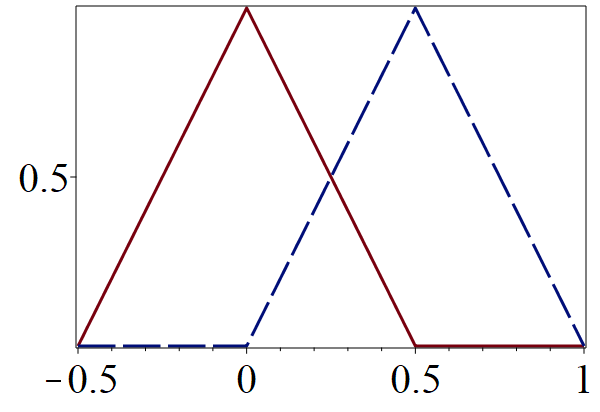}
			\caption{$\phi$}
			\label{fig:PhiOrg}
		\end{subfigure}
		 \begin{subfigure}[b]{0.3\textwidth}
			\centering
			 \includegraphics[width=\textwidth]{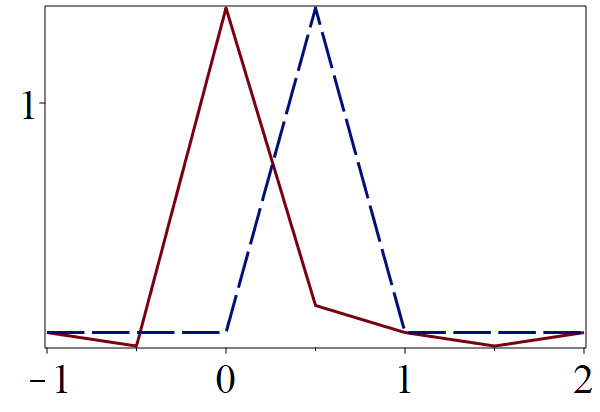}
			\caption{$\mrphi$}
			\label{fig:PhiNew}
		\end{subfigure}
		 \begin{subfigure}[b]{0.3\textwidth}
			\centering
			 \includegraphics[width=\textwidth]{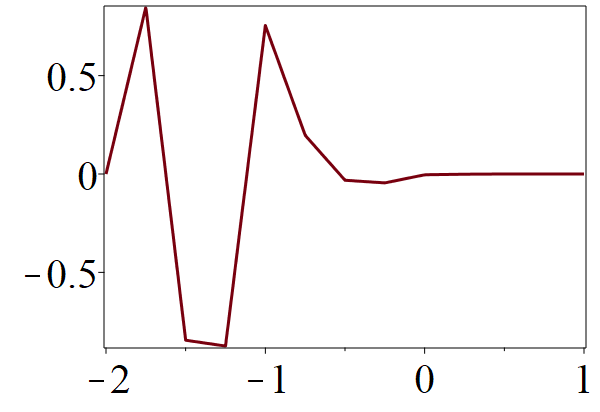}
			\caption{$\psi^1$}
			\label{fig:Psi1}
		\end{subfigure}
		 \begin{subfigure}[b]{0.3\textwidth}
			\centering
			 \includegraphics[width=\textwidth]{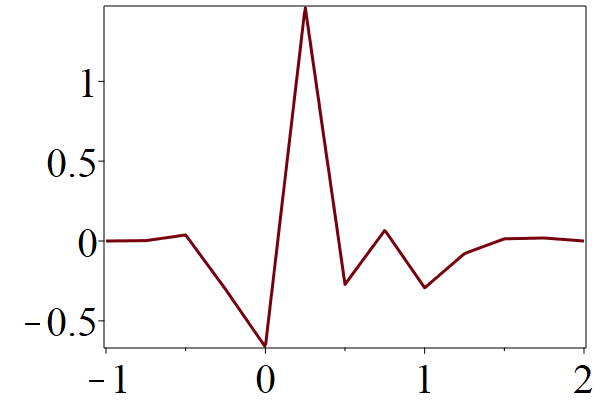}
			\caption{$\psi^2$}
			\label{fig:Psi2}
		\end{subfigure}
		 \begin{subfigure}[b]{0.3\textwidth}
			\centering
			 \includegraphics[width=\textwidth]{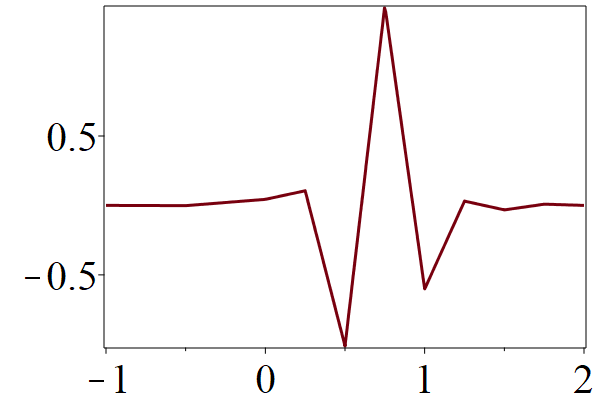}
			\caption{$\psi^3$}
			\label{fig:Psi3}
		\end{subfigure}
		 \begin{subfigure}[b]{0.3\textwidth}
			\centering
			 \includegraphics[width=\textwidth]{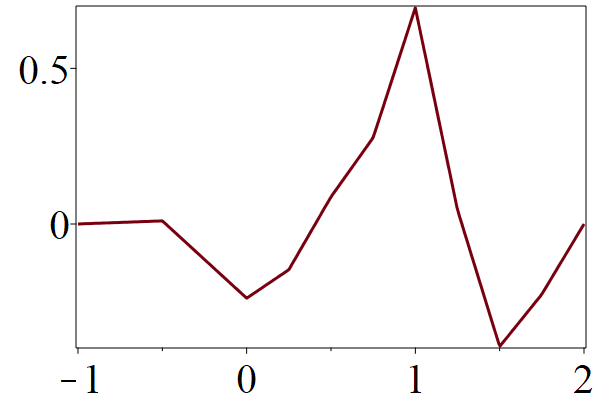}
			\caption{$\psi^4$}
			\label{fig:Psi4}
		\end{subfigure}
		\caption{Graphs of $\phi=[\varphi(2\cdot),\varphi(2\cdot-1)]^\tp$ and the new refinable vector function $\mrphi$, together with graphs of $\psi^1,,\dots,\psi^4$ constructed from $\phi$ in Example~\ref{B2B2}. A graph with a solid (resp. dash) line denotes the first (resp. second) component of a function vector.
$\{\mathring{\phi}; (\psi^1,\ldots,\psi^4)^\tp\}$ is a compactly supported tight $2$-framelet in $\Lp{2}$ with balanced vanishing moments $2$.
}
		\label{fig:psi}
	\end{figure}

\begin{exmp}\label{B2B2}Let $\dm=r=2$, and $\varphi:=B_2(\cdot-1)$ where $B_2$ is the B-spline of order $2$ in \er{bspline}. Then $\varphi$ satsfies $\wh{\varphi}(2\xi)=\wh{A}(\xi)\wh{\varphi}(\xi)$ where $\wh{A}$ is given by \er{B2A}. Note that $\sr(A,2)=2$.
Define $\phi:=[\varphi(2\cdot),\varphi(2\cdot-1)]^\tp$. By \cite[Proposition 6.2]{han09}, $\phi$ satisfies $\wh{\phi}(2\xi)=\wh{a}(\xi)\wh{\phi}(\xi)$, where
$$\wh{a}(\xi)=\frac{1}{4}\begin{bmatrix}2 & 1+e^{i\xi}\\
	2e^{-i\xi} &1+e^{-i\xi}\end{bmatrix},\quad\xi\in\R.$$
	Moreover, $\sr(a,2)=2$ with a matching filter $\vgu\in\lrs{0}{1}{2}$ satisfying	 $\wh{\vgu}(\xi)=[1,e^{i\xi/2}]+\bo(|\xi|^2)$ as $\xi\to 0$.
Now applying the general construction steps presented above, we obtain a desired strongly invertible filter $\theta\in \lrs{0}{2}{2}$ satisfying items (i) and (ii) of Theorem~\ref{thm:momfilter}:
$$
\wh{\theta}(\xi)=\frac{\sqrt{2}}{24}\begin{bmatrix}1 & -e^{i\xi}+2-e^{-i\xi}\\
	0 & 1\end{bmatrix}.
$$
Direct computation shows that \er{vgu:bo} and \er{phi:bo} hold with $m=2$ and
	 $$\wh{c}(\xi)=\wh{d}(\xi)=\frac{\sqrt{2}}{2}.$$
We obtain $b\in\lrs{0}{4}{2}$ such that $\{\mra;\mrb\}$, with $\mra$ and $\mrb$ being defined in \eqref{mab},
is a finitely supported tight $2$-framelet filter bank
with $2$ balancing orders.
For simplicity of presentation, we write $$
\wh{b}(\xi)=\wh{D}(2\xi)\wh{E}(2\xi)\wh{F}(\xi),
$$
with
	\begin{itemize}
		\item $\wh{D}(\xi)$ is the $4\times 4$ diagonal matrix of $2\pi$-periodic trigonometric polynomials given by $\wh{D}(\xi)=\DG\left(d_1(\xi),d_2,d_3,d_4\right)$ with
		 $$d_1(\xi)=\frac{\sqrt{70}}{43680\lambda_1}\left(-52\sqrt{5249}+364\sqrt{105}+\left(364\sqrt{105}+52\sqrt{5249}\right)e^{i\xi}\right),$$
		 $$d_2=-\frac{1}{96\lambda_1\lambda_2},\quad d_3=-\frac{\sqrt{6}}{576\lambda_3\lambda_4},\quad d_4=-\frac{\sqrt{2}}{48\lambda_3},$$
		where
		 $$\lambda_1=\sqrt{154299444795192054502909},\quad \lambda_2=\sqrt{392716620870},$$
		 $$\lambda_3=\sqrt{28366141},\quad \lambda_4=\sqrt{1870079147}.$$

		\item $\wh{E}(\xi)$ is the $4\times 4$ matrix of $2\pi$-periodic trigonometric polynomials given by $\wh{E}(\xi)=E_{-1}e^{i\xi}+E_0+E_1e^{-i\xi}$ with
		$$E_{-1}=\begin{bmatrix}0 & 0 & -116055812828 &-116541733616\\[0.5cm]
		0 & 0 &31811040936954791 & 514403328290664092\\[0.5cm]
		0 & 0 &-24580769 & -98323076\\[0.5cm]
		0 & 0 &449 & 1796\end{bmatrix},$$
		 $$E_0=\small\begin{bmatrix}2591577536 & 18434692032 & 6760552105 &0\\[0.5cm]
		-2573825532553116272 & -13507206835285209804 &-7105187917736625576 & -3908121799627104752\\[0.5cm]
		393292304 & 2064784596 &6589364832 & -16304613472\\[0.5cm]
		-7184 & -37716 &-81521 & -117450\end{bmatrix},$$
		$$E_{1}=\begin{bmatrix}0 & 0 & 0 &0\\[0.5cm]
		86906212475244992 & 618190750766215104 &226709010062624185 & 0\\[0.5cm]
		-205796480 & -398928912 & 239683385 & 0\\[0.5cm]
		-128302 & -76394 &-20406 & 0\end{bmatrix}.$$

		\item $\wh{F}$ is the $4\times 2$ matrix of $2\pi$-periodic trigonometric polynomials given by \er{Fbo}.
		
	\end{itemize}
	The filter $b$ is supported on $[-4,3]$, i.e., $b(k)=0$ whenever $k\notin \Z\cap[-4,3]$. Define $\psi=[\psi^1,\psi^2,\psi^3,\psi^4]^{\tp}$ via $\wh{\psi}(\xi)=\wh{b}(\xi/2)\wh{\phi}(\xi/2)$.
Define a new refinable vector function $\wh{\mathring{\phi}}(\xi):=\wh{\theta}(\xi) \wh{\phi}(\xi)$.
Then $\|\wh{\mathring{\phi}}(\xi)\|^2=1+\bo(|\xi|^4)$ as $\xi\to 0$ and
$\{\mrphi;\psi\}$ is a compactly supported tight $2$-framelet in $\Lp{2}$
such that all the desired properties in items (1)--(4) of Theorem~\ref{thm:qtf} are satisfied. Note that $\psi$ has $2$ vanishing moments. See Figure~\ref{fig:psi} for graphs of $\phi,\mathring{\phi},\psi^1,\psi^2,\psi^3,\psi^4$.
\end{exmp}

As we discussed in item (3) of Lemma~\ref{lem:conv},
a $2\pi$-periodic trigonometric polynomial $\wh{\Theta}$ is strongly invertible (i.e., $1/\wh{\Theta}$ is also a $2\pi$-periodic trigonometric polynomial) if and only if $\wh{\Theta}(\xi)=ce^{-im\xi}$ for some $m\in \Z$ and $c\in \C\bs\{0\}$.
Thus, as mentioned in Section~\ref{sec:ffrt}, for framelets constructed from scalar refinable functions to have high vanishing moments, usually it is inevitable to sacrifice the compactness of the associated discrete (scalar) framelet transform, because a non-trivial scalar filter $\Theta$ is not strongly invertible.
Examples~\ref{B20} and \ref{B2B2} demonstrate that this difficulty can be easily resolved by simply vectorizing the scalar refinable function, and do the constructions by using the new refinable vector function.

\begin{figure}[htb]
		\centering
		 \begin{subfigure}[b]{0.24\textwidth}
			\centering
			 \includegraphics[width=\textwidth]{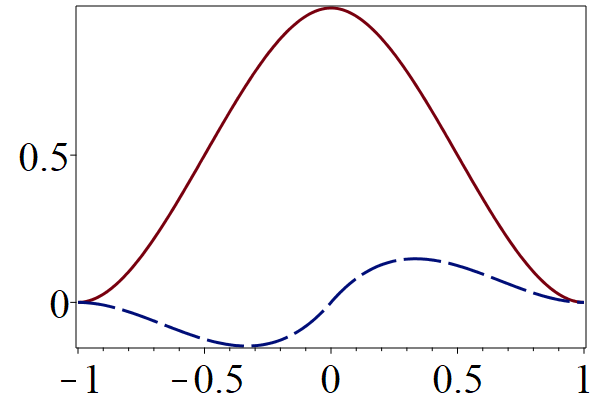}
			\caption{$\phi$}
			\label{fig:PhiOrg:HCS}
		\end{subfigure}
		 \begin{subfigure}[b]{0.24\textwidth}
			\centering
			 \includegraphics[width=\textwidth]{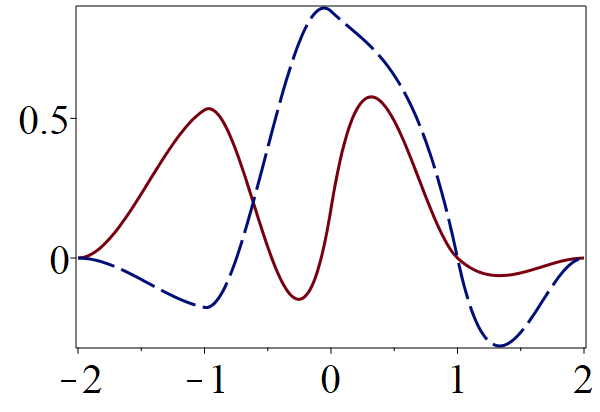}
			\caption{$\mrphi$}
			\label{fig:PhiNew:HCS}
		\end{subfigure}
		 \begin{subfigure}[b]{0.24\textwidth}
			\centering
			 \includegraphics[width=\textwidth]{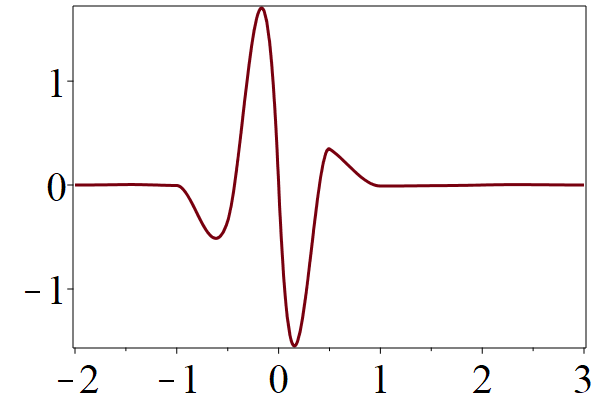}
			\caption{$\psi^1$}
			\label{fig:Psi1:HCS}
		\end{subfigure}
		 \begin{subfigure}[b]{0.24\textwidth}
			\centering
			 \includegraphics[width=\textwidth]{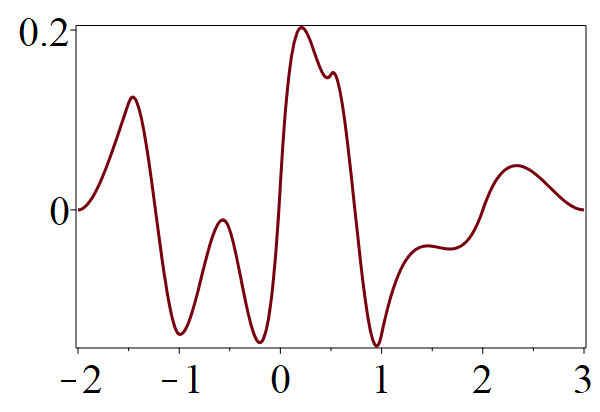}
			\caption{$\psi^2$}
			\label{fig:Psi2:HCS}
		\end{subfigure}\\
		 \begin{subfigure}[b]{0.24\textwidth}
			\centering
			 \includegraphics[width=\textwidth]{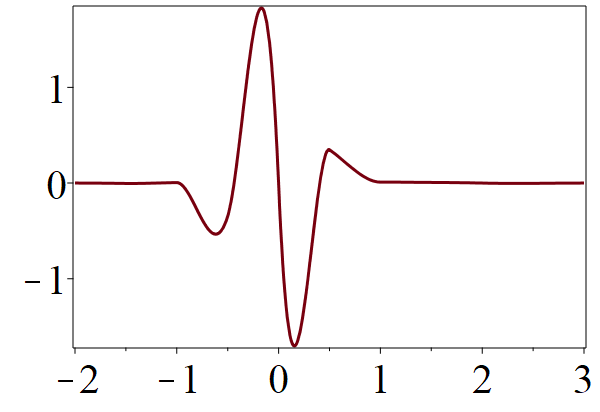}
			\caption{$\psi^3$}
			\label{fig:Psi3:HCS}
		\end{subfigure}
		 \begin{subfigure}[b]{0.24\textwidth}
			\centering
			 \includegraphics[width=\textwidth]{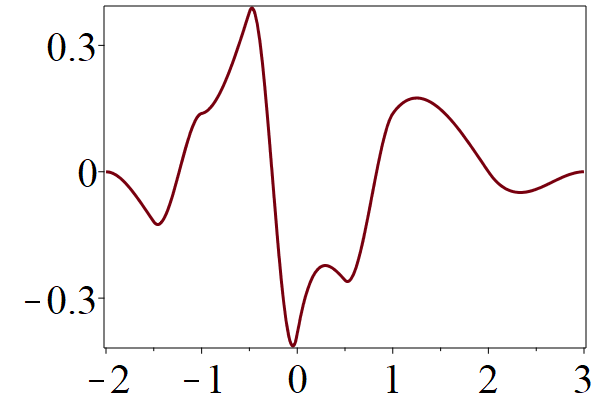}
			\caption{$\psi^4$}
			\label{fig:Psi4:HCS}
		\end{subfigure}
		 \begin{subfigure}[b]{0.24\textwidth}
			\centering
			 \includegraphics[width=\textwidth]{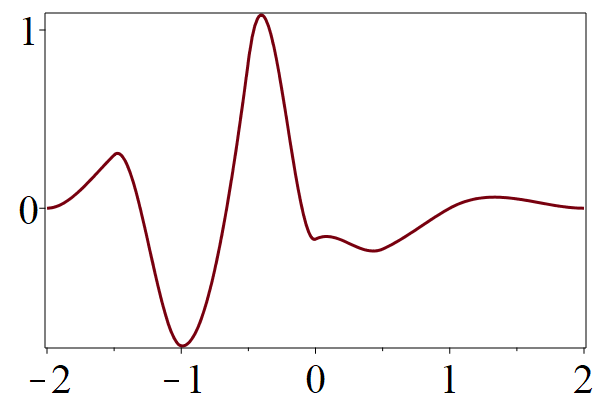}
			\caption{$\psi^5$}
			\label{fig:Psi5:HCS}
		\end{subfigure}
		 \begin{subfigure}[b]{0.24\textwidth}
			\centering
			 \includegraphics[width=\textwidth]{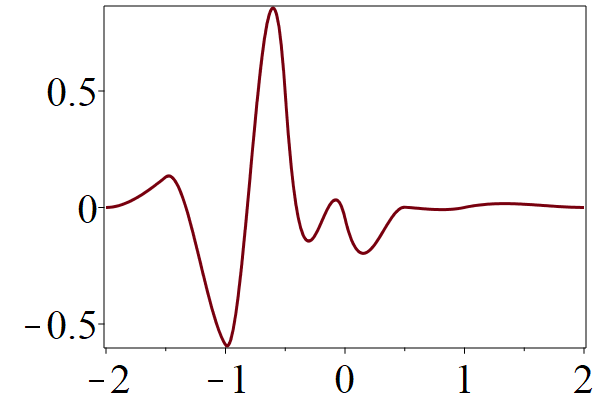}
			\caption{$\psi^6$}
			\label{fig:Psi6:HCS}
		\end{subfigure}
		\caption{Graphs of $\phi$ and the new refinable vector function $\mathring{\phi}$, together with graphs of $\psi^1,,\dots,\psi^6$ constructed from the Hermite cubic splines $\phi$ defined as \er{HCS} in Example~\ref{expl:hermite}.
$\{\mathring{\phi}; [\psi^1,\ldots,\psi^6]^\tp\}_{(-1,-1,1,1,1,1)}$ is a compactly supported quasi-tight $2$-framelet in $\Lp{2}$ with balanced vanishing moments $2$.
}\label{fig:psi:HCS}
	\end{figure}

\begin{exmp}\label{expl:hermite}
Let $\phi=[\phi_1,\phi_2]^{\tp}$ be the Hermite cubic splines as defined in \er{HCS}. We have $\wh{\phi}(2\cdot)=\wh{a}\wh{\phi}$, where $a\in\lrs{0}{2}{2}$ is given in \er{mask:HCS:0}. It is well known that $\sr(a,2)=4$ with a matching filter $\vgu\in\lrs{0}{1}{2}$ satisfying $\wh{\vgu}(\xi)=[1,i\xi]+\bo(|\xi|^4)$ as $\xi\to 0$.
Thus there exist quasi-tight $2$-framelets derived from $\phi$ which satisfy all claims of Theorem~\ref{thm:qtf}, with the maximum possible choice $m=4$ for balanced vanishing moments.
For simplicity of presentation, here we present an example of quasi-tight framelets with $m=2$ instead. Following the construction guidelines, we first construct
a desired strongly invertible filter $\theta\in \lrs{0}{2}{2}$ as follows:
$$
\wh{\theta}(\xi)=\frac{\sqrt{2}}{256}\begin{bmatrix}96e^{i\xi}+32 &-77e^{-i\xi}+506+51e^{i\xi}\\
	160e^{i\xi}-32 &-385e^{-i\xi}-78-17e^{i\xi}\end{bmatrix}.
$$
Direct computation shows that \er{vgu:bo} and \er{phi:bo} hold with $m=2$ and
	 $$\wh{c}(\xi)=\frac{\sqrt{2}}{2}-\frac{i\sqrt{2}}{4}\xi+\bo(|\xi|^2),\quad \wh{d}(\xi)=\frac{\sqrt{2}}{2}+\frac{i\sqrt{2}}{4}\xi+\bo(|\xi|^2),\quad\xi\to 0.$$

We obtain $b\in\lrs{0}{6}{2}$ such that $\{a;b\}_{\Theta;(\epsilon_1,\dots,\epsilon_6)}$ is a finitely supported quasi-tight $2$-multiframelet filter bank with $\wh{\Theta}(\xi)=\ol{\wh{\theta}(\xi)}^{\tp}\wh{\theta}(\xi)$, $\epsilon_1=\epsilon_2=-1$, and $\epsilon_3=\dots=\epsilon_6=1$.
For simplicity of presentation, we write
	 $$\wh{b}(\xi)=\begin{bmatrix}\DG(-1,1)-\tfrac{1}{4}\DG(-1,1)N(2\xi) &\pmb{0}_{2\times2}\\[0.3cm]
	 \DG(-1,1)+\tfrac{1}{4}\DG(-1,1)N(2\xi) &\pmb{0}_{2\times2}\\[0.3cm]
	\pmb{0}_{2\times2}& De^{2i\xi}
	 \end{bmatrix}\tilde{D}\wh{E}(2\xi)\wh{F}(\xi)\wh{\theta}(\xi),$$
	where $\pmb{0}_{q\times t}$ denotes the $q\times t$ zero matrix and
\begin{itemize}
\item $N(\xi)$ is the $2\times 2$ matrix of $2\pi$-periodic trigonometric polynomials given by
		 $$N(\xi)=\left(N_{-1}e^{i\xi}+N_0+N_1e^{-i\xi}\right)\DG\left(\frac{28831}{932734773870005846016},\frac{104831}{70129776762814464}\right),$$
		where
		 $$N_{-1}=\DG\left(\frac{2219987}{92160},\frac{8071987}{74502}\right)\begin{bmatrix}-280650717637&32961105478501\\
		 1406309548267&-205511772233035\end{bmatrix}\DG\left(47,\frac{5}{53759}\right),$$
		 $$N_{0}=\DG\left(\frac{28831}{230400},\frac{104831}{37251}\right)\begin{bmatrix}45834164001503531&-3709367687537217\\
		 -3709367687537217&2370098256094979\end{bmatrix}\DG\left(1,\frac{25}{53759}\right),$$
		 $$N_{1}=\DG\left(\frac{104339389}{18432},\frac{8071987}{74502}\right)\begin{bmatrix}-280650717637&1406309548267\\
		 164805527392505&-1027558861165175\end{bmatrix}\DG\left(\frac{1}{5},\frac{1}{53759}\right).$$
		
		\item $D=\DG(d_1,d_2)$ where
		 $$d_1=\frac{1761312\sqrt{13212226268199396309514273}}{51195503191172527},\quad d_2=\frac{\sqrt{547545488675642}}{37192694}.$$

		\item $\tilde{D}=\DG(d_3,d_4,d_5,d_6)$ where
		 $$d_3=\frac{1}{16174191},\quad d_4=\frac{1}{1505540889600},$$
		 $$d_5=\frac{1}{24064994385271402392453120},\quad d_6=\frac{1}{269129558593851555840}.$$

		\item $\wh{E}$ is the $4\times 4$ matrix of $2\pi$-periodic trigonometric polynomials given by
		 $$\wh{E}(\xi)=D_{-1}E_{-1}e^{i\xi}+D_0E_0+D_1E_1e^{-i\xi}+D_2E_2e^{-2i\xi}+E_3e^{-3i\xi},$$
		where
		 $$D_{-1}=\DG\left(187,1505540889600,263278430792912227,370729762969181\right),$$
		 $$D_0=\DG\left(1,1,187,18596347\right),$$
		 $$D_1=\DG\left(318703,318703,14399,18596347\right),$$
		 $$D_2=\DG\left(16174191,4589004497,1108723,1431918719\right),$$
		
		$$E_{-1}=\begin{bmatrix}
		0 & -3713 & -18565 & 5417\\
		0 & 0 & 0 & 0\\
		0 & 3713 & 18565 & -5417\\
		0 & 3713 & 18565 & -5417\\
		\end{bmatrix}$$
		 $$E_0=\footnotesize\begin{bmatrix}16174191 & 7619899 & 4457703 & 956109\\
		0 & -204477398651 & -262538115175 & -230405377741\\
		-121774017724980520053 & -210030256538769846979 & 192032829613853577201 &  427912364665047103611\\
		-1724291840139 & -4565648413949 & -9050015332017 &  -535702395579\end{bmatrix},$$
		 $$E_1=\small\begin{bmatrix}-1 & 0 & 0 & 0\\
		-505023 & -344267 & -201399 & -43197\\
		5110725327108891443 & 1450280624007721927 & -383860742942300317 &  -1583445964329957343\\
		-824326243735 & -1270122003803 & -1307289140551 &  -967745336173\end{bmatrix},$$
		 $$E_2=\small\begin{bmatrix}0 & 0 & 0 & 0\\
		1 & 0 & 0 & 0\\
		-19768729169137143 & -12286216971251881 & -7187537033154957 &  -1541616578141871\\
		-8001370677 & -5625792235 & -3291134295 & -705897885\end{bmatrix},$$
		$$E_3=\begin{bmatrix}0 & 0 & 0 & 0\\
		0 & 0 & 0 & 0\\
		569741919122396546911 & 0 & 0 &  0\\
		336929465078003105 & 0 & 0 & 0\end{bmatrix}.$$
		
		\item $\wh{F}$ is the $4\times 2$ matrix of $2\pi$-periodic trigonometric polynomials given by \er{Fbo}.
	\end{itemize}
	The filter $b$ is supported on $[-3,5]$. Define $\psi=[\psi^1,\psi^2,\psi^3,\psi^4,\psi^5,\psi^6]^{\tp}$ via $\wh{\psi}(\xi)=\wh{b}(\xi/2)\wh{\phi}(\xi/2)$.
Define a new refinable vector function $\wh{\mathring{\phi}}(\xi):=\wh{\theta}(\xi) \wh{\phi}(\xi)$.
Then $\|\wh{\mathring{\phi}}(\xi)\|^2=1+\bo(|\xi|^4)$ as $\xi\to 0$ and $\{\mrphi;\psi\}_{(\epsilon_1,\dots,\epsilon_6)}$ is a compactly supported quasi-tight $2$-framelet in $\Lp{2}$
such that all the desired properties in items (1)--(4) of Theorem~\ref{thm:qtf} are satisfied with $m=2$.
Note that $\psi$ has $2$ vanishing moments. See Figure~\ref{fig:psi:HCS} for graphs of $\phi,\mathring{\phi},\psi^1,\dots,\psi^6.$
\end{exmp}

\end{document}